\documentclass[10pt,journal,compsoc]{IEEEtran}
\usepackage{times}
\usepackage{epsfig}
\usepackage{graphicx}
\usepackage{amsmath}
\usepackage{amssymb}
\usepackage[english]{babel}
\usepackage[utf8x]{inputenc}
\usepackage[english]{babel}
\usepackage[OT1]{fontenc}
\usepackage{amssymb}
\usepackage{amsmath}
\usepackage{amsthm}
\usepackage{subfig}
\usepackage{epsfig}
\usepackage{epstopdf}
\usepackage{xcolor,colortbl}
\usepackage{graphicx}
\usepackage{algorithm}
\usepackage{algorithmic}
\usepackage{amsfonts}
\usepackage{fancyhdr}
\usepackage{dsfont}
\usepackage{color}
\usepackage{bm}
\usepackage{afterpage}
\usepackage{multirow}
\usepackage{diagbox}
\usepackage[font=small,skip=0pt]{caption}
\usepackage{epsfig}
\usepackage{epstopdf}
\usepackage{xcolor,colortbl}
\usepackage{graphicx}
\usepackage{algorithm}
\usepackage{algorithmic}
\usepackage{amsfonts}
\usepackage{fancyhdr}

\usepackage[breaklinks=true,bookmarks=false]{hyperref}

\definecolor{colortwo}{RGB}{255,113,31}
\definecolor{colorone}{RGB}{218,218,217}

\newcommand{\bbb}[1]{\boldsymbol{\mathbf{#1}}}
\def\cone{\textcolor[rgb]{0.9961,0,0}}
\def\ctwo{\textcolor[rgb]{0,0.5820,0.9}}
\def\cthree{\textcolor[rgb]{0,0.7,0}}

\newtheorem{lemma}{Lemma}
\newtheorem{theorem}{Theorem}

\def\beq{\begin{eqnarray}}
\def\eeq{\end{eqnarray}}
\def\noi{\noindent}
\def\nn{\nonumber}
\def\la{\langle}
\def\ra{\rangle}

\def\ghs{\hspace{0.05cm}} 

\def\L{\mathcal{L}}
\def\J{\mathcal{J}}
\newcommand{\figsizetwo}{\fontsize{4.8}{4.8}\selectfont}
\newcommand{\tablefont}{\fontsize{3.7}{13.8}\selectfont}

\newcommand{\figurewidth}{0.14\textwidth}
\newcommand{\figureheight}{0.13\textwidth}

\newenvironment{psmallmatrix}
  {\left(\begin{smallmatrix}}
  {\end{smallmatrix}\right)}
\ifCLASSOPTIONcompsoc
  \usepackage[nocompress]{cite}
\else
  \usepackage{cite}
\fi

\begin{document}

\title{$\ell_0$TV: A Sparse Optimization Method for Impulse Noise Image Restoration}

\author{Ganzhao~Yuan,~Bernard~Ghanem

\IEEEcompsocitemizethanks{\IEEEcompsocthanksitem Ganzhao~Yuan (corresponding author) is with School of Data and Computer Science, Sun Yat-sen University (SYSU), China, and also with Key Laboratory
of Machine Intelligence and Advanced Computing, Ministry of Education, China. E-mail: yuanganzhao@gmail.com.}

\IEEEcompsocitemizethanks{\IEEEcompsocthanksitem Bernard Ghanem is with Visual Computing Center, King Abdullah University of Science and Technology (KAUST), Saudi Arabia. E-mail: bernard.ghanem@kaust.edu.sa.}


}

\markboth{IEEE Transactions on Pattern Analysis and Machine Intelligence}%
{Shell \MakeLowercase{\textit{et al.}}: Bare Demo of IEEEtran.cls for Computer Society Journals}

\IEEEtitleabstractindextext{%
\begin{abstract}
Total Variation (TV) is an effective and popular prior model in the field of regularization-based image processing. This paper focuses on total variation for removing impulse noise in image restoration. This type of noise frequently arises in data acquisition and transmission due to many reasons, e.g. a faulty sensor or analog-to-digital converter errors. Removing this noise is an important task in image restoration. State-of-the-art methods such as Adaptive Outlier Pursuit(AOP) \cite{yan2013restoration}, which is based on TV with $\ell_{02}$-norm data fidelity, only give sub-optimal performance. In this paper, we propose a new sparse optimization method, called $\ell_0TV$-PADMM, which solves the TV-based restoration problem with $\ell_0$-norm data fidelity. To effectively deal with the resulting non-convex non-smooth optimization problem, we first reformulate it as an equivalent biconvex Mathematical Program with Equilibrium Constraints (MPEC), and then solve it using a proximal Alternating Direction Method of Multipliers (PADMM). Our $\ell_0TV$-PADMM method finds a desirable solution to the original $\ell_0$-norm optimization problem and is proven to be convergent under mild conditions. We apply $\ell_0TV$-PADMM to the problems of image denoising and deblurring in the presence of impulse noise. Our extensive experiments demonstrate that $\ell_0TV$-PADMM outperforms state-of-the-art image restoration methods.

\end{abstract}

}

\maketitle
\IEEEdisplaynontitleabstractindextext
\IEEEpeerreviewmaketitle

\renewcommand{\arraystretch}{1}
\IEEEraisesectionheading{\section{Introduction}\label{sec:introduction}}

Image restoration is an inverse problem, which aims at estimating the original \emph{clean} image $\bbb{u}$ from a blurry and/or noisy observation $\bbb{b}$. Mathematically, this problem is formulated as:
\beq
 \textstyle \bbb{b}  = \left((\bbb{Ku}) \odot { \bbb{\varepsilon}}_m \right) + \bbb{\varepsilon}_a, \label{eq:inverse}
\eeq

\begin{figure}[!t]
\begin{center}
\includegraphics[height=0.96in,width=3.1in]{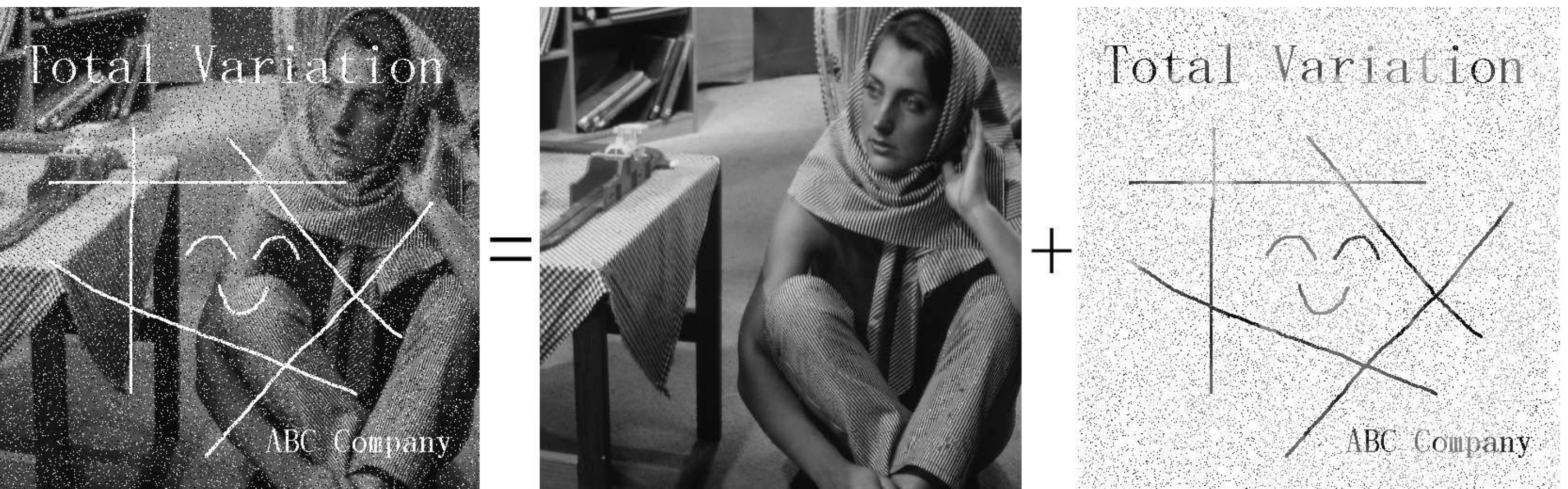}\vspace{0.05cm}
\includegraphics[height=0.96in,width=3.1in]{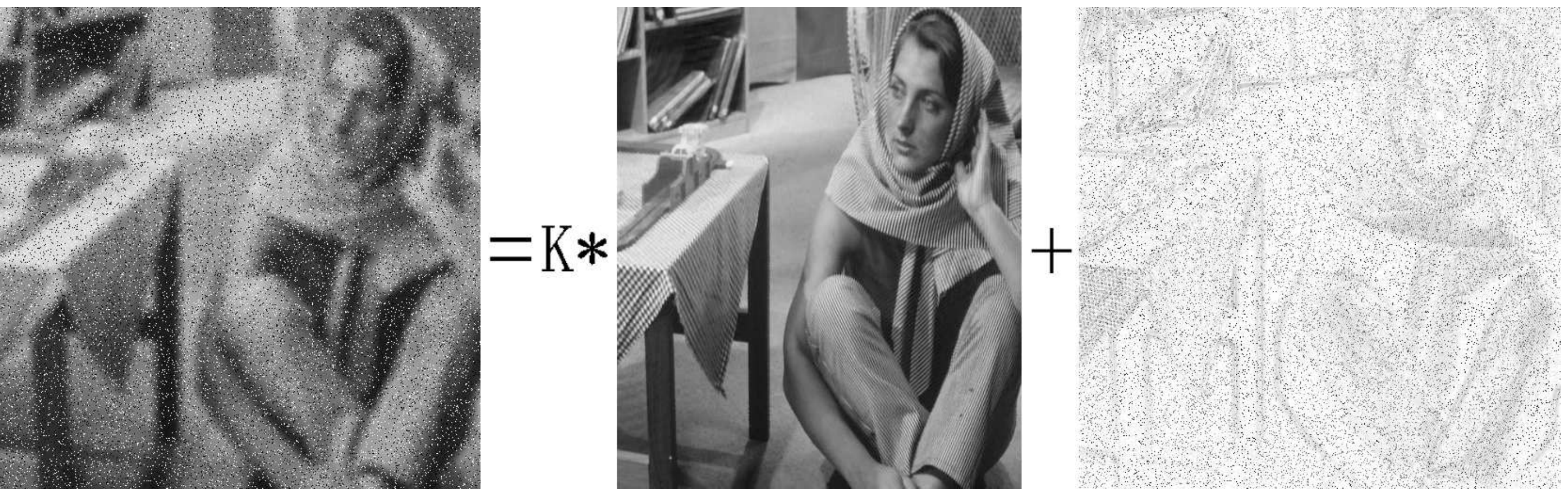}
\end{center}
\caption{An example of an image recovery result using our proposed $\ell_0$TV-PADMM method. Left column: corrupted image. Middle column: recovered image. Right column: absolute residual between these two images. 
}
\label{fig:intro:pic}
\end{figure}

\noi where $\bbb{K}$ is a linear operator, $\bbb{\varepsilon}_m$ and $\bbb{\varepsilon}_a$ are the noise vectors,
and $\odot$ denotes an elementwise product. Let $\bbb{1}$ and $\bbb{0}$ be column vectors of all entries equal to one and zero, respectively. When $\bbb{\varepsilon}_m=\bbb{1}$ and $\bbb{\varepsilon}_a \neq \bbb{0}$ (or $\bbb{\varepsilon}_m \neq \bbb{0}$ and $\bbb{\varepsilon}_a = \bbb{0}$), (\ref{eq:inverse}) corresponds to the additive (or multiplicative) noise model. For convenience, we adopt the vector representation for images, where a 2D $M\times N$ image is column-wise stacked into a vector $\bbb{u}\in\mathbb{R}^{n \times 1}$ with $n=M\times N$. So, for completeness, we have $\bbb{1}, \bbb{0},\bbb{b},\bbb{u},\bbb{\varepsilon}_a, \bbb{\varepsilon}_m \in \mathbb{R}^{n}$, and $\bbb{K}\in \mathbb{R}^{n\times n}$. Before proceeding, we present an image restoration example on the well-known `barbara' image using our proposed method for solving impulse noise removal in Figure \ref{fig:intro:pic}.

 In general image restoration problems, $\bbb{K}$ represents a certain linear operator, e.g. convolution,
 wavelet transform, etc., and recovering $\bbb{u}$ from $\bbb{b}$ is known as image deconvolution or image deblurring.
 When $\bbb{K}$ is the identity operator, estimating $\bbb{u}$ from $\bbb{b}$ is referred to as image denoising \cite{wang2008new}.
 The problem of estimating $\bbb{u}$ from $\bbb{b}$ is called a linear inverse problem which, for most scenarios of practical interest,
is ill-posed due to the singularity and/or the ill-conditioning of $\bbb{K}$. Therefore, in order to stabilize the recovery of $\bbb{u}$, it is necessary to incorporate prior-enforcing regularization on the solution. Therefore, image restoration can be modelled globally as the following  optimization problem:

\begin{align}
\textstyle \min_{ \bbb{u}  } ~\ell(\bbb{Ku},\bbb{b})+ \lambda~\Omega(\bbb{\nabla}_x \bbb{u},\bbb{\nabla}_y \bbb{u}), \label{eq:inverse:opt}
\end{align}

\noi where $\ell(\bbb{Ku},\bbb{b})$ measures the data fidelity between $\bbb{Ku}$ and the observation $\bbb{b}$, $\bbb{\nabla}_x\in \mathbb{R}^{n\times n}$ and $\bbb{\nabla}_y\in \mathbb{R}^{n\times n}$ are two suitable linear transformation matrices such that $\bbb{\nabla}_x \bbb{u} \in \mathbb{R}^{n}$ and $\bbb{\nabla}_y \bbb{u} \in \mathbb{R}^{n}$ compute the discrete gradients of the image $\bbb{u}$ along the $x$-axis and $y$-axis, respectively\footnote{In practice, one does not need to compute and store the matrices $\bbb{\nabla}_x$ and $\bbb{\nabla}_y$ explicitly. Since the adjoint of the gradient operator ${\bbb{\nabla}}$ is the negative divergence operator $-\bbb{\text{div}}$, i.e.,
 $\la \bbb{r}, \bbb{\nabla}_x \bbb{u} \ra = \la - \bbb{\text{div}}_x \bbb{r},\bbb{u} \ra, \la \bbb{s}, \bbb{\nabla}_y \bbb{u} \ra = \la - \bbb{\text{div}}_y \bbb{s},\bbb{u} \ra$ for any $\bbb{r}, \bbb{s}\in \mathbb{R}^{n}$,
 the inner product between vectors can be evaluated efficiently. Fore more details on the computation of $\bbb{\nabla}$
 and $\bbb{\text{div}}$ operators, please refer to \cite{chambolle2004algorithm,Weiss2006some,aujol2009some}.}, $\Omega(\bbb{\nabla}_x \bbb{u},\bbb{\nabla}_y \bbb{u})$ is the regularizer on $\bbb{\nabla}_x \bbb{u}$ and $\bbb{\nabla}_y \bbb{u}$, and $\lambda$ is a positive parameter used to balance the two terms for minimization. Apart from regularization, other prior information such as bound constraints \cite{beck2009fast,zuo2011generalized} or hard constraints can be incorporated into the general optimization framework in (\ref{eq:inverse:opt}).

{
\setlength\extrarowheight{4pt}
\begin{table}[!t]
\caption{Data Fidelity Models}\label{tab:dist}
\begin{center}
\scalebox{0.75}{\begin{tabular}{|c|c|}
  \hline
  Data Fidelity Function          & Noise and References \\
  \hline
   $\ell_{2}(\bbb{Ku},\bbb{b}) = \|\bbb{Ku}-\bbb{b}\|_2^2$                                 &add. Gaussian noise         \cite{rudin1992nonlinear,chambolle2004algorithm} \\
\hline
   $\ell_{1}(\bbb{Ku},\bbb{b}) = \|\bbb{Ku}-\bbb{b}\|_1$                                   &add. Laplace noise      \cite{yang2009efficient,Clason2010Duality} \\
\hline
   $\ell_{\infty}(\bbb{Ku},\bbb{b}) =\|\bbb{Ku}-\bbb{b}\|_{\infty}$                         &add. uniform noise    \cite{Clason2012fitting,Weiss2006some}\\
\hline
   $\ell_{p}(\bbb{Ku},\bbb{b}) = \la \bbb{Ku} - \bbb{b} \odot \log(\bbb{Ku})  , \bbb{1}\ra$     & mul. Poisson noise       \cite{Le2007Variational,steidl2010removing} \\
\hline
   $\ell_{g}(\bbb{Ku},\bbb{b})= \la \log(\bbb{Ku}) + \bbb{b} \odot \frac{1}{\bbb{Ku}} , \bbb{1} \ra $ & mul. Gamma noise  \cite{aubert2008variational,woo2013proximal}\\
\hline
   $\ell_{r}(\bbb{Ku},\bbb{b})= \la \log(\bbb{Ku}) + \bbb{b} \odot \bbb{b} \odot \frac{1}{\bbb{2Ku}} , \bbb{1} \ra $ & mul. Rayleigh noise  \cite{seabra2008convex,afonso2015blind}\\
\hline
   $\ell_{02}(\bbb{Ku},\bbb{b})= \|\bbb{Ku}-\bbb{b}+\bbb{z}\|_2^2,s.t.\|\bbb{z}\|_0 \leq k$ & mixed Gaussian impulse noise   \cite{yan2013restoration}\\
\hline
   $\ell_{0}(\bbb{Ku},\bbb{b})= \|\bbb{Ku}-\bbb{b}\|_0 $ & add./mul. impulse noise  \bbb{[ours]} \\
  \hline
\end{tabular}}
\end{center}
\end{table}
}

\subsection{Related Work}

 This subsection presents a brief review of  existing TV methods, from the viewpoint of data fidelity models, regularization models and optimization algorithms.

\vspace{3pt}\noindent\textbf{Data Fidelity Models:} The fidelity function $\ell(\cdot,\cdot)$ in (\ref{eq:inverse:opt})  usually penalizes the difference between $\bbb{Ku}$ and $\bbb{b}$ by using different norms/divergences. Its form
 depends on the assumed distribution of the noise model. Some typical noise models and
 their corresponding fidelity terms are listed in Table \ref{tab:dist}. The classical TV model
 \cite{rudin1992nonlinear} only considers TV minimization involving the squared $\ell_2$-norm fidelity term for recovering
 images corrupted by additive Gaussian noise. However, this model is far from optimal when the noise is not Gaussian. Other works \cite{yang2009efficient,Clason2010Duality} extend classical TV to use the $\ell_1$-norm in the fidelity term. Since the $\ell_1$-norm fidelity term coincides with the probability density function of Laplace distribution, it is suitable for image restoration in the presence of Laplace noise. Moreover, additive uniform noise \cite{Clason2012fitting,Weiss2006some}, multiplicative Poisson noise \cite{Le2007Variational}, and multiplicative Gamma noise
 \cite{woo2013proximal} have been considered in the literature. Some extensions have been made to deal with mixed Rayleigh impulse noise and mixed Poisson impulse noise in \cite{afonso2015blind}. Recently, a sparse noise model using an $\ell_{02}$-norm for data fidelity has been investigated in \cite{yan2013restoration} to remove impulse and mixed Gaussian impulse noise. In this paper, we consider $\ell_0$-norm data fidelity and show that it is particularly suitable for reconstructing images corrupted with additive/multiplicative \footnote{The impulse noise has a discrete nature (corrupted or uncorrupted), thus it can be viewed as additive noise or multiplicative noise.} impulse noise.

\vspace{3pt}\noindent \textbf{Regularization Models:} Several regularization models have been studied in the literature
 (see Table \ref{tab:reg}). The Tikhonov-like regularization \cite{tikhonov1977solution} function $\Omega_{\text{tik}}$ is quadratic and smooth, therefore it is relatively inexpensive to minimize with  first-order smooth optimization methods. However,
 since this method tends to overly smooth images, it often erodes strong edges and texture details. To address this issue, the total variation (TV) regularizer was proposed by Rudin, Osher and Fatemi in \cite{rudin1992nonlinear} for image denoising. Several other variants of TV have been extensively studied. The original TV norm $\Omega_{\text{tv}_{2}}$  in \cite{rudin1992nonlinear} is isotropic, while an anisotropic variation $\Omega_{\text{tv}_{1}}$ is also used.  From a numerical point of view, $\Omega_{\text{tv}_{2}}$ and $\Omega_{\text{tv}_{1}}$ cannot be directly minimized since they are not differentiable. A popular method is to use their smooth approximation $\Omega_{\text{stv}}$ and $\Omega_{\text{hub}}$ (see \cite{nikolova2005analysis} for details). Very recently, the Potts model $\Omega_{\text{pot}}$ \cite{geman1984,Mumford1989,BoykovVZ01}, which is based on the $\ell_0$-norm, has received much attention. It has been shown to be particularly effective for image smoothing \cite{xu2011image} and motion deblurring \cite{xu2013Unnatural}. 



{
\setlength\extrarowheight{4pt}
\begin{table}[!t]
\caption{Regularization Models}\label{tab:reg}
\begin{center}
\scalebox{0.745}{\begin{tabular}{|c|c|}
  \hline
  Regularization Function           & Description and References \\
  \hline
   $ \Omega_{\text{tik}}(\bbb{g},\bbb{h})=\sum_{i=1}^{n} \bbb{g}_i^2 +  \bbb{h}_i^2 $                                & Tikhonov-like       \cite{tikhonov1977solution} \\
  \hline
  $ \Omega_{\text{tv}_{2}}(\bbb{g},\bbb{h})= \sum_{i=1}^{n} {(\bbb{g}_i^2 + \bbb{h}_i^2)}^{\frac{1}{2}}$          & Isotropic      \cite{rudin1992nonlinear,woo2013proximal} \\
    \hline
     $ \Omega_{\text{tv}_{1}}(\bbb{g},\bbb{h})= \sum_{i=1}^{n} {|\bbb{g}_i| + |\bbb{h}_i|}$          & Anisotropic      \cite{wang2008new,yang2009efficient} \\
     \hline
     $ \Omega_{\text{stv}}(\bbb{g},\bbb{h})=\sum_{i=1}^{n}  {(\bbb{g}_i^2 + \bbb{h}_i^2 +\varepsilon^2)}^{\frac{1}{2}}$          & smooth TV     \cite{chan1999nonlinear,Weiss2006some} \\
     \hline
     $\Omega_{\text{pot}}(\bbb{g},\bbb{h})= \sum_{i=1}^{n} { |\bbb{g}_i|_0 + |\bbb{h}_i|_0 }$          & Potts model   \cite{xu2011image,xu2013Unnatural} \\
     \hline
     $ \begin{array}{c}
       \Omega_{\text{hub}}(\bbb{g},\bbb{h}) = \sum_{i=1}^{n} \varphi(\bbb{g}_i;\bbb{h}_i),~ \\\varphi(\bbb{g}_i;\bbb{h}_i) = {\tiny \begin{cases} \varepsilon \| \bbb{g}_i;\bbb{h}_i\|_2^2/2;~\|\bbb{g}_i;\bbb{h}_i\|_2 \leq {1}/{\varepsilon}\\\|\bbb{g}_i;\bbb{h}_i\|_2-\varepsilon/2;~\text{otherwise}\end{cases}}\end{array}$ & Huber-Like     \cite{nikolova2005analysis} \\
  \hline
\end{tabular}}
\end{center}
\end{table}
}

\vspace{3pt}\noindent\textbf{Optimization Algorithms:} The optimization problems involved in TV-based image restoration  are usually difficult due to the non-differentiability of the TV norm and the high dimensionality of the image data.  In the past several decades, a plethora of approaches have been proposed, which include PDE methods based on the Euler-Lagrange equation \cite{rudin1992nonlinear},
 the interior-point method \cite{chan1999nonlinear}, the semi-smooth Newton method \cite{ng2007semismooth},
 the second-order cone optimization method \cite{goldfarb2005second}, the splitting Bregman method \cite{GoldsteinO09,zhang2010bregmanized}, the fixed-point iterative method \cite{chen2012fast}, Nesterov's first-order optimal method \cite{Nesterov03,beck2009fast}, and alternating direction methods \cite{wang2008new,chen2011matrix,woo2013proximal}.
 Among these methods, some solve the TV problem in its primal form \cite{wang2008new}, while others consider its dual or primal-dual forms \cite{chan1999nonlinear,Clason2010Duality}. In this paper, we handle the TV problem with $\ell_0$-norm data fidelity using a primal-dual formulation, where the resulting equality constrained optimization is solved using proximal Alternating Direction Method of Multipliers (PADMM). It is worthwhile to note that the Penalty Decomposition Algorithm (PDA) in \cite{LuZ13} can also solve our problem, however, it lacks numerical stability. This motivates us to design a new $\ell_0$-norm optimization algorithm in this paper.

\subsection{Contributions and Organization}
 The main contributions of this paper are two-fold. \textbf{(1)} $\ell_0$-norm data fidelity is proposed to address the TV-based image restoration problem\footnote{We are also aware of Ref. \cite{chartrand2007quasi} where $\ell_0$-norm data fidelity is considered. However, their interpretation from the MAP viewpoint is not correct. }. Compared with existing models, our model is particularly suitable for image restoration in the presence of impulse noise. \textbf{(2)} To deal with the resulting NP-hard \footnote{The $\ell_0$ norm  problem is known to be NP-hard \cite{Natarajan1995}, since it is equivalent to NP-complete subset selection problems.} $\ell_0$ norm optimization, we propose a proximal ADMM to solve an equivalent MPEC form of the problem. A preliminary version of this paper appeared in \cite{yuan2015l0tv}.

The rest of the paper is organized as follows. Section \ref{sec:motivation} presents the motivation and formulation of the problem for impulse noise removal. Section \ref{sec:optimization} presents the equivalent MPEC problem and our proximal ADMM solution. Section \ref{sec:connect} discusses the connection between our method and prior work. Section \ref{sec:exp} provides extensive and comparative results in favor of our $\ell_0$TV method. Finally, Section \ref{sec:conc} concludes the paper.

\renewcommand{\arraystretch}{1}
\section{Motivation and Formulations}\label{sec:motivation}

\subsection{Motivation}
\label{sect:motivation}
 This work focuses on image restoration in the presence of impulse noise, which is very common in data acquisition and transmission due to faulty sensors or analog-to-digital converter errors, etc. Moreover, scratches in photos and video sequences can be also viewed as a special type of impulse noise. However, removing this kind of noise is not easy, since corrupted pixels are randomly distributed in the image and the intensities at corrupted pixels are usually indistinguishable from those of their neighbors.  There are two main types of impulse noise in the literature \cite{Clason2010Duality,ji2011robust}: random-valued and salt-and-pepper impulse noise. Let $[u_{\min}, u_{\max}]$ be the dynamic range of an image, where $u_{\min} = 0$ and $u_{\max} = 1$ in this paper. We also denote the original and corrupted intensity values at position $i$ as $\bbb{u}_i$ and $\mathcal{T}(\bbb{u}_i)$, respectively. 

 \noi \textbf{Random-valued impulse noise}: A certain percentage of pixels are altered to take on a uniform random number $d_i\in[u_{\min}, u_{\max}]$:
\beq \label{eq:def:rv}
\textstyle \mathcal{T}(\bbb{u}_i) = \begin{cases}
d_{i}, &\text{with probability}~r_{rv};\\
(\bbb{Ku})_{i}, &\text{with probability}~1-r_{rv}.\\
\end{cases}
\eeq
\noi \textbf{Salt-and-pepper impulse noise}: A certain percentage of pixels are altered to be either $u_{\min}$ or $u_{\max}$:
\beq \label{eq:def:sp}
\textstyle \mathcal{T}(\bbb{u}_i)= \begin{cases}
u_{\min}, &\text{with probability}~r_{sp}/2;\\
u_{\max}, &\text{with probability}~r_{sp}/2;\\
(\bbb{Ku})_{i}, &\text{with probability}~1-r_{sp}.\\
\end{cases}
\eeq



\noi  The above definition means that impulse noise corrupts a portion of pixels in the image while keeping other pixels unaffected. Expectation maximization could be used to find the MAP estimate of $\bbb{u}$ by maximizing the conditional posterior probability $p(\bbb{u}|\mathcal{T}(\bbb{u}))$, the probability that $\mathbf{u}$ occurs when $\mathcal{T}(\mathbf{u})$ is observed. By the Bayes' theorem, we have that
\beq 
\textstyle p(\bbb{u}|\mathcal{T}(\bbb{u}))={p(\bbb{u})\cdot p(\mathcal{T}(\bbb{u})|\bbb{u})}~/~{p(\mathcal{T}(\bbb{u}))}. \nn
\eeq
\noi Taking the negative logarithm of the above equation, the estimate is a solution of the following
minimization problem: 
\beq \label{eq:inter:log}
\textstyle \max_{\bbb{u}}~\log p(\mathcal{T}(\bbb{u})|\bbb{u}) + \log p(\bbb{u}).
\eeq
We now focus on the two terms in (\ref{eq:inter:log}). (i) The expression $p(\mathcal{T}(\bbb{u})|\bbb{u})$ can be viewed as a fidelity term measuring the discrepancy between the estimate $\bbb{u}$ and the noisy image $\mathcal{T}(\mathbf{u})$. The choice of the likelihood $p(\mathcal{T}(\mathbf{u})|\mathbf{u})$ depends upon the property of noise. From the definition of impulse noise given above, we have that
\beq \label{eq:puu}
\textstyle p(\mathcal{T}(\bbb{u})|\bbb{u}) &=& 1-r = 1 - \|\mathcal{T}(\bbb{u}) - \bbb{b}\|_0 /{n}, \nn
\eeq
\noi where $r$ is the noise density level as defined in (\ref{eq:def:rv}) and (\ref{eq:def:sp}) and $\|\cdot\|_0$ counts the number of non-zero elements in a vector. (ii) The term $p(\bbb{u})$ in (\ref{eq:inter:log}) is used to regularize a solution that has a low probability. We use a prior which has the Gibbs form: $p(\bbb{u}) = \frac{1}{\vartheta} \exp(-E(\bbb{u}))$ with $E(\bbb{u}) = \sigma \cdot \Omega_{\text{tv}}(\bbb{\nabla}_x \mathbf{u},\bbb{\nabla}_y \mathbf{u})$. Here, $E(\bbb{u})$ is the TV prior energy functional, $\vartheta$ is a normalization factor such that the TV prior is a probability, and $\sigma$ is the free parameter of the Gibbs measure. Replacing $p(\mathcal{T}(\mathbf{u})|\mathbf{u})$ and $p(\mathbf{u})$ into (\ref{eq:inter:log}) and ignoring a constant, we obtain the following $\ell_0TV$ model:
 \beq \label{eq:l0tv:1}
\textstyle \min_{\bbb{u}}~ \|\mathbf{Ku}-\mathbf{b}\|_0 +  \lambda \sum_{i=1}^{n} \Big[|(\bbb{\nabla}_x \mathbf{u})_i|^p + |(\bbb{\nabla}_y \mathbf{u})_i|^p\Big]^{1/p},\nn
 \eeq
\noi where $\lambda$ is a positive number related to $n$, $\sigma$ and $r$. The parameter $p$ can be 1 (anisotropic TV) or $2$ (isotropic TV), and $(\bbb{\nabla}_x \mathbf{u})_i$ and $(\bbb{\nabla}_y \mathbf{u})_i$ denote the $i$th component of the vectors $\bbb{\nabla}_x \bbb{u}$ and $\bbb{\nabla}_y \bbb{u}$, respectively. For convenience, we define $\forall \bbb{x}\in\mathbb{R}^{2n}$:
 \beq
  \textstyle \|\bbb{x}\|_{p,1}\triangleq\sum_{i=1}^n(|\bbb{x}_i|^p+|\bbb{x}_{n+i}|^p)^{\frac{1}{p}};~\bbb{\nabla} \triangleq \left[\bbb{\nabla}_x \atop \bbb{\nabla}_y\right]\in \mathbb{R}^{2n \times n}.   \nn
\eeq
\noi In order to make use of more prior information, we consider the following box-constrained model:
\beq \label{eq:l0tv:2}
\textstyle \min_{\bbb{0} \leq \bbb{u} \leq \bbb{1}}~ \|\bbb{o}\odot \left(\bbb{Ku}-\bbb{b}\right)\|_0 + \lambda \| \bbb{\nabla} \bbb{u}\|_{p,1},
\eeq
where $\bbb{o}\in\{0,1\}^n$ is specified by the user. When $\bbb{o}_i$ is 0, it indicates the pixel in position $i$ is an outlier, while when $\bbb{o}_i$ is 1, it indicates the pixel in position $i$ is a potential outlier. For example, in our experiments, we set $\bbb{o}=\mathbf{1}$ for the random-valued impulse noise and $ \mathbf{o}_i\scriptsize = \begin{cases}
0, &{\bbb{b}_i=u_{\min}~\text{or}~u_{\max}}\\
1, &{\text{otherwise}}\\
\end{cases}$ for the salt-and-pepper impulse noise. In what follows, we focus on optimizing the general formulation in (\ref{eq:l0tv:2}).  





\subsection{Equivalent MPEC Reformulations}

 In this section, we reformulate the problem in (\ref{eq:l0tv:2}) as an equivalent MPEC from a primal-dual viewpoint. First, we provide the variational characterization of the $\ell_0$-norm using the following lemma. 
 \begin{lemma}\label{le:mpec}
  For any given $\bbb{w}\in \mathbb{R}^{n}$, it holds that
  \beq\label{eq:mpec}
  \|\bbb{w}\|_0=\min_{\bbb{0}\le \bbb{v}\le \bbb{1}}~\la \bbb{1},\bbb{1}-\bbb{v}\ra,~s.t.~\bbb{v}\odot|\bbb{w}|=\bbb{0},
 \eeq
 and $\bbb{v}^*=\bbb{1}-{\rm sign}(|\bbb{w}|)$ is the unique optimal solution of the problem in (\ref{eq:mpec}). Here, the standard signum function sign is applied componentwise, and ${\rm sign}(0) = 0$.

%
\begin{proof}


The total number of zero elements in $\bbb{w}$ can be computed as $n - \|\bbb{w}\|_0  = \max_{\bbb{v}\in\{0,1\}}~\sum_{i=1}^n \bbb{v}_i,~s.t.~\bbb{v} \in \Phi$, where $\Phi \triangleq \{ \bbb{v}~|~\bbb{v}_i \cdot|\bbb{w}_i|=0,~\forall i \in [n]\}\nn$. Note that when $\bbb{w}_i=0$, $\bbb{v}_i=1$ will be achieved by maximization, when $\bbb{w}_i\neq0$, $\bbb{v}_i=0$ will be enforced by the constraint. Thus, $\bbb{v}^*_i=1-{\rm sign}(|\bbb{w}_i|)$. Since the objective function is linear, maximization is always achieved at the boundaries of the feasible solution space. Thus, the constraint of $\bbb{v}_i\in\{0,1\}$ can be relaxed to $0\leq \bbb{v}_i\leq1$, we have: $\|\bbb{w}\|_0 = n -\max_{ \mathbf{0} \leq \bbb{v} \leq \mathbf{1},~\bbb{v} \in \Phi}~\sum_{i=1}^n \bbb{v}_i = \min_{ \mathbf{0} \leq \bbb{v} \leq \mathbf{1},~\bbb{v} \in \Phi }~{\la \mathbf{1},\mathbf{1}-\bbb{v}\ra}$.

  \end{proof}
 \end{lemma}


The result of Lemma \ref{le:mpec} implies that the $\ell_0$-norm minimization problem in (\ref{eq:l0tv:2}) is equivalent to
\beq
\begin{split}
\label{eq:main:mpec}
  &\textstyle \min_{\bbb{0} \leq \bbb{u},\bbb{v} \leq \bbb{1}}~  \la \bbb{1}, \bbb{1}-\bbb{v}\ra + \lambda \| \bbb{\nabla} \bbb{u}\|_{p,1} \\
  &\textstyle \quad {\rm s.t.}\quad \bbb{v}\odot |\bbb{o}\odot (\bbb{Ku}-\bbb{b})| =\bbb{0}.
\end{split}
  \eeq
 If $\bbb{u}^*$ is a global optimal solution of (\ref{eq:l0tv:2}), then
 $(\bbb{u}^*,\bbb{1}-{\rm sign}(|\bbb{Ku}^*-\bbb{b}|))$ is globally optimal to (\ref{eq:main:mpec}). Conversely, if $(\bbb{u}^*,\bbb{1}-{\rm sign}(|\bbb{Ku}^*-\bbb{b}|))$ is a global optimal solution of (\ref{eq:main:mpec}),
 then $\bbb{u}^*$ is globally optimal to (\ref{eq:l0tv:2}).

 Although the MPEC problem in (\ref{eq:main:mpec}) is obtained by increasing the dimension of the original $\ell_0$-norm problem in (\ref{eq:l0tv:2}), this does not lead to additional local optimal solutions. Moreover, compared with (\ref{eq:l0tv:2}), (\ref{eq:main:mpec}) is a non-smooth non-convex minimization problem and its non-convexity is only caused by the complementarity
 constraint $\bbb{v}\odot |\bbb{o}\odot (\bbb{Ku}-\bbb{b})| =\bbb{0}$.

Such a variational characterization of the $\ell_0$-norm is proposed in \cite{dAspremont2003,Hu2008,feng2013complementarity,BiLP14,Bi2014}, but it is not used to develop any optimization algorithms for $\ell_0$-norm problems. We argue that, from a practical perspective, improved solutions to (\ref{eq:l0tv:2}) can be obtained by reformulating the $\ell_0$-norm in terms of complementarity constraints \cite{luo1996mathematical,yuan2015l0tv,yuan2016proximal,yuan2016binary,yuan2016sparsity,yuan2017exact}. In the following section, we will develop an algorithm to solve (\ref{eq:main:mpec}) based on proximal ADMM and show that such a ``lifting'' technique can achieve a desirable solution of the original $\ell_0$-norm optimization problem.

\section{Proposed Optimization Algorithm} \label{sec:optimization}

This section is devoted to the solution of (\ref{eq:main:mpec}). This problem is rather difficult to solve, because it is neither convex nor smooth. Our solution is based on the proximal ADM method, which iteratively updates the primal and dual variables of the augmented Lagrangian function of (\ref{eq:main:mpec}).

First, we introduce two auxiliary vectors $\bbb{x}\in\mathbb{R}^{2n}$ and
 $\bbb{y}\in\mathbb{R}^n$ to reformulate (\ref{eq:main:mpec}) as:
 \beq \label{eq:opt:sub}
  & \min_{\bbb{0} \leq \bbb{u},\bbb{v} \leq \bbb{1},~\bbb{x},~\bbb{y} }~\la \bbb{1}, \bbb{1}-\bbb{v}\ra + \lambda \|\bbb{x}\|_{p,1} \\
  & \quad\ \ {\rm s.t.} \quad \bbb{\nabla} \bbb{u} = \bbb{x},~\bbb{Ku} - \bbb{b} = \bbb{y},~\bbb{v}\odot \bbb{o}\odot |\bbb{y}| =\bbb{0}. \nn
 \eeq
 Let $\mathcal{L}:\mathbb{R}^n\times\mathbb{R}^n\times\mathbb{R}^{2n}\times\mathbb{R}^{n}\times\mathbb{R}^{2n}\times\mathbb{R}^{n}\times\mathbb{R}^{n}\to\mathbb{R}$
 be the augmented Lagrangian function of (\ref{eq:opt:sub}).
\beq
\textstyle \mathcal{L}(\bbb{u},\bbb{v},\bbb{x},\bbb{y},\bbb{\xi},\bbb{\zeta},\bbb{\pi}) := \la \bbb{1}, \bbb{1}-\bbb{v}\ra + \lambda\|\bbb{x}\|_{p,1}+\nn~~~~~~~~~~~~\\
\textstyle \la \bbb{\nabla} \bbb{u} - \bbb{x}, \bbb{\xi}\ra  +\frac{\beta}{2} \|\bbb{\nabla} \bbb{u} - \bbb{x}\|^2  + \la \bbb{Ku}-\bbb{b}-\bbb{y},\bbb{\zeta} \ra+\nn~~~~~~~~~~~\\
\textstyle \frac{\beta}{2} \|\bbb{Ku}-\bbb{b}-\bbb{y}\|^2 +\la \bbb{v}\odot \bbb{o}\odot |\bbb{y}|,\bbb{\pi} \ra + \frac{\beta}{2} \|\bbb{v}\odot \bbb{o}\odot |\bbb{y}|\|^2,\nn~~~~~
\eeq
\noi where $\bbb{\xi}$, $\bbb{\zeta}$  and $\bbb{\pi}$ are the Lagrange multipliers associated with the constraints
 ${\bbb{\nabla u}} = \bbb{x}$, $\bbb{Ku} - \bbb{b} = \bbb{y}$ and $\bbb{v}\odot \bbb{o}\odot |\bbb{y}|=0$, respectively, and $\beta>0$ is the penalty parameter.
 The detailed iteration steps of the proximal ADM for (\ref{eq:opt:sub}) are described in Algorithm \ref{alg:proximal:admm}. In simple terms, ADM updates are performed by optimizing for a set of primal variables at a time, while keeping all other primal and dual variables fixed. The dual variables are updated by gradient ascent on the resulting dual problem. 

\begin{algorithm} [!t]
\caption{\label{alg:proximal:admm} {\bf ($\ell_0TV$-ADMM) A Proximal ADMM for Solving the Biconvex MPEC Problem (\ref{eq:main:mpec})}}
\textbf{(S.0)}~Choose a starting point $(\bbb{u}^{0},\bbb{v}^{0},\bbb{x}^{0},\bbb{y}^{0},\bbb{\xi}^{0},\bbb{\zeta}^{0})$. Set $k=0$. Select step size $\gamma \in (0,2)$, $\mu>0$, $\beta=1$, and $L={\mu+\beta\|\bbb{\nabla}\|^2+\beta\|\bbb{K}\|^2}$.

\vspace{3pt}\textbf{(S.1)}~Solve the following minimization problems with $\bbb{D}:=L\bbb{I} - (\beta\bbb{\nabla}^T\bbb{\nabla}+\beta \bbb{K}^T\bbb{K})$ and $\bbb{E}:=\mu\bbb{I}$:
     \beq
\begin{split}
\left[\begin{array}{c}
           \bbb{u}^{k+1}\\
           \bbb{v}^{k+1}\\
        \end{array}\right]= \mathop{\arg\min}_{\bbb{0}\le \bbb{u},\bbb{v}\le \bbb{1}}\mathcal{L}(\bbb{u},\bbb{v},\bbb{x}^k,\bbb{y}^k,\bbb{\xi}^k,\bbb{\zeta}^k,\bbb{\pi}^k)\\
+\tfrac{1}{2}\|\bbb{u}-\bbb{u}^{k}\|^2_{\bbb{D}} + \tfrac{1}{2}\|\bbb{v}-\bbb{v}^{k}\|^2_{\bbb{E}}  ~~~\label{subprob2}
\end{split}~~~~~~
\eeq
\beq
\left[\begin{array}{c}
           \bbb{x}^{k+1}\\
           \bbb{y}^{k+1}\\
        \end{array}\right]=\mathop{\arg\min}_{\bbb{x},\bbb{y}} \mathcal{L}(\bbb{u}^{k+1},\bbb{v}^{k+1},\bbb{x},\bbb{y},\bbb{\xi}^k,\bbb{\zeta}^k,\bbb{\pi}^k)\label{subprob3}
 \eeq

\vspace{3pt}\textbf{(S.2)}~Update the Lagrange multipliers:
                \beq
                &&\bbb{\xi}^{k+1}=\bbb{\xi}^{k} + \gamma \beta(\bbb{\nabla} \bbb{u}^{k} - \bbb{x}^{k}),\label{update-multi1}\\
                &&\bbb{\zeta}^{k+1}=\bbb{\zeta}^{k} + \gamma \beta(\bbb{K} \bbb{u}^{k} -\bbb{b}- \bbb{y}^{k}),\label{update-multi2}\\
                &&\bbb{\pi}^{k+1}=\bbb{\pi}^{k} + \gamma \beta(\bbb{o}\odot \bbb{v}^{k}\odot  |\bbb{y}^{k}|).\label{update-multi3}
                \eeq
\vspace{3pt}\textbf{(S.3)}~if $\left(k \text{~is a multiple of~} 30\right)$, then $\beta=\beta \times\sqrt{10}$
\vspace{3pt}\textbf{(S.4)}~ Set $k:=k+1$ and then go to Step (S.1).
 \end{algorithm}

Next, we focus our attention on the solutions of the subproblems in (\ref{subprob2}) and (\ref{subprob3}) arising in Algorithm \ref{alg:proximal:admm}. We will show that the computation required in each iteration of Algorithm \ref{alg:proximal:admm} is insignificant. 





\noi (i) $(\bbb{u},\bbb{v})$-subproblem. Proximal ADM introduces a convex proximal term to the objective. The specific form of $\bbb{D}$ is chosen to expedite the computation of the closed form solution. The introduction of $\mu$ is to guarantee strongly convexity of the subproblems.

 $\bbb{u}$-subproblem in (\ref{subprob2}) reduces to the following minimization problem:

 \beq\label{eq:subprob:u}
 \begin{split}
 \textstyle  \bbb{u}^{k+1}  = \mathop{\arg \min}_{\bbb{0}\leq \bbb{u} \leq \bbb{1}}\ \tfrac{\beta}{2}\|\bbb{\nabla} \bbb{u} - \bbb{x}^{k}+\bbb{\xi}^k/\beta\|^2+  \\
 \textstyle\tfrac{\beta}{2} \|\bbb{Ku}-\bbb{b} - \bbb{y}^{k}+\bbb{\zeta}^{k}/\beta\|^2+\tfrac{1}{2}\|\bbb{u}-\bbb{u}^{k}\|_{\bbb{D}}^2.
 \end{split}
 \eeq
  After an elementary calculation, subproblem (\ref{eq:subprob:u}) can be simplified as
  \begin{align}\label{uk}\textstyle
    \bbb{u}^{k+1} =\mathop{\arg \min}_{\bbb{0}\leq \bbb{u} \leq \bbb{1}}\tfrac{1}{2}\|\bbb{u}-(\bbb{u}^{k}-\bbb{g}^{k}/L)\|^2\nn
  \end{align}
  with $\bbb{g}^{k}=\bbb{\nabla}^T\xi^{k}+\bbb{K}^T\bbb{\zeta}^{k}+\beta\bbb{\nabla}^T(\bbb{x}^{k} - \bbb{\nabla} \bbb{u}^{k})+\beta \bbb{K}^T(\bbb{b}+\bbb{y}^{k}-\bbb{K}\bbb{u}^{k})$. Then, the solution $\bbb{u}^{k}$ of (\ref{subprob2}) has the following closed form expression:
  \beq\textstyle
   \bbb{u}^{k+1} = \min (\bbb{1},\max (\bbb{0},\bbb{u}^{k}-\bbb{g}^{k}/L )). \nn
  \eeq
\noi Here the parameter $L$ depends on the spectral norm of the linear matrices $\bbb{\nabla}$ and $\bbb{K}$. Using the definition of $\bbb{\nabla}$ and the classical finite differences that $\| \bbb{\nabla}_y\| \leq 2$ and $\| \bbb{\nabla}_y\| \leq 2$ (see \cite{aujol2009some,chambolle2004algorithm,zuo2011generalized}), the spectral norm of $\bbb{\nabla}$ can be computed by: $\|\bbb{\nabla}\| = \|  \begin{psmallmatrix}
                          \bbb{\nabla}_x \\
                          \bbb{0} \\
                        \end{psmallmatrix}
+  \begin{psmallmatrix}
                          \bbb{0} \\
                          \bbb{\nabla}_y \\
                        \end{psmallmatrix} \|
                        \leq  \|  \begin{psmallmatrix}
                          \bbb{\nabla}_x \\
                          \bbb{0} \\
                        \end{psmallmatrix} \|
+ \|\begin{psmallmatrix}
                          \bbb{0} \\
                          \bbb{\nabla}_y \\
                        \end{psmallmatrix} \|
                        = \| \bbb{\nabla}_x\| + \| \bbb{\nabla}_y\| \leq 4$.

$\bbb{v}$-subproblem in (\ref{subprob2}) reduces to the following minimization problem:
 \beq\label{eq:subprob:v} \textstyle
  \bbb{v}^{k+1} =\mathop{\arg\min}_{\bbb{0}\leq \bbb{v} \leq \bbb{1}}~ \textstyle \tfrac{1}{2} \sum_{i=1}^n \bbb{s}_i^{k} \bbb{v}_i^2 +\la \bbb{v}, \bbb{c}^{k} \ra, \nn
 \eeq
 \noi where $\bbb{c}^{k}= o\odot \bbb{\pi}^{k}\odot  |\bbb{y}^{k}| - \bbb{1}- \mu \bbb{v}^k$, $\bbb{s}^{k}=\beta \bbb{o}\odot \bbb{y}^{k} \odot  \bbb{y}^{k} +\mu$. Therefore, the solution $\bbb{v}^{k}$ can be computed as:
 \beq \textstyle
\bbb{v}^{k+1} = \min(\bbb{1},\max(\bbb{0},-\tfrac{\bbb{c}^{k} }{\bbb{s}^{k}})). \nn
 \eeq
 \noi (iii) $(\bbb{x},\bbb{y})$-subproblem. Variable $\bbb{x}$ in (\ref{subprob3}) is updated by solving the following problem:
 \beq\label{eq:subprob:x}
 \textstyle \bbb{x}^{k+1} =\mathop{\arg\min}_{\bbb{x}\in\mathbb{R}^{2n}}\ \tfrac{\beta}{2}\|\bbb{x} -\bbb{h}^{k}\|^2 + \lambda \|\bbb{x}\|_{p,1}, \nn
 \eeq
  where $\bbb{h}^{k}:=\bbb{\nabla} \bbb{u}^{k+1}+\bbb{\xi}^{k}/\beta$. It is not difficult to check that for $p=1$,
 \beq\label{eq:update:x1}
    \bbb{x}^{k+1} ={\rm sign}\big(\bbb{h}^{k}\big)\odot \max\big(|\bbb{h}^{k}|-\lambda/\beta,0\big), \nn
 \eeq
  and when $p=2$,
 \beq\label{eq:update:x2}
  \textstyle  \left[\begin{array}{c}
           \bbb{x}_i^{k+1}\\
           \bbb{x}_{i+n}^{k+1}\\
        \end{array}\right]
     =\big(\max(0, 1-\frac{\lambda/\beta}{\|(\bbb{h}_i^{k};\bbb{h}_{i+n}^{k})\|})\big)
    \left[\begin{array}{c}
           \bbb{h}_i^{k}\\
           \bbb{h}_{i+n}^{k}\\
        \end{array}\right]\nn
 \eeq
 Variable $\bbb{y}$ in (\ref{subprob3}) is updated by solving the following problem:
  \beq \label{eq:subprob:z}
  \textstyle  \bbb{y}^{k+1} =\arg \min_{\bbb{y}}~ \tfrac{\beta}{2}\|\bbb{y} - \bbb{q}^{k}  \|^2 + \tfrac{\beta}{2} \| \bbb{w}^{k}\odot |\bbb{y}| + \bbb{\pi}^k/\beta\|^2, \nn
  \eeq
\noi where $\bbb{q}^{k}= \bbb{Ku}^{k+1} - \bbb{b} + \bbb{\zeta}^{k}/\beta  $ and $\bbb{w}^{k} = \bbb{o} \odot \bbb{v}^{k+1}$. A simple computation yields that the solution $\bbb{y}^{k}$ can be computed in closed form as:
 \beq \label{eq:update:y}
 \textstyle  \bbb{y}^{k+1}  = \text{sign}( \bbb{q}^{k} )\odot \max\big(0,\frac{ |\bbb{q}^{k}| - \bbb{\pi}^{k} \odot \bbb{w}^{k}/\beta   }{1 + \bbb{v}^{k} \odot \bbb{w}^{k}}\big),\nn
 \eeq


Proximal ADM has excellent convergence in practice. The global convergence of ADM for convex problems was given by He and Yuan in \cite{HeY12,chen2011matrix} under the variation inequality framework. However, since our optimization problem in (\ref{eq:main:mpec}) is non-convex, the convergence analysis for ADM needs additional conditions. By imposing some mild conditions, Wen et al. \cite{wen2012alternating}  managed to show that the sequence generated by ADM converges to a KKT point. Along a similar line, we establish the convergence property of proximal ADM. Specifically, we have the following convergence result.



\begin {theorem}
\label{LowRankDPConvergence}
\textbf{Convergence of Algorithm \ref{alg:proximal:admm}.} Let $X\triangleq(\bbb{u,v,x,y})$, $Y\triangleq(\bbb{\xi,\zeta,\pi})$ and $\{X^{k},Y^{k}\}_{k=1} ^{\infty}$ be the sequence generated by Algorithm \ref{alg:proximal:admm}. Assume that $\{Y^{k}\}_{k=1} ^{\infty}$ is bounded and satisfies $\sum_{k=0}^{\infty}\|Y^{k+1}-Y^k\|_F^2 < \infty$. Then any accumulation point of sequence satisfies the KKT conditions of (\ref{eq:opt:sub}).
\begin{proof}
Please refer to \textbf{Appendix A}.
\end{proof}
\end{theorem}

\textbf{Remark 1.} The condition $\sum_{k=0}^{\infty}\|Y^{k+1}-Y^k\|_F^2 < \infty$ holds when the multiplier does not change in two consecutive iterations. By the boundedness of the penalty parameter $\beta$ and Eqs (\ref{update-multi1}-\ref{update-multi3}), this condition also indicates that the equality constraints in (\ref{eq:opt:sub}) are satisfied. This assumption can be checked by measuring the violation of the equality constraints. Theorem \ref{LowRankDPConvergence} indicates that when the equality constraint holds, PADMM converges to a KKT point. Though not satisfactory, it provides some assurance on the convergence of Algorithm \ref{alg:proximal:admm}.

\textbf{Remark 2.} Two reasons explain the good performance of our method. (i) It targets a solution to the original problem in (\ref{eq:l0tv:2}). (ii) It has monotone and self-penalized properties owing to the complimentarity constraints brought on by the MPEC. Our method directly handles the complimentary constraints in (\ref{eq:opt:sub})$: \bbb{v} \odot \bbb{o} \odot |\bbb{y}| = \bbb{0}$ with $\bbb{v}\geq \bbb{0}$. These constraints are the only sources of non-convexity for the optimization problem and they characterize the optimality of the KKT solution of (\ref{eq:l0tv:2}). These special properties of MPEC distinguish it from general nonlinear optimization \cite{yuan2016proximal,yuan2016sparsity,yuan2016binary,yuan2017exact}. We penalize the complimentary error of $\bbb{v} \odot \bbb{o} \odot |\bbb{y}|$ (which is always non-negative) and ensure that the error is decreasing in every iteration.


\section{Connection with Existing Work}\label{sec:connect}
In this section, we discuss the connection between the proposed method $\ell_0TV$-PADM and prior work.

\subsection{Sparse Plus Low-Rank Matrix Decomposition}
 Sparse plus low-rank matrix decomposition \cite{wright2009robust,ji2011robust} is becoming a powerful tool that effectively corrects large errors
 in structured data in the last decade. It aims at decomposing a given corrupted image $\bbb{B}$ (which is of matrix form)
 into its sparse component ($\bbb{S}$) and low-rank component ($\bbb{L}$) by solving: $\min_{\bbb{B},\bbb{L}}~\|\bbb{S}\|_0+\lambda ~rank(\bbb{L}),~s.t.~\bbb{B}=\bbb{L}+\bbb{S}$. Here the sparse component represents the foreground of an image which can be treated as outliers or impulse noise,
 while the low-rank component corresponds to the background, which is highly correlated. This is equivalent to the following optimization problem:
 \beq
 \textstyle \min_{\bbb{L}} ~\|\bbb{B}-\bbb{L}\|_0+\lambda ~rank(\bbb{L}),\nn
 \eeq
 \noi which is also based on $\ell_0$-norm data fidelity. While they consider the low-rank prior in their objective function, we consider the Total Variation (TV) prior in ours.

\subsection{Convex Optimization Method $\ell_1TV$}
The goal of image restoration in the presence of impulse noise has been pursued by a number of authors (see, e.g., \cite{yang2009efficient,Clason2010Duality}) using $\ell_1TV$, which can be formulated as follows:
\beq \label{eq:l1}
\textstyle \min_{\bbb{0} \leq \bbb{u} \leq \bbb{1}} ~\|\bbb{K}\bbb{u}-\bbb{b}\|_1 + \lambda ~ \|\bbb{\nabla} \bbb{u}\|_{p,1}.
\eeq
It is generally believed that $\ell_1TV$ is able to remove the impulse noise properly. This is because $\ell_1$-norm provides the tightest convex relaxation for the $\ell_0$-norm over the unit ball in the sense of $\ell_\infty$-norm. It is shown in \cite{CandesT05} that the problem of minimizing $\|\bbb{Ku}-\bbb{b}\|_1$ is equivalent to $\|\bbb{Ku}-\bbb{b}\|_0$ with high probability under the assumptions that (i) $\bbb{Ku}-\bbb{b}$ is sparse at the optimal solution $\bbb{u}^*$ and (ii) $\bbb{K}$ is a random Gaussian matrix and sufficiently ``incoherent'' (i.e., number of rows in $\bbb{K}$ is greater than its number of columns). However, these two assumptions required in \cite{CandesT05} do not necessarily hold true for our $\ell_0TV$ optimization problem. Specifically, when the noise level of the impulse noise is high, $\bbb{Ku}-\bbb{b}$ may not be sparse at the optimal solution $\bbb{u}^*$. Moreover, the matrix $\bbb{K}$ is a square identity or ill-conditioned matrix. Generally, $\ell_1TV$ will only lead to a sub-optimal solution.


\subsection{Adaptive Outlier Pursuit Algorithm}
Very recently, Yan \cite{yan2013restoration} proposed the following new model for image restoration in the presence of impulse noise and mixed Gaussian impulse noise:
\beq \label{eq:l02tv}
\textstyle \min_{\bbb{u},\bbb{z}}~\chi\|\bbb{Ku}-\bbb{b}-\bbb{z}\|_2^2 + \|\bbb{\nabla} \bbb{u} \|_{p,1}, ~s.t.~\|\bbb{z}\|_0\leq k,~
\eeq
\noi where $\chi>0$ is the regularization parameter. They further reformulate the problem above into $\min_{\bbb{u},~\bbb{v}}~\|\bbb{v}\odot(\bbb{K}\bbb{u}-\bbb{b})\|_2^2 +   \lambda~\|\bbb{\nabla} \bbb{u} \|_{p,1},~s.t.~\bbb{0}\leq \bbb{v} \leq \bbb{1},~\la \bbb{v},\bbb{1} \ra \leq n-k$ and then solve this problem using an Adaptive Outlier Pursuit(AOP) algorithm. The AOP algorithm is actually an alternating minimization method, which separates the minimization problem over $\bbb{u}$ and $\bbb{v}$ into two steps. By iteratively restoring the images and updating the set of damaged pixels, it is shown that AOP algorithm outperforms existing state-of-the-art methods for impulse noise denoising, by a large margin.

Despite the merits of the AOP algorithm, we must point out that it incurs three drawbacks, which are unappealing in practice. First, the formulation in (\ref{eq:l02tv}) is only suitable for mixed Gaussian impulse noise, i.e. it produces a sub-optimal solution when the observed image is corrupted by pure impulse noise. (ii) Secondly, AOP is a multiple-stage algorithm. Since the minimization sub-problem over $\bbb{u}$\footnote{It actually reduces to the $\ell_2TV$ optimization problem.} needs to be solved exactly in each stage, the algorithm may suffer from slow convergence. (iii) As a by-product of (i), AOP inevitably introduces an additional parameter (that specifies the Gaussian noise level), which is not necessarily readily available in practical impulse denoising problems.

In contrast, our proposed $\ell_0$TV method is free from these problems. Specifically, (i) as have been analyzed in Section 2, i.e. our $\ell_0$-norm model is optimal for impulse noise removal. Thus, our method is expected to produce higher quality image restorations, as seen in our results. (ii) Secondly, we have integrated $\ell_0$-norm minimization into a unified proximal ADM optimization framework, it is thus expected to be faster than the multiple stage approach of AOP. (iii) Lastly, while the optimization problem in (\ref{eq:l02tv}) contains two parameters, our model only contains one single parameter.

\subsection{Other $\ell_0$-Norm Optimization Techniques}

Actually, the optimization technique for the $\ell_0$-norm regularization problem is the key to removing impulse noise. However, existing solutions are not appealing. The $\ell_0$-norm problem can be reformulated as a $0\text{-}1$ mixed integer programming \cite{bienstock1996computational}problem which can be solved by a tailored branch-and-bound algorithm but it involves high computational complexity. The simple projection methods are inapplicable to our model since they assume the objective function is smooth. Similar to the $\ell_1$ relaxation, the convex methods such as $k$-support norm relaxation \cite{mcdonald2014spectral}, $k$-largest norm relaxation \cite{yu2011adversarial}, QCQP and SDP relaxations \cite{Chan2007} only provide loose approximation of the original problem. The non-convex methods such as Schatten $\ell_p$ norm \cite{GeJY11,LuTYL16}, re-weighted $\ell_1$ norm \cite{candes2008enhancing}, $\ell_{\text{1-2}}$ norm DC (difference of convex) approximation \cite{YinLHX15}, the Smoothly Clipped Absolute Deviation (SCAD) penalty method\cite{zhang2010nearly}, the Minimax Concave Plus (MCP) penalty method \cite{fan2001variable} only produce sub-optimal results since they give approximate solutions for the $\ell_0TV$ problem or incur high computational overhead.

We take $\ell_p$ norm approximation method for example and it may suffer two issues. First, it involves an additional hyper-parameter $p$ which may not be appealing in practice. Second, the $\ell_p$ regularized norm problem for general $p$ could be difficult to solve. This includes the iterative re-weighted least square method \cite{Lu2014} and proximal point method. The former approximates $\|\bbb{x}\|_p^p$ by $\sum_{i=1}^n(\bbb{x}_i^2+\epsilon)^{p/2}$ with a small parameter $\epsilon$ and solves the resulting re-weighted least squares sub-problem which reduces to a weighted $\ell_2TV$ problem. The latter needs to evaluate a relatively expensive proximal operator $\Pi(\bbb{a}) = \min_{\bbb{x}}~\frac{1}{2}\|\bbb{x}-\bbb{a}\|_2^2 + \lambda \|\bbb{x}\|_p^p$ in general, except that it has a closed form solution for some special values such as $p=\frac{1}{2}$ and $p=\frac{2}{3}$ \cite{xu2012l}.

Recently, Lu et al. propose a Penalty Decomposition Algorithm (PDA) for solving the $\ell_0$-norm optimization algorithm \cite{LuZ13}. As has been remarked in \cite{LuZ13}, direct ADM on the $\ell_0$ norm problem can also be used for solving $\ell_0TV$ minimization simply by replacing the quadratic penalty functions in the PDA by augmented Lagrangian functions. Nevertheless, as observed in our preliminary experiments and theirs, the practical performance of direct ADM is worse than that of PDA.

Actually, in our experiments, we found PDA is unstable. The penalty function can reach very large values $(\geq 10^8)$, and the solution can be degenerate when the minimization problem of the augmented Lagrangian function in each iteration is not exactly solved. This motivates us to design a new $\ell_0$-norm optimization algorithm in this paper. We consider a proximal ADM algorithm to the MPEC formulation of $\ell_0$-norm since it has a primal-dual interpretation. Extensive experiments have demonstrated that proximal ADM for solving the ``lifting'' MPEC formulation for $\ell_0TV$ produces better image restoration qualities.

\begin{figure*}[!ht]
\begin{center}
\includegraphics[width=0.9\textwidth,height=1.5in]{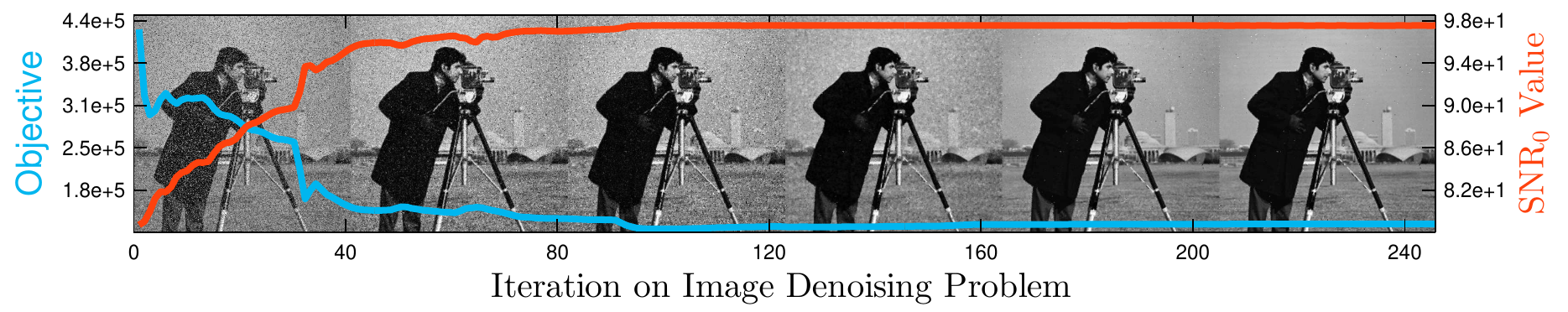}
\includegraphics[width=0.9\textwidth,height=1.5in]{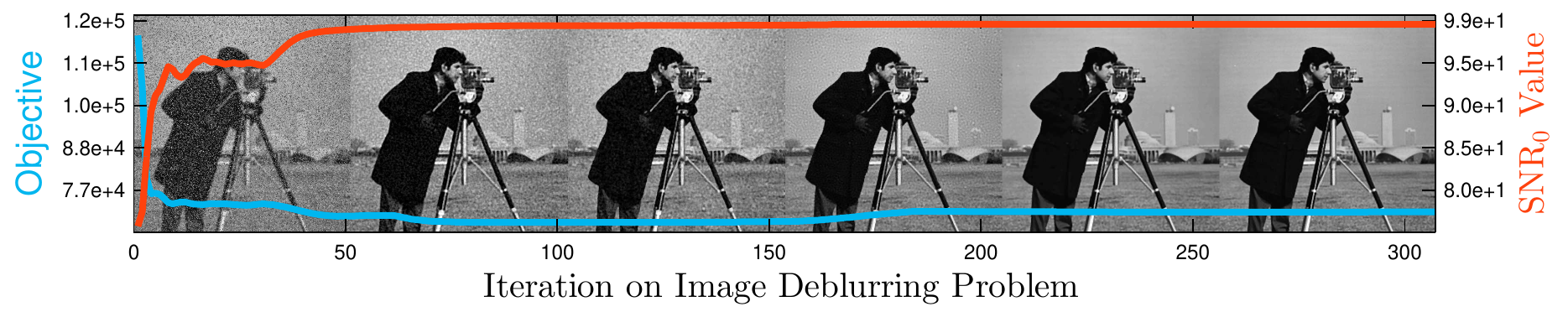}
\end{center}
\caption{Asymptotic behavior for optimizing (\ref{eq:l0tv:2}) to denoise and deblur the corrupted 'cameraman' image. We plot the value of the objective function (solid blue line) and the SNR value (dashed red line) against the number of optimization iterations. At specific iterations (i.e. 1, 10, 20, 40, 80, and 160), we also show the denoised and deblurred image. Clearly, the corrupting noise is being effectively removed throughout the optimization process.}
\label{fig:exp:convergence}
\end{figure*}
\section{Experimental Validation} \label{sec:exp}
 In this section, we provide empirical validation for our proposed  $\ell_0TV$-PADMM method by conducting extensive image denoising experiments and performing a thorough comparative analysis with the state-of-the-art. 


 In our experiments, we use 5 well-known test images of size $512\times 512$. All code is implemented in MATLAB using a 3.20GHz CPU and 8GB RAM. Since past studies \cite{cai2012image,chen2012fast} have shown that the isotropic TV model performs better than the anisotropic one, we choose $p=2$ as the order of the TV norm here. In our experiments, we apply the following algorithms:

\vspace{6pt}\noi \textbf{(i) BM3D} is an image denoising strategy based on an enhanced sparse representation in transform-domain. The enhancement of the sparsity is achieved by grouping similar 2D image blocks into 3D data arrays \cite{dabov2007image}.

\vspace{6pt}\noi \textbf{(ii) MFM}, Median Filter Methods. We utilize adaptive median filtering to remove salt-and-pepper impulse noise and adaptive center-weighted median filtering to remove random-valued impulse noise.

\vspace{6pt}\noi \textbf{(iii) $\ell_1TV$-SBM}, the Split Bregman Method (SBM) of \cite{GoldsteinO09}, which has been implemented in \cite{Getreuer10}. We use this convex optimization method as our baseline implementation.

\vspace{6pt}\noi \textbf{(iv) TSM}, the Two Stage Method\cite{ChanHN05,ChanHN04,Cai2010}. The method first detects the damaged pixels by MFM and then solves the TV image inpainting problem.

\vspace{6pt}\noi \textbf{(v) $\ell_pTV$-ADMM (direct)}. We directly use ADMM (Alternating Direction Method of Multipliers) to solve the non-smooth non-convex $\ell_p$ problem with proximal operator being computed analytically. We only consider $p=\frac{1}{2}$ in our experiments \cite{xu2012l}.

\vspace{6pt}\noi \textbf{(vi) $\ell_{02}TV$-AOP}, the Adaptive Outlier Pursuit (AOP) method described in \cite{yan2013restoration}. We use the implementation provided by the author. Here, we note that AOP iteratively calls the $\ell_1TV$-$SBM$ procedure, mentioned above.

\vspace{6pt}\noi \textbf{(vii) $\ell_0TV$-PDA}, the Penalty Decomposition Algorithm (PDA) \cite{LuZ13} for solving the $\ell_0TV$ optimization problem in (\ref{eq:l0tv:2}). 

\vspace{6pt}\noi \textbf{(viii) $\ell_0TV$-PADMM}, the proximal ADMM described in Algorithm \ref{alg:proximal:admm} for solving the $\ell_0TV$ optimization problem in (\ref{eq:l0tv:2}). We set the relaxation parameter to 1.618 and the strongly convex parameter $\mu$ to $0.01$. All MATLAB codes to reproduce the experiments of this paper are available online at the authors' research webpages.

%




\subsection{Experiment Setup}\label{sec4.1}
For the denoising and deblurring test, we use the following strategies to generate artificial noisy images.

\noi \textbf{(a) Denoising problem}. We corrupt the original image by injecting random-value, salt-and-pepper noise, and mixed noise (half random-value and half salt-and-pepper)  with different densities (10\% to 90\%) to the images. 

\noi \textbf{(b) Deblurring problem}. Although blurring kernel estimation has been pursued by many studies (e.g. \cite{xu2010two}), here we assume that the blurring kernel is known beforehand. We blur the original images with a $9 \times 9$ Gaussian blurring kernel and add impulse noise with different densities (10\% to 90\%). We use the following MATLAB scripts to generate a blurring kernel of radius $r$ ($r$ is set to 7 in the experiments):
\beq \label{eq:varyingr}
\begin{split}
\text{[x,y] = meshgrid~($-$r:r,~$-$r:r)},~~~~~~~~~~~~~ \\
\text{K=double(x.}~\hat{}~\text{2}+\text{y.}~\hat{}~2<=\text{r.}~\hat{}~\text{2)},~\text{P=K/sum(K(:))}.
\end{split}
\eeq 


We run all the previously mentioned algorithms on the generated noisy and blurry images. For $\ell_{02}TV$-AOP, we adapt the author's image denoising implementation to the image deblurring setting. Since both BM3D and Median Filter Methods (MFM) are not convenient to solve the deblurring problems, we do not test them in the deblurring problem. We terminate $\ell_0TV$-\text{PADMM} whenever $\|\bbb{\nabla} \bbb{u}^{k} - \bbb{x}^{k}\|_{2} \le \frac{1}{255}$ and $\|\bbb{K}\bbb{u}^{k}-\bbb{b}-\bbb{y}^{k}\|_2\le\frac{1}{255}$ and $\|\bbb{o}\odot \bbb{v}^{k}\odot |\bbb{y}^{k}|\|_2 \le \frac{1}{255}$. For $\ell_pTV$-\text{PADMM}, $\ell_0TV$-\text{PDA}, and $\ell_0TV$-\text{PADMM}, we use the same stopping criterion to terminate the optimization. For $\ell_1TV$-\text{SBM} and $\ell_{02}TV$-\text{AOP}, we adopt the default stopping conditions provided by the authors. For the regularization parameter $\lambda$, we swept over $\{0.1,0.6,1.1,...,9.6\}$. For the regularization parameter $\chi$ in $\ell_{02}TV$-\text{AOP}, we swept over $\{10,50,100,500,1000,5000,10000,50000\}$ and set $k$ to the number of corrupted pixels.

To evaluate these methods, we compute their Signal-to-Noise Ratios (SNRs). Since the corrupted pixels follow a Bernoulli-like distribution, it is generally hard to measure the data fidelity between the original images and the recovered images. Therefore, we consider three ways to measure SNR.
{ \beq \label{eq:snr}
\textstyle  SNR_{0}(\bm{u}) &\triangleq& \frac{n-\|\bm{u}^0 - \bm{u}\|_{0\text{-}\epsilon}}{n-\|\bm{u}^0 - \bm{u}^0\|_{0\text{-}\epsilon}} \times 100,\nn\\
\textstyle  SNR_{1}(\bm{u}) &\triangleq& 10 \log_{10} \frac{\|\bm{u}^0 - \bar{\bm{u}}\|_1}{\|\bm{u} - \bar{\bm{u}}\|_1},\nn\\
\textstyle SNR_{2}(\bm{u}) &\triangleq& 10 \log_{10} \frac{\|\bm{u}^0 - \bar{\bm{u}}\|_2^2}{\|\bm{u} - \bar{\bm{u}}\|_2^2},\nn
\eeq}{\normalsize}\noi where $\bm{u}^0$ is the original clean image and $\bar{\bm{u}}$ is the mean intensity value of $\bm{u}^0$, and $\|\cdot\|_{0\text{-}\epsilon}$ is the soft $\ell_0$-norm which counts the number of elements whose magnitude is greater than a threshold $\epsilon$. We adopt $\epsilon=\frac{20}{255}$ in our experiments.



\begin{table*}[!t]
\tiny
\begin{center}
\caption{General denoising problems. The results separated by `/' are $SNR_{0}$, $SNR_{1}$ and $SNR_{2}$, respectively. The $1^{st}$, $2^{nd}$, and $3^{rd}$ best results are colored with \textbf{\cone{red}}, \ctwo{blue} and \cthree{green}, respectively. }
\label{denoising:rv}
\scalebox{0.99}{\begin{tabular}{|p{1.7cm}|p{1.34cm}|p{1.34cm}|p{1.34cm}|p{1.34cm}|p{1.34cm}|p{1.34cm}|p{1.42cm}|p{1.42cm}|}
\hline
\diagbox{Img.}{Alg.} &   BM3D &$\ell_1TV$-$SBM$ &  $MFM$   & $TSM$&$\ell_{02}TV$-$AOP$& { \tablefont  $\ell_PTV$-$PADMM$ }& $\ell_0TV$-$PDA$    & { \tablefont $\ell_0TV$-$PADM$}                \\
\hline
\hline
\multicolumn{9}{|c|}{Random-Value Impulse Noise} \\
\hline
walkbridge+10\%  &  93/7.1/11.0  &  \cthree{95}/12.3/15.6  &  92/7.7/12.3  &  \cthree{95}/11.8/12.9  &  \ctwo{96}/\cthree{12.8}/\ctwo{16.6}  &  \cthree{95}/12.1/13.8  &  \textbf{\cone{97}}/\textbf{\cone{14.1}}/\textbf{\cone{16.9}}  &  \textbf{\cone{97}}/\ctwo{13.8}/\cthree{15.9}  \\
walkbridge+30\%  &  76/3.7/7.1  &  \ctwo{89}/\cthree{8.6}/11.0  &  82/6.1/10.3  &  \cthree{85}/5.8/7.8  &  \ctwo{89}/8.4/\ctwo{12.1}  &  \ctwo{89}/7.8/11.5  &  \textbf{\cone{91}}/\textbf{\cone{9.6}}/\textbf{\cone{12.8}}  &  \textbf{\cone{91}}/\ctwo{9.5}/\cthree{11.9}  \\
walkbridge+50\%  &  59/2.2/4.3  &  76/\cthree{4.9}/5.7  &  67/4.1/7.0  &  69/2.7/4.8  &  76/\ctwo{5.4}/8.1  &  \cthree{79}/\ctwo{5.4}/\cthree{8.7}  &  \ctwo{84}/\textbf{\cone{7.0}}/\textbf{\cone{10.1}}  &  \textbf{\cone{85}}/\textbf{\cone{7.0}}/\ctwo{9.2}  \\
walkbridge+70\%  &  42/1.0/1.9  &  56/2.0/1.7  &  45/2.0/3.3  &  50/1.3/2.2  &  53/2.5/4.0  &  \cthree{59}/\cthree{3.0}/\cthree{5.0}  &  \ctwo{65}/\ctwo{4.0}/\ctwo{6.2}  &  \textbf{\cone{76}}/\textbf{\cone{5.1}}/\textbf{\cone{7.0}}  \\
walkbridge+90\%  &  26/-0.1/-0.1  &  \cthree{32}/-0.2/-1.1  &  28/0.3/0.5  &  30/0.0/-0.0  &  31/\cthree{0.4}/\cthree{0.8}  &  30/\cthree{0.4}/\cthree{0.8}  &  \ctwo{34}/\ctwo{0.7}/\ctwo{1.3}  &  \textbf{\cone{57}}/\textbf{\cone{2.7}}/\textbf{\cone{3.9}}  \\
\hline
pepper+10\%  &  67/5.0/9.9  &  \ctwo{99}/\ctwo{19.1}/\cthree{21.5}  &  \ctwo{99}/15.0/\ctwo{22.2}  &  \cthree{97}/13.5/15.8  &  74/5.4/11.3  &  \ctwo{99}/13.6/20.3  &  \textbf{\cone{100}}/\textbf{\cone{20.2}}/\textbf{\cone{24.6}}  &  \ctwo{99}/\cthree{18.0}/21.0  \\
pepper+30\%  &  55/3.7/7.0  &  \ctwo{96}/\cthree{12.3}/13.6  &  \ctwo{96}/11.4/16.3  &  \cthree{87}/6.3/9.5  &  72/5.2/10.7  &  \textbf{\cone{98}}/12.0/\cthree{16.8}  &  \textbf{\cone{98}}/\textbf{\cone{15.1}}/\textbf{\cone{19.7}}  &  \textbf{\cone{98}}/\ctwo{14.6}/\ctwo{18.3}  \\
pepper+50\%  &  44/2.4/4.5  &  \cthree{85}/6.7/6.7  &  \cthree{85}/7.0/9.7  &  71/3.5/5.5  &  65/4.5/8.9  &  \ctwo{94}/\cthree{9.7}/\cthree{13.1}  &  \textbf{\cone{96}}/\textbf{\cone{11.8}}/\textbf{\cone{15.7}}  &  \textbf{\cone{96}}/\ctwo{11.6}/\ctwo{14.4}  \\
pepper+70\%  &  33/1.2/2.1  &  63/2.8/2.1  &  59/3.1/4.4  &  52/1.6/2.4  &  51/2.7/4.7  &  \cthree{79}/\cthree{5.2}/\cthree{6.2}  &  \ctwo{84}/\ctwo{6.8}/\ctwo{8.9}  &  \textbf{\cone{93}}/\textbf{\cone{9.0}}/\textbf{\cone{11.4}}  \\
pepper+90\%  &  24/0.2/0.1  &  \cthree{35}/0.1/-1.0  &  30/0.6/0.6  &  31/0.3/0.1  &  28/0.7/\cthree{1.1}  &  \cthree{35}/\cthree{0.9}/1.0  &  \ctwo{39}/\ctwo{1.3}/\ctwo{1.7}  &  \textbf{\cone{76}}/\textbf{\cone{4.2}}/\textbf{\cone{4.8}}  \\
\hline
mandrill+10\%  &  74/3.3/6.0  &  89/8.1/9.0  &  \cthree{92}/6.9/6.9  &  \ctwo{93}/\cthree{9.6}/\cthree{9.6}  &  84/3.7/7.4  &  \ctwo{93}/\cthree{9.6}/\cthree{9.6}  &  \textbf{\cone{95}}/\textbf{\cone{11.1}}/\textbf{\cone{11.5}}  &  \textbf{\cone{95}}/\ctwo{10.8}/\ctwo{10.3}  \\
mandrill+30\%  &  63/2.0/3.6  &  83/\cthree{5.9}/\cthree{6.6}  &  76/3.8/5.9  &  83/4.7/4.9  &  73/3.0/5.5  &  \cthree{85}/5.8/\ctwo{6.8}  &  \textbf{\cone{87}}/\textbf{\cone{6.8}}/\textbf{\cone{7.4}}  &  \ctwo{86}/\ctwo{6.4}/6.5  \\
mandrill+50\%  &  50/1.1/2.2  &  73/\cthree{3.6}/3.7  &  65/2.9/\cthree{4.6}  &  69/2.0/3.4  &  61/2.2/4.0  &  \cthree{74}/\cthree{3.6}/\ctwo{5.0}  &  \ctwo{77}/\textbf{\cone{4.6}}/\textbf{\cone{5.6}}  &  \textbf{\cone{78}}/\ctwo{4.4}/\cthree{4.6}  \\
mandrill+70\%  &  36/0.4/0.8  &  57/1.4/0.6  &  51/1.5/2.4  &  52/0.9/1.5  &  47/1.2/2.2  &  \cthree{62}/\cthree{2.3}/\cthree{3.4}  &  \ctwo{64}/\ctwo{2.9}/\textbf{\cone{3.9}}  &  \textbf{\cone{70}}/\textbf{\cone{3.1}}/\ctwo{3.5}  \\
mandrill+90\%  &  28/-0.3/-0.6  &  36/-0.6/-1.9  &  37/0.2/0.4  &  34/-0.1/-0.4  &  33/0.1/0.3  &  \cthree{39}/\cthree{0.5}/\cthree{0.9}  &  \ctwo{42}/\ctwo{0.8}/\ctwo{1.2}  &  \textbf{\cone{58}}/\textbf{\cone{1.9}}/\textbf{\cone{2.5}}  \\
\hline
lake+10\%  &  92/6.9/12.5  &  \textbf{\cone{98}}/\ctwo{16.9}/\textbf{\cone{21.3}}  &  \cthree{96}/11.3/17.7  &  \ctwo{97}/14.0/15.0  &  \ctwo{97}/8.7/16.1  &  \textbf{\cone{98}}/14.3/19.2  &  \textbf{\cone{98}}/\textbf{\cone{17.2}}/\ctwo{21.1}  &  \textbf{\cone{98}}/\cthree{16.7}/\cthree{19.5}  \\
lake+30\%  &  75/4.3/8.1  &  \ctwo{93}/\cthree{11.3}/13.9  &  91/9.3/\cthree{14.4}  &  86/7.1/10.0  &  \cthree{92}/7.9/13.9  &  \textbf{\cone{95}}/10.5/\ctwo{15.0}  &  \textbf{\cone{95}}/\textbf{\cone{12.7}}/\textbf{\cone{16.7}}  &  \textbf{\cone{95}}/\ctwo{12.0}/14.3  \\
lake+50\%  &  58/2.6/4.9  &  79/6.5/7.2  &  71/5.9/9.4  &  69/3.7/5.9  &  78/6.2/10.2  &  \cthree{88}/\cthree{8.3}/\ctwo{11.7}  &  \textbf{\cone{91}}/\textbf{\cone{10.0}}/\textbf{\cone{13.7}}  &  \ctwo{90}/\ctwo{9.5}/\cthree{11.5}  \\
lake+70\%  &  41/1.3/2.3  &  54/2.9/2.6  &  42/2.5/4.1  &  47/1.8/2.8  &  43/2.8/4.6  &  \cthree{60}/\cthree{4.7}/\cthree{7.0}  &  \ctwo{68}/\ctwo{5.8}/\ctwo{8.6}  &  \textbf{\cone{84}}/\textbf{\cone{7.4}}/\textbf{\cone{9.0}}  \\
lake+90\%  &  24/0.3/0.3  &  \ctwo{26}/0.5/-0.4  &  \cthree{25}/0.6/0.8  &  \ctwo{26}/0.5/0.4  &  24/0.6/1.0  &  13/\cthree{0.7}/\cthree{1.1}  &  \ctwo{26}/\ctwo{1.1}/\ctwo{1.7}  &  \textbf{\cone{62}}/\textbf{\cone{4.2}}/\textbf{\cone{5.3}}  \\
\hline
jetplane+10\%  &  \cthree{39}/2.5/6.1  &  \textbf{\cone{99}}/\textbf{\cone{17.5}}/\textbf{\cone{21.0}}  &  \ctwo{98}/11.5/17.5  &  \ctwo{98}/12.8/13.3  &  \cthree{39}/3.4/8.3  &  \textbf{\cone{99}}/13.1/\cthree{19.1}  &  \textbf{\cone{99}}/\ctwo{17.0}/\ctwo{20.0}  &  \ctwo{98}/\cthree{15.6}/17.0  \\
jetplane+30\%  &  32/0.7/2.6  &  \ctwo{95}/10.3/11.5  &  \cthree{94}/9.0/\cthree{13.3}  &  87/5.0/7.3  &  38/3.2/7.5  &  \textbf{\cone{97}}/\cthree{10.4}/\ctwo{15.0}  &  \textbf{\cone{97}}/\textbf{\cone{12.4}}/\textbf{\cone{15.7}}  &  \textbf{\cone{97}}/\ctwo{11.5}/12.6  \\
jetplane+50\%  &  27/-0.6/-0.1  &  \cthree{80}/4.5/4.0  &  75/4.2/6.7  &  69/1.5/2.8  &  34/2.4/5.2  &  \ctwo{92}/\cthree{7.9}/\ctwo{10.6}  &  \textbf{\cone{94}}/\textbf{\cone{9.3}}/\textbf{\cone{12.2}}  &  \textbf{\cone{94}}/\ctwo{9.0}/\cthree{10.0}  \\
jetplane+70\%  &  22/-1.7/-2.4  &  53/0.6/-0.7  &  42/0.2/0.9  &  47/-0.5/-0.5  &  23/-0.6/-0.3  &  \cthree{67}/\cthree{3.2}/\cthree{4.8}  &  \ctwo{74}/\ctwo{4.4}/\ctwo{6.4}  &  \textbf{\cone{90}}/\textbf{\cone{6.7}}/\textbf{\cone{7.4}}  \\
jetplane+90\%  &  18/-2.5/-4.1  &  \cthree{25}/-1.8/-3.6  &  \cthree{25}/-1.7/-2.5  &  \ctwo{26}/-1.8/-2.9  &  18/-2.3/-3.4  &  14/\cthree{-1.6}/\cthree{-2.2}  &  \ctwo{26}/\ctwo{-1.2}/\ctwo{-1.5}  &  \textbf{\cone{74}}/\textbf{\cone{3.4}}/\textbf{\cone{3.7}}  \\
\hline\hline
\multicolumn{9}{|c|}{Salt-and-Pepper Impulse Noise} \\
\hline
walkbridge+10\%  &  90/5.4/9.9  &  \cthree{96}/12.9/17.3  &  90/7.6/12.4  &  \ctwo{98}/15.8/19.9  &  \ctwo{98}/\cthree{16.3}/\cthree{20.7}  &  \ctwo{98}/15.8/19.9  &  \textbf{\cone{99}}/\ctwo{17.2}/\ctwo{22.7}  &  \textbf{\cone{99}}/\textbf{\cone{17.5}}/\textbf{\cone{23.2}}  \\
walkbridge+30\%  &  71/3.0/4.5  &  \cthree{94}/10.4/14.3  &  83/6.3/9.8  &  \ctwo{96}/\cthree{11.7}/\cthree{16.4}  &  \cthree{94}/10.5/15.2  &  \ctwo{96}/\cthree{11.7}/\cthree{16.4}  &  \ctwo{96}/\ctwo{12.0}/\ctwo{17.1}  &  \textbf{\cone{97}}/\textbf{\cone{12.3}}/\textbf{\cone{17.5}}  \\
walkbridge+50\%  &  51/-0.1/-1.7  &  \cthree{89}/8.1/11.4  &  71/4.0/5.4  &  \ctwo{92}/\ctwo{9.3}/\ctwo{14.0}  &  88/7.8/11.8  &  \ctwo{92}/\ctwo{9.3}/\cthree{13.9}  &  \ctwo{92}/\cthree{9.2}/13.8  &  \textbf{\cone{93}}/\textbf{\cone{9.5}}/\textbf{\cone{14.3}}  \\
walkbridge+70\%  &  32/-2.0/-4.6  &  \cthree{82}/6.1/8.7  &  49/1.4/2.7  &  \textbf{\cone{87}}/\ctwo{7.3}/\ctwo{11.5}  &  69/4.4/6.9  &  \textbf{\cone{87}}/\ctwo{7.3}/\ctwo{11.5}  &  \ctwo{85}/\cthree{6.9}/\cthree{11.0}  &  \textbf{\cone{87}}/\textbf{\cone{7.4}}/\textbf{\cone{11.6}}  \\
walkbridge+90\%  &  15/-3.2/-6.2  &  \cthree{67}/\ctwo{3.7}/5.1  &  26/0.2/0.6  &  \ctwo{73}/\textbf{\cone{4.8}}/\textbf{\cone{7.8}}  &  36/0.9/1.6  &  \ctwo{73}/\textbf{\cone{4.8}}/\ctwo{7.7}  &  56/\cthree{3.3}/\cthree{5.8}  &  \textbf{\cone{74}}/\textbf{\cone{4.8}}/\textbf{\cone{7.8}}  \\
\hline
pepper+10\%  &  68/4.9/9.6  &  \ctwo{99}/14.8/20.1  &  \ctwo{99}/15.0/21.8  &  \textbf{\cone{100}}/\cthree{20.5}/\cthree{24.9}  &  \cthree{74}/5.4/11.4  &  \textbf{\cone{100}}/\cthree{20.5}/\cthree{24.9}  &  \textbf{\cone{100}}/\ctwo{23.2}/\ctwo{30.5}  &  \textbf{\cone{100}}/\textbf{\cone{23.9}}/\textbf{\cone{31.0}}  \\
pepper+30\%  &  52/3.1/4.8  &  \cthree{98}/14.6/18.3  &  95/10.8/13.6  &  \ctwo{99}/\cthree{16.8}/\cthree{22.9}  &  73/5.4/11.2  &  \ctwo{99}/\cthree{16.8}/\cthree{22.9}  &  \ctwo{99}/\ctwo{17.7}/\ctwo{24.8}  &  \textbf{\cone{100}}/\textbf{\cone{18.5}}/\textbf{\cone{25.6}}  \\
pepper+50\%  &  38/0.3/-1.1  &  \ctwo{97}/12.9/16.1  &  \cthree{84}/6.1/7.0  &  \textbf{\cone{99}}/\ctwo{14.9}/\ctwo{21.5}  &  71/5.2/10.6  &  \textbf{\cone{99}}/\cthree{14.8}/\ctwo{21.5}  &  \textbf{\cone{99}}/14.5/\cthree{21.1}  &  \textbf{\cone{99}}/\textbf{\cone{15.4}}/\textbf{\cone{22.4}}  \\
pepper+70\%  &  25/-1.5/-3.9  &  \cthree{95}/10.6/13.3  &  57/2.1/3.4  &  \textbf{\cone{98}}/\ctwo{12.5}/\ctwo{18.5}  &  61/3.9/7.4  &  \textbf{\cone{98}}/\ctwo{12.5}/\ctwo{18.5}  &  \ctwo{96}/\cthree{11.4}/\cthree{16.9}  &  \textbf{\cone{98}}/\textbf{\cone{12.7}}/\textbf{\cone{18.7}}  \\
pepper+90\%  &  14/-2.7/-5.5  &  \ctwo{89}/\cthree{7.2}/8.5  &  27/0.4/0.6  &  \textbf{\cone{93}}/\ctwo{8.8}/\ctwo{12.7}  &  32/1.2/1.9  &  \textbf{\cone{93}}/\ctwo{8.8}/\cthree{12.5}  &  \cthree{75}/4.8/7.9  &  \textbf{\cone{93}}/\textbf{\cone{9.0}}/\textbf{\cone{12.9}}  \\
\hline
mandrill+10\%  &  77/2.7/4.9  &  \cthree{93}/9.8/11.3  &  90/4.5/6.9  &  \ctwo{97}/\cthree{13.1}/\cthree{14.3}  &  87/4.2/9.2  &  \ctwo{97}/\cthree{13.1}/\cthree{14.3}  &  \textbf{\cone{98}}/\ctwo{14.4}/\ctwo{17.1}  &  \textbf{\cone{98}}/\textbf{\cone{14.5}}/\textbf{\cone{17.2}}  \\
mandrill+30\%  &  61/1.5/2.3  &  \cthree{90}/7.8/9.0  &  75/4.0/5.9  &  \ctwo{92}/\cthree{8.9}/\cthree{10.7}  &  79/3.6/7.2  &  \ctwo{92}/\cthree{8.9}/\cthree{10.7}  &  \textbf{\cone{93}}/\ctwo{9.3}/\ctwo{11.8}  &  \textbf{\cone{93}}/\textbf{\cone{9.4}}/\textbf{\cone{11.9}}  \\
mandrill+50\%  &  44/-0.9/-2.8  &  \cthree{84}/5.7/\cthree{6.6}  &  67/2.7/3.3  &  \ctwo{87}/\cthree{6.6}/\ctwo{8.5}  &  68/2.8/5.2  &  \ctwo{87}/\cthree{6.6}/\ctwo{8.5}  &  \ctwo{87}/\ctwo{6.7}/\textbf{\cone{8.8}}  &  \textbf{\cone{88}}/\textbf{\cone{6.8}}/\textbf{\cone{8.8}}  \\
mandrill+70\%  &  27/-2.7/-5.6  &  \cthree{76}/\cthree{3.8}/\cthree{4.3}  &  48/1.1/1.9  &  \textbf{\cone{80}}/\textbf{\cone{4.9}}/\ctwo{6.5}  &  54/2.0/3.6  &  \textbf{\cone{80}}/\textbf{\cone{4.9}}/\ctwo{6.5}  &  \ctwo{79}/\ctwo{4.8}/\textbf{\cone{6.6}}  &  \textbf{\cone{80}}/\textbf{\cone{4.9}}/\ctwo{6.5}  \\
mandrill+90\%  &  10/-3.8/-7.2  &  \ctwo{63}/\cthree{2.0}/1.9  &  36/0.3/0.6  &  \textbf{\cone{69}}/\textbf{\cone{3.1}}/\ctwo{4.3}  &  35/0.4/0.8  &  \textbf{\cone{69}}/\textbf{\cone{3.1}}/\ctwo{4.3}  &  \cthree{59}/\ctwo{2.4}/\cthree{3.8}  &  \textbf{\cone{69}}/\textbf{\cone{3.1}}/\textbf{\cone{4.4}}  \\
\hline
lake+10\%  &  91/6.6/11.9  &  \ctwo{99}/16.4/22.9  &  \cthree{96}/11.3/17.6  &  \ctwo{99}/\cthree{19.6}/\cthree{25.9}  &  \ctwo{99}/9.0/17.2  &  \ctwo{99}/\cthree{19.6}/25.7  &  \textbf{\cone{100}}/\ctwo{20.3}/\ctwo{27.5}  &  \textbf{\cone{100}}/\textbf{\cone{20.6}}/\textbf{\cone{27.9}}  \\
lake+30\%  &  71/3.9/5.6  &  \cthree{97}/13.6/18.7  &  90/9.1/12.8  &  \ctwo{98}/\cthree{15.0}/\cthree{21.4}  &  \cthree{97}/8.6/16.0  &  \ctwo{98}/\cthree{15.0}/21.3  &  \ctwo{98}/\ctwo{15.1}/\ctwo{21.7}  &  \textbf{\cone{99}}/\textbf{\cone{15.4}}/\textbf{\cone{22.3}}  \\
lake+50\%  &  52/1.2/-0.4  &  \cthree{94}/11.2/15.3  &  76/5.7/6.8  &  \textbf{\cone{97}}/\ctwo{12.5}/\ctwo{18.3}  &  91/7.7/13.6  &  \textbf{\cone{97}}/\ctwo{12.5}/\cthree{18.2}  &  \ctwo{96}/\cthree{12.2}/17.9  &  \textbf{\cone{97}}/\textbf{\cone{12.7}}/\textbf{\cone{18.6}}  \\
lake+70\%  &  33/-0.5/-3.0  &  90/\cthree{9.0}/\cthree{12.1}  &  52/2.4/3.7  &  \ctwo{93}/\textbf{\cone{10.4}}/\textbf{\cone{15.2}}  &  63/5.0/8.2  &  \ctwo{93}/\textbf{\cone{10.4}}/\textbf{\cone{15.2}}  &  \cthree{91}/\ctwo{9.7}/\ctwo{14.4}  &  \textbf{\cone{94}}/\textbf{\cone{10.4}}/\textbf{\cone{15.2}}  \\
lake+90\%  &  18/-1.6/-4.5  &  \cthree{80}/\cthree{6.2}/\cthree{7.5}  &  26/0.5/0.9  &  \textbf{\cone{84}}/\ctwo{7.3}/\ctwo{10.1}  &  25/1.1/1.9  &  \ctwo{83}/\ctwo{7.3}/\ctwo{10.1}  &  51/4.3/7.3  &  \textbf{\cone{84}}/\textbf{\cone{7.4}}/\textbf{\cone{10.2}}  \\
\hline
jetplane+10\%  &  \cthree{49}/2.5/6.0  &  \textbf{\cone{100}}/17.0/23.4  &  \ctwo{98}/11.6/17.3  &  \textbf{\cone{100}}/\cthree{20.4}/\cthree{26.8}  &  39/3.4/8.5  &  \textbf{\cone{100}}/\cthree{20.4}/\cthree{26.8}  &  \textbf{\cone{100}}/\ctwo{20.7}/\ctwo{28.0}  &  \textbf{\cone{100}}/\textbf{\cone{21.3}}/\textbf{\cone{29.2}}  \\
jetplane+30\%  &  39/0.6/1.2  &  \ctwo{98}/13.6/17.9  &  \cthree{93}/8.3/10.4  &  \textbf{\cone{99}}/\ctwo{15.5}/\ctwo{21.9}  &  40/3.4/8.3  &  \textbf{\cone{99}}/\ctwo{15.5}/\ctwo{21.9}  &  \textbf{\cone{99}}/\cthree{15.3}/\cthree{21.6}  &  \textbf{\cone{99}}/\textbf{\cone{15.9}}/\textbf{\cone{22.7}}  \\
jetplane+50\%  &  33/-1.4/-4.1  &  \ctwo{96}/10.9/14.1  &  \cthree{79}/4.0/5.1  &  \textbf{\cone{98}}/\ctwo{12.7}/\ctwo{18.4}  &  39/3.1/7.2  &  \textbf{\cone{98}}/\ctwo{12.7}/\ctwo{18.4}  &  \textbf{\cone{98}}/\cthree{12.1}/\cthree{17.3}  &  \textbf{\cone{98}}/\textbf{\cone{12.9}}/\textbf{\cone{18.5}}  \\
jetplane+70\%  &  30/-2.8/-6.4  &  \cthree{93}/8.5/\cthree{10.5}  &  53/0.3/1.2  &  \textbf{\cone{96}}/\ctwo{10.2}/\textbf{\cone{14.6}}  &  32/1.2/3.0  &  \textbf{\cone{96}}/\ctwo{10.2}/\textbf{\cone{14.6}}  &  \ctwo{94}/\cthree{9.2}/\ctwo{13.3}  &  \textbf{\cone{96}}/\textbf{\cone{10.3}}/\textbf{\cone{14.6}}  \\
jetplane+90\%  &  28/-3.7/-7.9  &  \ctwo{87}/\cthree{5.6}/\cthree{6.0}  &  26/-1.7/-2.1  &  \textbf{\cone{89}}/\ctwo{6.6}/\ctwo{8.6}  &  29/-1.9/-2.8  &  \textbf{\cone{89}}/\ctwo{6.6}/\ctwo{8.6}  &  \cthree{54}/2.4/4.8  &  \textbf{\cone{89}}/\textbf{\cone{6.8}}/\textbf{\cone{8.7}}  \\
\hline\hline
\multicolumn{9}{|c|}{Mixed Impulse Noise (Half Random-Value Noise and Half Salt-and-Pepper Noise)} \\
\hline
walkbridge+10\%  &  91/6.1/10.1  &  \cthree{93}/10.6/14.7  &  91/7.5/12.3  &  \ctwo{96}/\cthree{12.6}/13.3  &  \ctwo{96}/12.5/\cthree{16.0}  &  \ctwo{96}/\cthree{12.6}/13.3  &  \textbf{\cone{98}}/\ctwo{14.8}/\ctwo{17.8}  &  \textbf{\cone{98}}/\textbf{\cone{15.1}}/\textbf{\cone{17.9}}  \\
walkbridge+30\%  &  73/3.6/6.7  &  \cthree{90}/\cthree{8.4}/11.8  &  83/6.3/10.3  &  88/6.6/8.3  &  89/\ctwo{8.6}/\cthree{12.2}  &  \ctwo{92}/\ctwo{8.6}/\cthree{12.2}  &  \textbf{\cone{93}}/\textbf{\cone{10.2}}/\textbf{\cone{13.5}}  &  \textbf{\cone{93}}/\textbf{\cone{10.2}}/\ctwo{12.9}  \\
walkbridge+50\%  &  55/1.5/1.9  &  81/\cthree{5.7}/7.0  &  70/4.3/6.8  &  76/3.5/5.7  &  78/\cthree{5.7}/8.7  &  \cthree{85}/\ctwo{6.3}/\cthree{10.0}  &  \ctwo{86}/\textbf{\cone{7.6}}/\textbf{\cone{10.8}}  &  \textbf{\cone{87}}/\textbf{\cone{7.6}}/\ctwo{10.1}  \\
walkbridge+70\%  &  37/-0.5/-1.8  &  63/2.4/1.9  &  50/2.0/2.9  &  58/1.9/3.2  &  56/2.8/\cthree{4.9}  &  \cthree{72}/\cthree{4.4}/\ctwo{7.2}  &  \ctwo{74}/\ctwo{5.1}/\textbf{\cone{7.9}}  &  \textbf{\cone{80}}/\textbf{\cone{5.7}}/\textbf{\cone{7.9}}  \\
walkbridge+90\%  &  21/-1.9/-4.0  &  34/-0.6/-2.1  &  30/0.1/0.4  &  34/0.3/0.5  &  31/0.6/1.3  &  \cthree{38}/\cthree{1.2}/\cthree{2.0}  &  \ctwo{40}/\ctwo{1.3}/\ctwo{2.3}  &  \textbf{\cone{63}}/\textbf{\cone{3.3}}/\textbf{\cone{4.9}}  \\
\hline
pepper+10\%  &  68/5.0/9.7  &  \cthree{98}/13.9/19.5  &  \ctwo{99}/\cthree{15.0}/\cthree{22.0}  &  \cthree{98}/14.3/16.0  &  74/5.4/11.3  &  \ctwo{99}/14.4/19.9  &  \textbf{\cone{100}}/\textbf{\cone{21.0}}/\textbf{\cone{25.6}}  &  \ctwo{99}/\ctwo{19.9}/\ctwo{23.4}  \\
pepper+30\%  &  54/3.7/6.8  &  \cthree{97}/12.7/16.0  &  96/11.4/15.4  &  91/7.5/10.8  &  72/5.3/10.8  &  \ctwo{98}/\cthree{12.8}/\ctwo{18.5}  &  \textbf{\cone{99}}/\textbf{\cone{15.8}}/\textbf{\cone{20.7}}  &  \ctwo{98}/\ctwo{14.9}/\cthree{18.4}  \\
pepper+50\%  &  41/1.8/2.3  &  \ctwo{92}/\cthree{8.5}/8.6  &  \cthree{86}/7.0/8.9  &  80/4.5/7.0  &  68/4.8/9.5  &  \textbf{\cone{97}}/\ctwo{11.2}/\ctwo{16.1}  &  \textbf{\cone{97}}/\textbf{\cone{12.6}}/\textbf{\cone{17.1}}  &  \textbf{\cone{97}}/\textbf{\cone{12.6}}/\cthree{15.7}  \\
pepper+70\%  &  29/-0.1/-1.2  &  73/3.6/2.4  &  62/3.0/3.6  &  63/2.5/3.8  &  54/3.3/5.9  &  \cthree{90}/\cthree{8.1}/\cthree{10.7}  &  \ctwo{92}/\ctwo{9.1}/\ctwo{12.5}  &  \textbf{\cone{94}}/\textbf{\cone{10.1}}/\textbf{\cone{12.8}}  \\
pepper+90\%  &  19/-1.4/-3.4  &  39/-0.2/-2.0  &  33/0.4/0.5  &  37/0.6/0.7  &  31/1.0/1.5  &  \ctwo{53}/\cthree{2.1}/\cthree{2.5}  &  \cthree{49}/\ctwo{2.2}/\ctwo{2.9}  &  \textbf{\cone{82}}/\textbf{\cone{5.6}}/\textbf{\cone{6.6}}  \\
\hline
mandrill+10\%  &  76/3.0/5.3  &  86/6.8/8.3  &  \cthree{91}/5.5/6.8  &  \ctwo{95}/10.4/10.1  &  83/3.6/7.3  &  \ctwo{95}/\cthree{10.5}/\cthree{10.3}  &  \textbf{\cone{96}}/\textbf{\cone{12.1}}/\textbf{\cone{12.4}}  &  \textbf{\cone{96}}/\ctwo{11.7}/\ctwo{11.2}  \\
mandrill+30\%  &  63/1.8/3.4  &  82/\cthree{5.4}/6.6  &  74/3.9/6.0  &  \cthree{85}/5.3/5.1  &  73/2.9/5.3  &  \ctwo{88}/\ctwo{6.5}/\cthree{7.4}  &  \textbf{\cone{89}}/\textbf{\cone{7.3}}/\textbf{\cone{8.1}}  &  \textbf{\cone{89}}/\textbf{\cone{7.3}}/\ctwo{7.5}  \\
mandrill+50\%  &  47/0.6/0.6  &  75/\cthree{3.7}/4.0  &  67/3.0/4.4  &  74/2.5/3.8  &  61/2.2/3.9  &  \cthree{78}/\ctwo{4.4}/\ctwo{5.6}  &  \ctwo{80}/\textbf{\cone{5.0}}/\textbf{\cone{5.9}}  &  \textbf{\cone{81}}/\textbf{\cone{5.0}}/\cthree{5.3}  \\
mandrill+70\%  &  32/-1.0/-2.6  &  60/1.3/0.2  &  53/1.5/1.8  &  58/1.3/2.1  &  48/1.4/2.6  &  \cthree{68}/\cthree{2.9}/\ctwo{4.2}  &  \ctwo{69}/\ctwo{3.3}/\textbf{\cone{4.3}}  &  \textbf{\cone{73}}/\textbf{\cone{3.5}}/\cthree{3.9}  \\
mandrill+90\%  &  20/-2.4/-4.9  &  35/-1.2/-3.3  &  36/0.3/0.5  &  37/0.2/0.1  &  33/0.3/0.6  &  \ctwo{46}/\ctwo{1.1}/\ctwo{1.8}  &  \cthree{45}/\cthree{1.0}/\cthree{1.3}  &  \textbf{\cone{62}}/\textbf{\cone{2.2}}/\textbf{\cone{2.8}}  \\
\hline
lake+10\%  &  91/6.8/12.0  &  \ctwo{98}/14.6/\cthree{20.5}  &  96/11.3/17.7  &  \ctwo{98}/\cthree{15.0}/15.5  &  \cthree{97}/8.7/16.1  &  \ctwo{98}/\cthree{15.0}/19.4  &  \textbf{\cone{99}}/\textbf{\cone{18.0}}/\textbf{\cone{22.2}}  &  \textbf{\cone{99}}/\ctwo{17.9}/\ctwo{21.2}  \\
lake+30\%  &  73/4.3/7.6  &  \ctwo{95}/\cthree{11.7}/\cthree{15.7}  &  91/9.3/13.7  &  90/8.0/10.5  &  \cthree{92}/7.9/13.8  &  \textbf{\cone{96}}/11.0/\ctwo{16.4}  &  \textbf{\cone{96}}/\textbf{\cone{13.1}}/\textbf{\cone{17.2}}  &  \textbf{\cone{96}}/\ctwo{12.8}/15.6  \\
lake+50\%  &  55/2.3/2.7  &  \ctwo{87}/7.9/9.0  &  75/6.1/8.6  &  78/4.8/7.2  &  \cthree{82}/6.6/11.0  &  \textbf{\cone{92}}/\cthree{9.3}/\ctwo{13.1}  &  \textbf{\cone{92}}/\textbf{\cone{10.4}}/\textbf{\cone{14.1}}  &  \textbf{\cone{92}}/\ctwo{10.0}/\cthree{12.2}  \\
lake+70\%  &  37/0.6/-0.6  &  66/3.7/3.1  &  44/2.8/3.7  &  58/2.6/4.1  &  48/3.7/6.2  &  \cthree{82}/\cthree{7.0}/\ctwo{9.8}  &  \ctwo{83}/\ctwo{7.7}/\textbf{\cone{10.8}}  &  \textbf{\cone{87}}/\textbf{\cone{7.9}}/\cthree{9.4}  \\
lake+90\%  &  22/-0.6/-2.7  &  \ctwo{34}/0.4/-1.1  &  20/0.6/0.7  &  30/0.8/1.0  &  24/0.8/1.5  &  22/\cthree{1.5}/\cthree{2.4}  &  \cthree{33}/\ctwo{2.0}/\ctwo{3.1}  &  \textbf{\cone{74}}/\textbf{\cone{5.3}}/\textbf{\cone{6.0}}  \\
\hline
jetplane+10\%  &  \cthree{44}/2.6/6.0  &  \textbf{\cone{99}}/\cthree{15.4}/\textbf{\cone{20.8}}  &  \ctwo{98}/11.6/17.5  &  \textbf{\cone{99}}/13.9/13.3  &  39/3.4/8.3  &  \textbf{\cone{99}}/13.9/\ctwo{19.3}  &  \textbf{\cone{99}}/\textbf{\cone{17.6}}/\textbf{\cone{20.8}}  &  \textbf{\cone{99}}/\ctwo{16.8}/\cthree{18.5}  \\
jetplane+30\%  &  36/0.8/2.5  &  \ctwo{97}/\cthree{11.6}/\cthree{14.2}  &  \cthree{94}/8.8/12.3  &  91/6.3/8.2  &  38/3.2/7.7  &  \textbf{\cone{98}}/11.0/\textbf{\cone{16.6}}  &  \textbf{\cone{98}}/\textbf{\cone{13.1}}/\ctwo{16.4}  &  \textbf{\cone{98}}/\ctwo{12.6}/14.1  \\
jetplane+50\%  &  30/-0.8/-1.6  &  \ctwo{90}/6.6/6.2  &  \cthree{79}/4.5/5.8  &  \cthree{79}/2.8/4.2  &  37/2.8/6.1  &  \textbf{\cone{95}}/\cthree{9.0}/\ctwo{12.7}  &  \textbf{\cone{95}}/\textbf{\cone{10.0}}/\textbf{\cone{13.0}}  &  \textbf{\cone{95}}/\ctwo{9.7}/\cthree{10.7}  \\
jetplane+70\%  &  25/-2.1/-4.5  &  68/1.7/-0.1  &  45/0.6/0.5  &  60/0.4/0.9  &  25/0.7/1.9  &  \ctwo{88}/\cthree{6.3}/\cthree{7.9}  &  \cthree{87}/\ctwo{6.6}/\textbf{\cone{8.8}}  &  \textbf{\cone{91}}/\textbf{\cone{7.3}}/\ctwo{8.0}  \\
jetplane+90\%  &  22/-3.1/-6.4  &  \ctwo{34}/-1.8/-4.4  &  19/-1.8/-2.4  &  30/-1.5/-2.3  &  16/-2.1/-3.0  &  19/\cthree{-0.8}/\cthree{-0.9}  &  \cthree{32}/\ctwo{-0.2}/\ctwo{-0.1}  &  \textbf{\cone{79}}/\textbf{\cone{4.2}}/\textbf{\cone{4.4}}  \\
\hline
\end{tabular}}
\end{center}
\end{table*}

\subsection{Convergence of $\ell_0TV$-\text{PADMM}}
Here, we verify the convergence property of our $\ell_0TV$-\text{PADMM} method on denoising and deblurring problems by considering the `cameraman' image subject to 30\% random-valued impulse noise. We set $\lambda=8$ for this problem. We record the objective and SNR values for $\ell_0TV$-\text{PADMM} at every iteration $k$ and plot these results in Figure \ref{fig:exp:convergence}.



We make two important observations from these results. \textbf{(i)}) The objective value (or the SNR value) does not necessarily decrease (or increase) monotonically, and we attribute this to the non-convexity of the optimization problem and the dynamic updates of the penalty factor in Algorithm \ref{alg:proximal:admm}. \textbf{(ii)} The objective and SNR values stabilize after the $120$th iteration, which means that our algorithm has converged, and the increase of the SNR value is negligible after the $80$th iteration. This implies that one may use a looser stopping criterion without sacrificing much restoration quality.

\begin{figure}[!t]
\begin{center}
\subfloat[\figsizetwo Random-Value]{\includegraphics[width=\figurewidth,height=\figureheight]{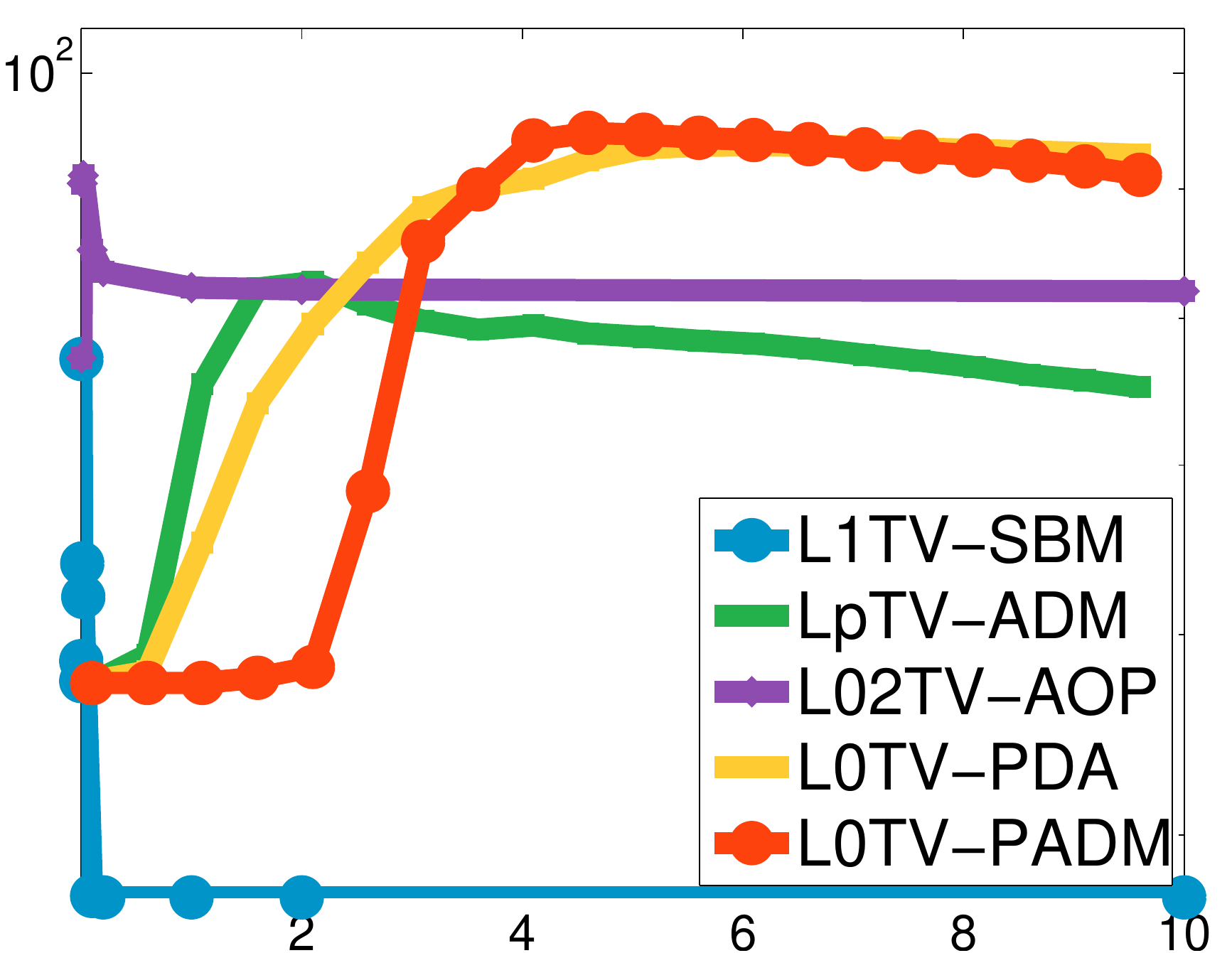}}\ghs
\subfloat[\figsizetwo Salt-and-Pepper]{\includegraphics[width=\figurewidth,height=\figureheight]{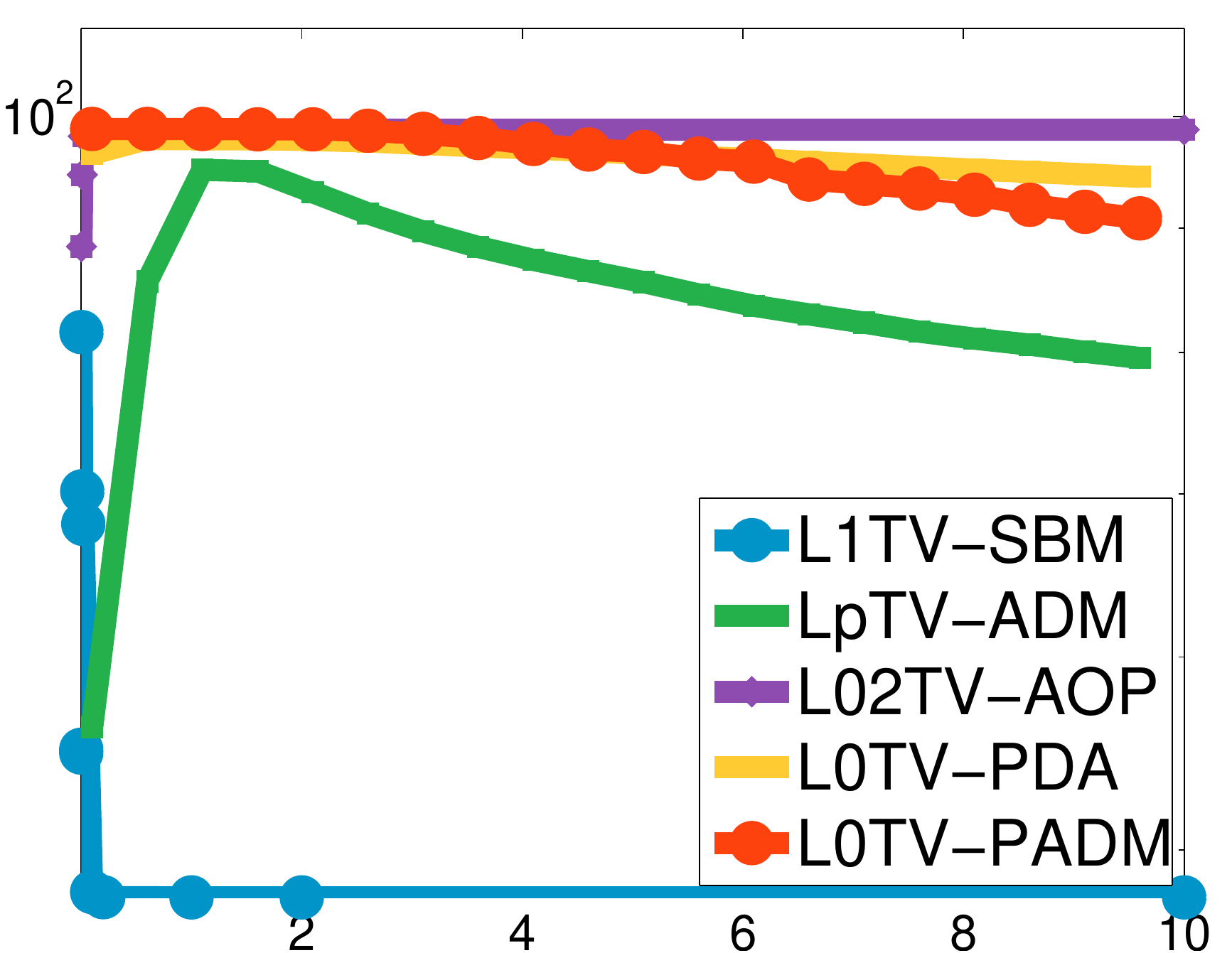}}\ghs
\subfloat[\figsizetwo Mixed]{\includegraphics[width=\figurewidth, height=\figureheight]{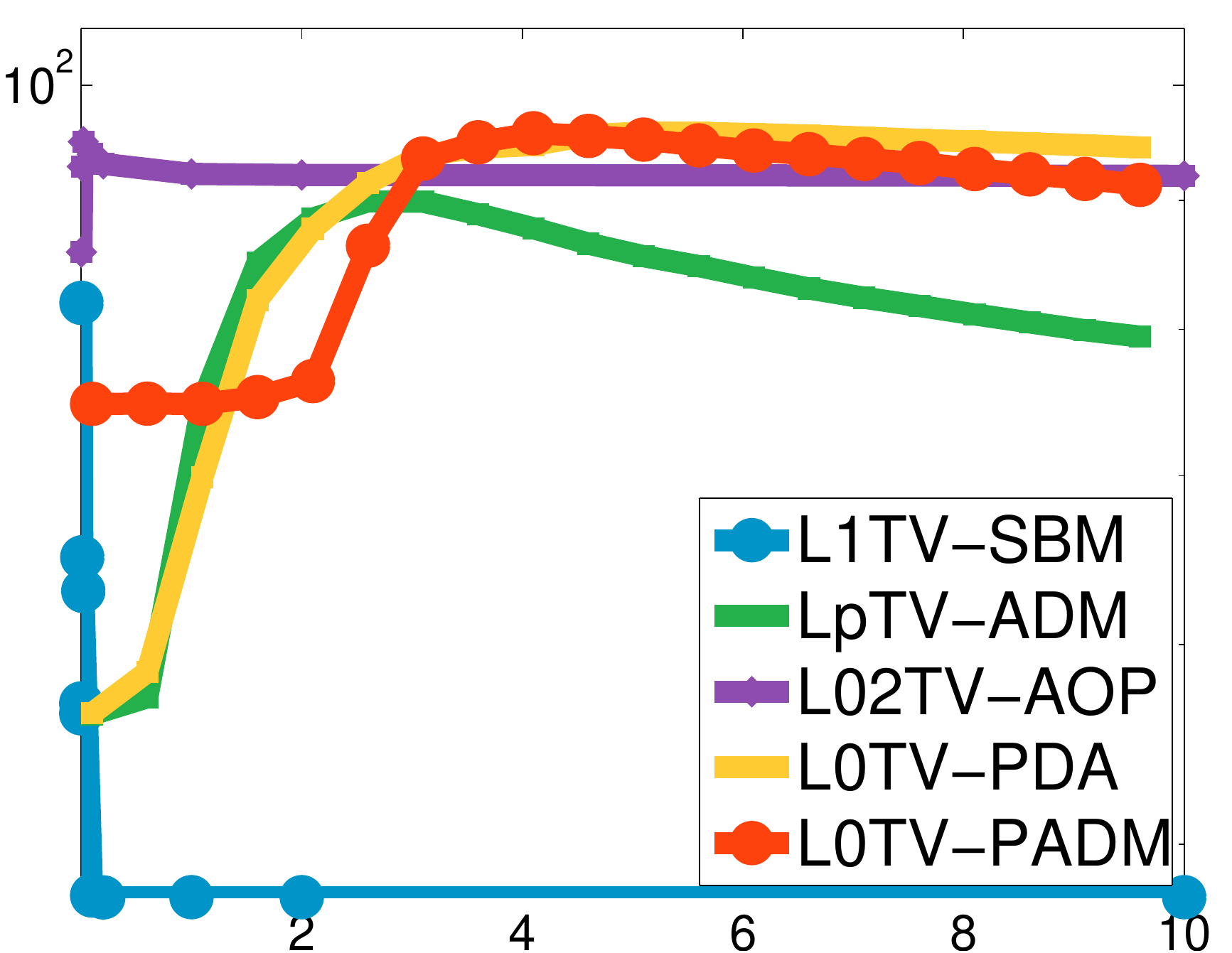}}
\\
\vspace{-8pt}
\subfloat[\figsizetwo Random-Value]{\includegraphics[width=\figurewidth, height=\figureheight]{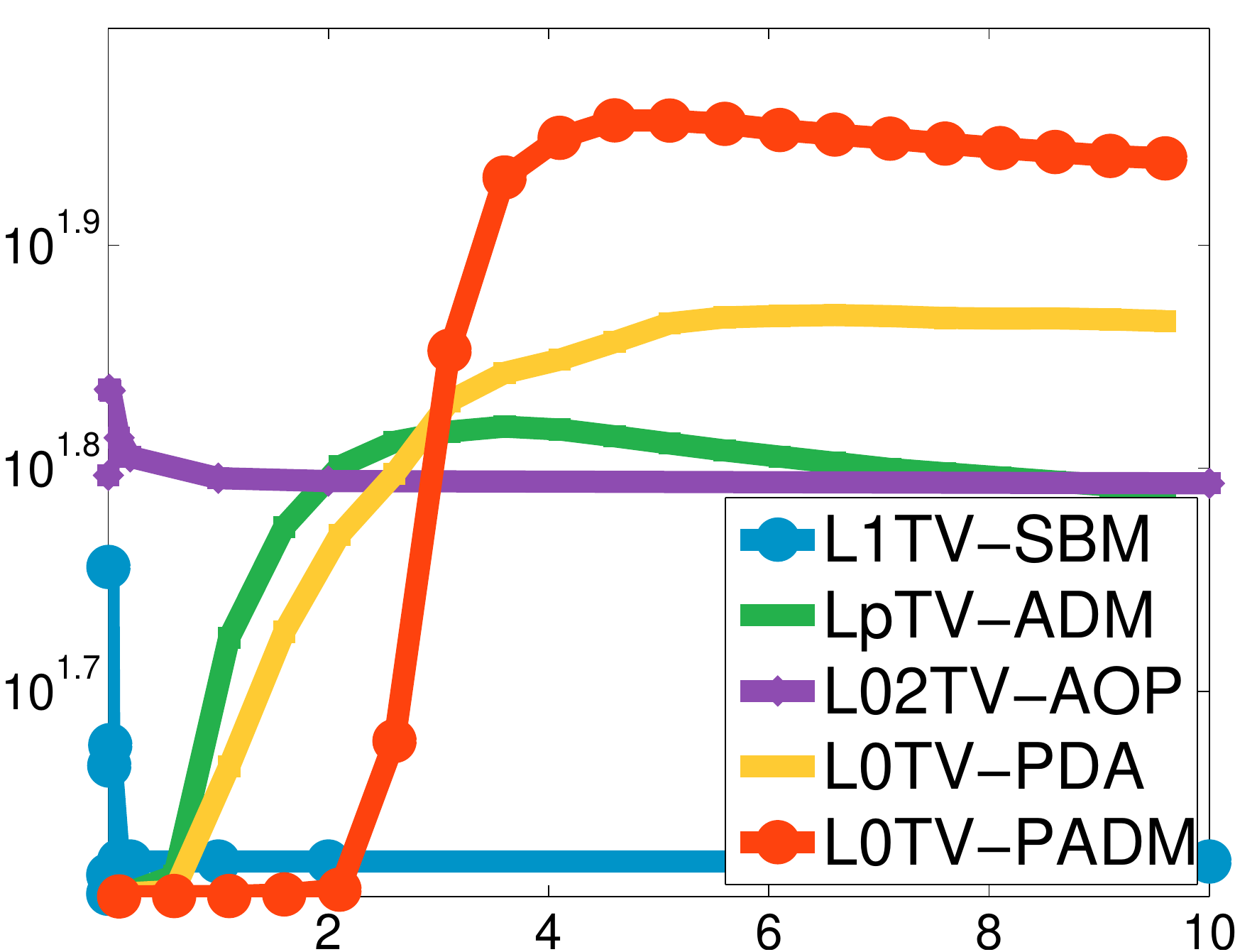}}\ghs
\subfloat[\figsizetwo Salt-and-Pepper]{\includegraphics[width=\figurewidth, height=\figureheight]{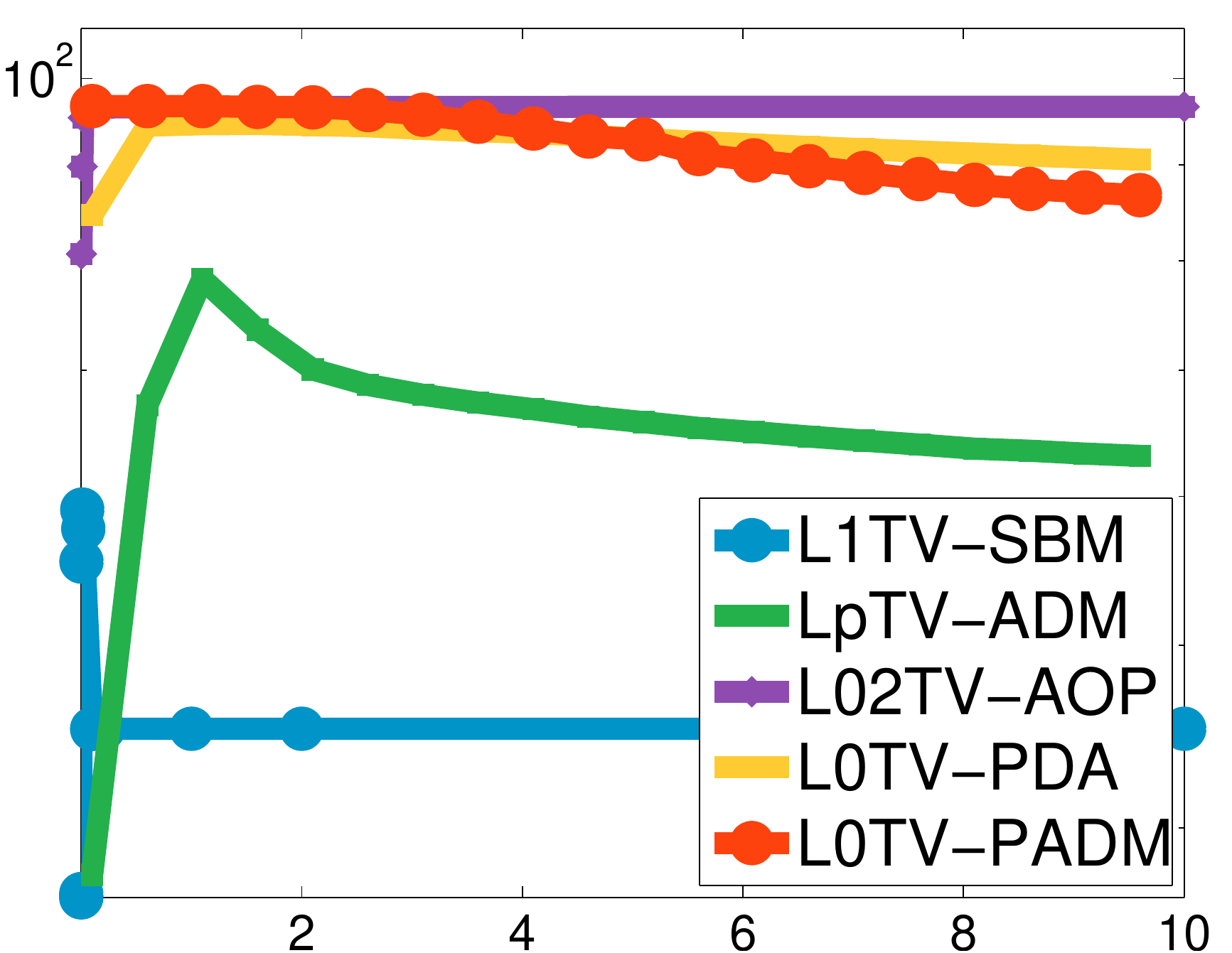}}\ghs
\subfloat[\figsizetwo Mixed]{\includegraphics[width=\figurewidth, height=\figureheight]{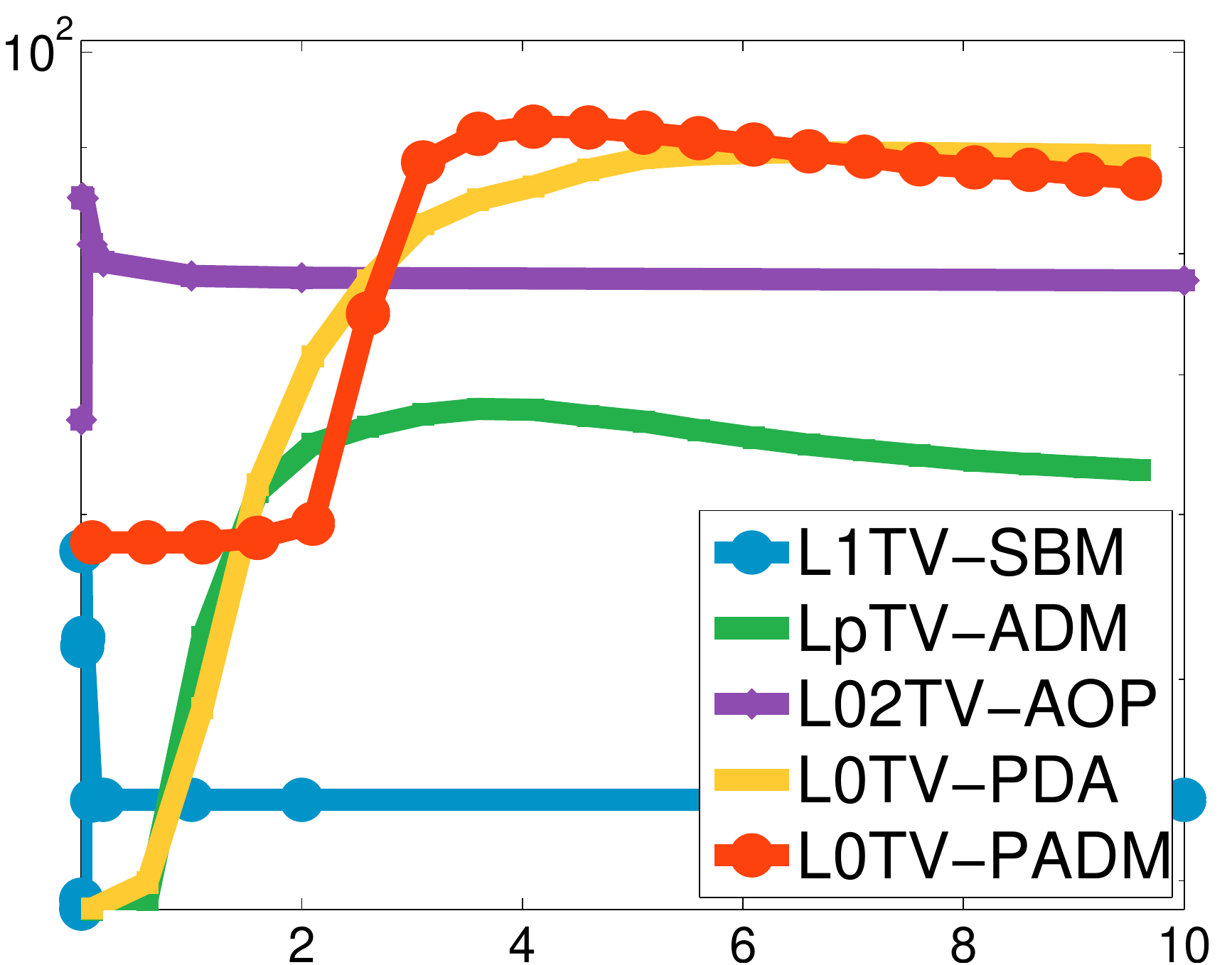}}
\\
\vspace{-8pt}
\subfloat[\figsizetwo Random-Value]{\includegraphics[width=\figurewidth, height=\figureheight]{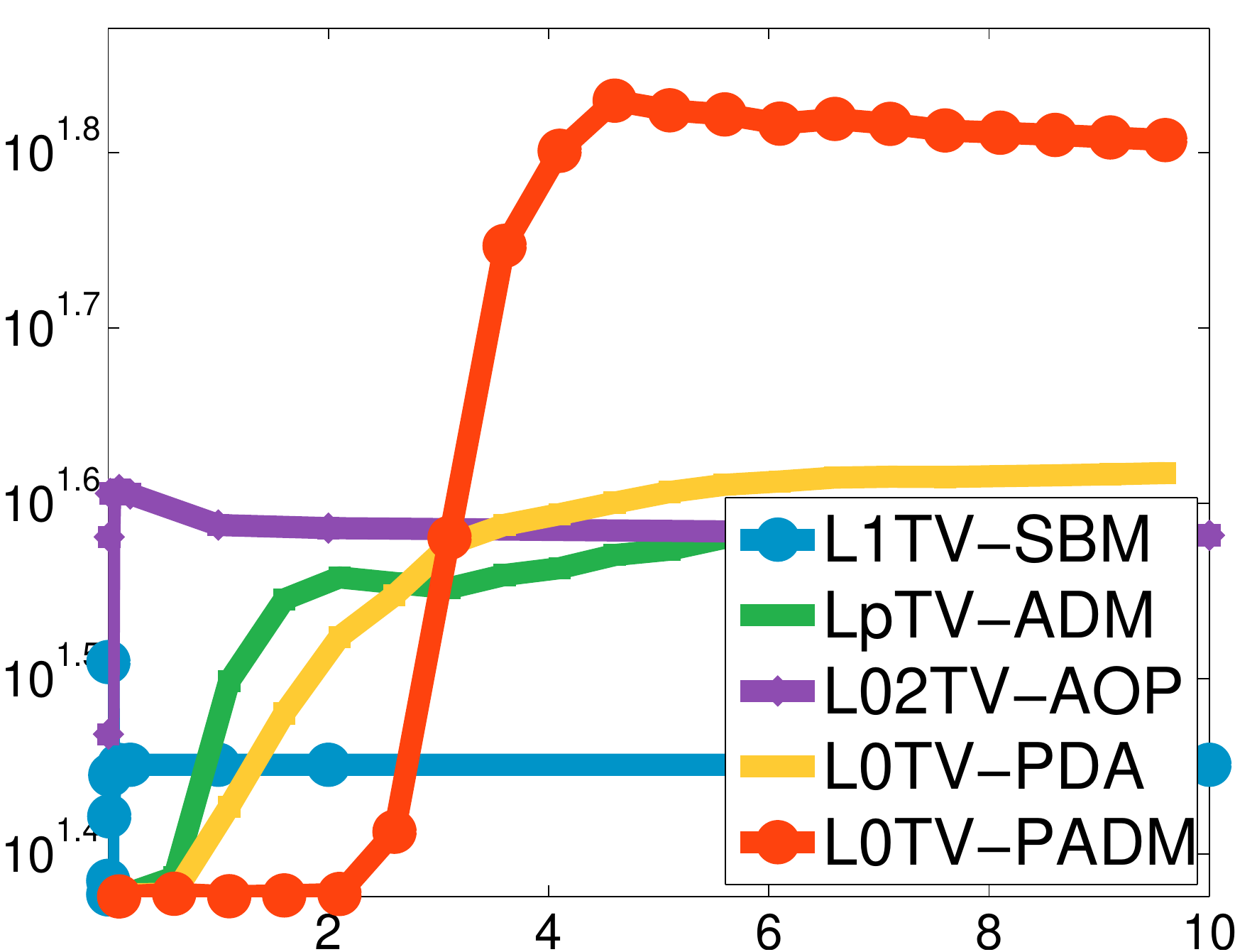}}\ghs
\subfloat[\figsizetwo Salt-and-Pepper]{\includegraphics[width=\figurewidth, height=\figureheight]{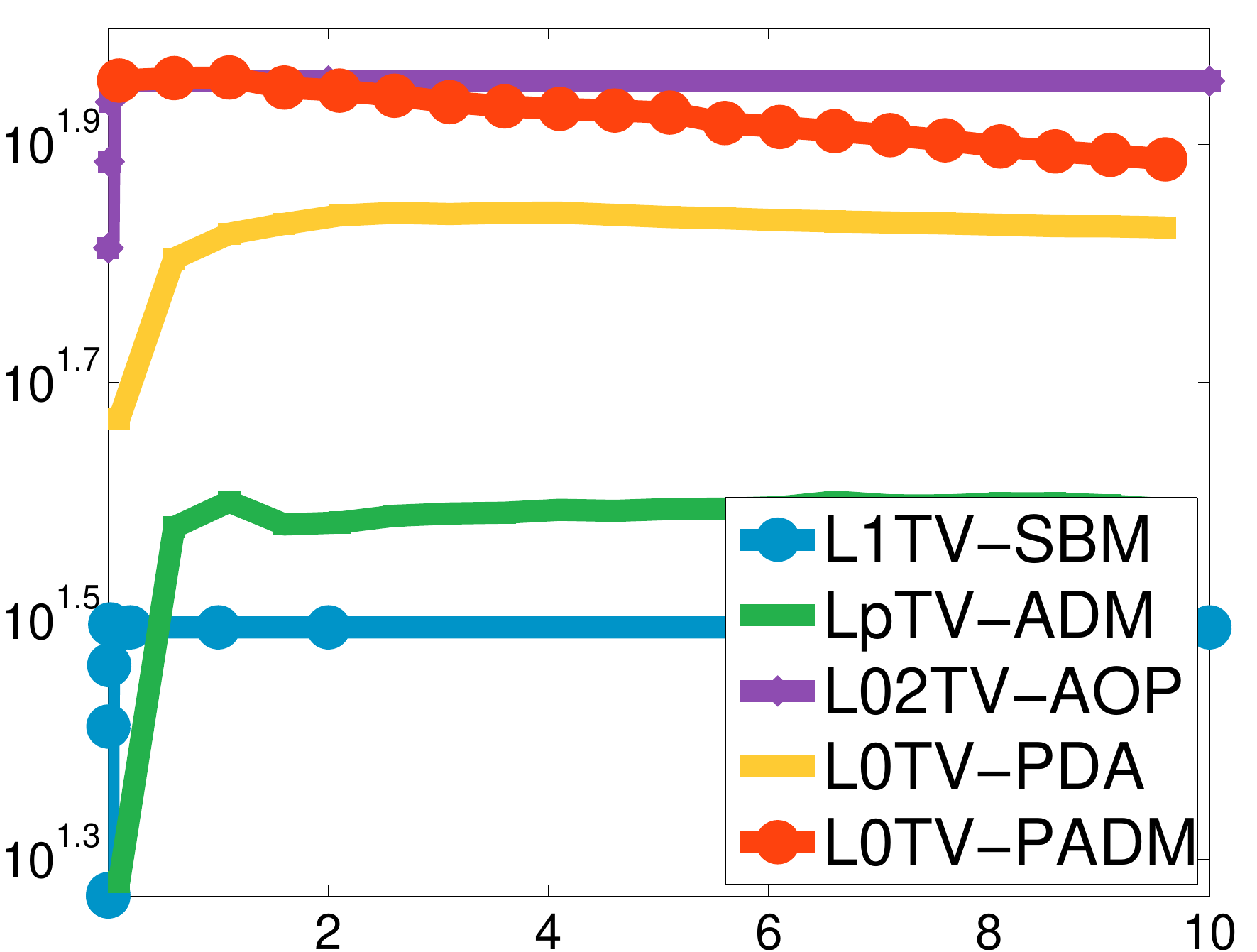}}\ghs
\subfloat[\figsizetwo Mixed]{\includegraphics[width=\figurewidth, height=\figureheight]{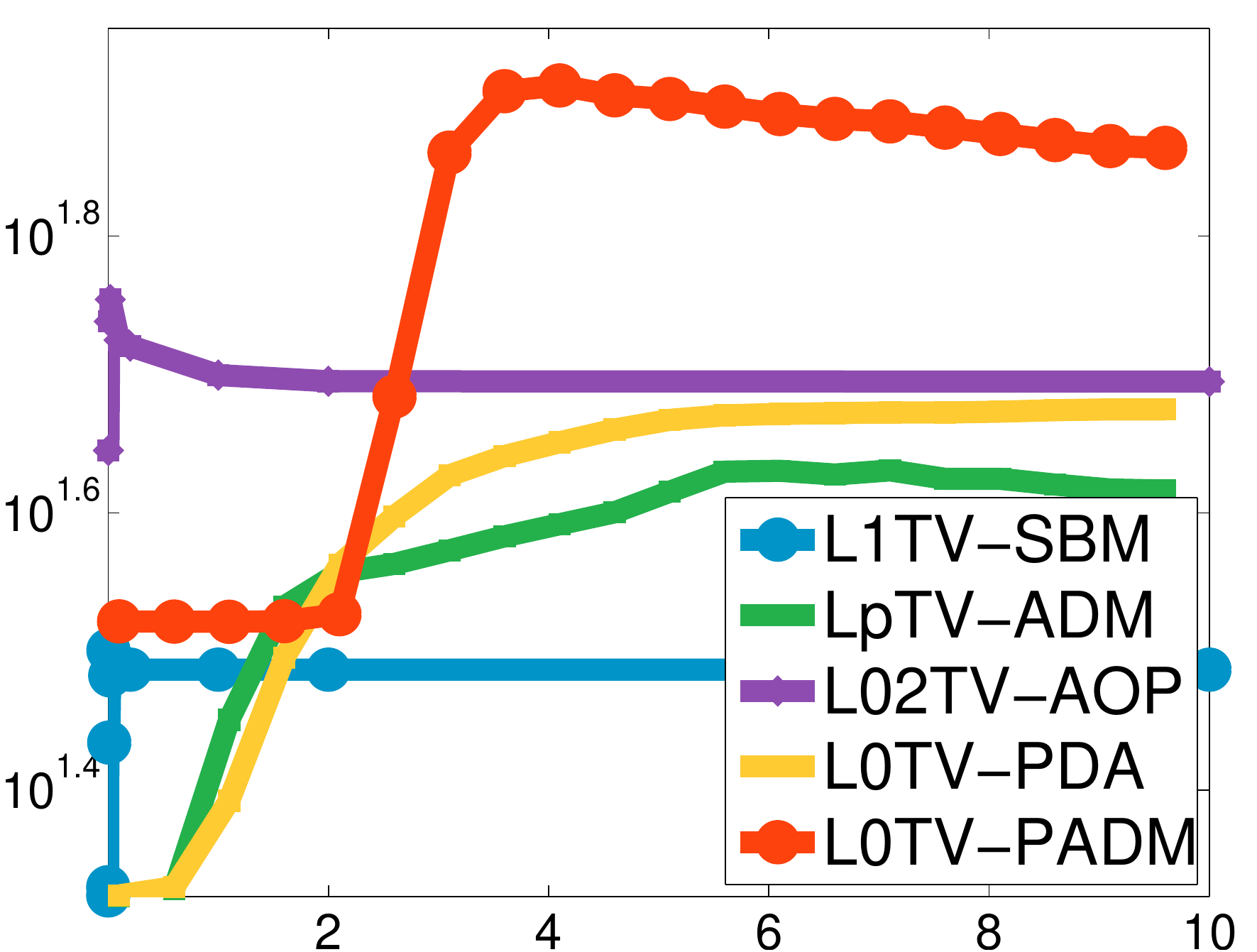}}

\caption{Image denoising with varying the tuning parameter $\lambda$ in (\ref{eq:l0tv:2}) on `cameraman' image. First row: noise level = 50\%. Second row: noise level = 70\%. Third row: noise level = 90\%.}
\label{fig:denoising:lambda}
\end{center}
\end{figure}

\subsection{General Image Denoising Problems}
In this subsection, we compare the performance of all 6 methods on general denoising problems. Table \ref{denoising:rv} shows image recovery results when random-value or salt-and-pepper or mixed impulse noise is added. Figure \ref{fig:denoising:lambda} shows image recovery results with varying the regularization parameter $\lambda$. For $\ell_{02}TV$ model in (\ref{eq:l02tv}), the parameter $\chi$ is scaled to the range $[0,10]$ for better visualization. We make the following interesting observations. \textbf{(i)} The $\ell_{02}TV$-AOP method greatly improves upon $\ell_1TV$-\text{SBM}, \text{MFM} and \text{TSM}, by a large margin. These results are consistent with the reported results in \cite{yan2013restoration}. \textbf{(ii)} The $\ell_0TV$-PDA method outperforms $\ell_{02}TV$-AOP in most test cases because it adopts the $\ell_0$-norm in the data fidelity term. \textbf{(iii)} In the case of random-value impulse noise, our $\ell_0TV$-PADMM method is better than $\ell_0TV$-PDA in $SNR_0$ value while it is comparable to $\ell_0TV$-PDA in $SNR_1$ and $SNR_2$. On the other hand, when salt-and-pepper impulse noise is added, we find that $\ell_0TV$-PADMM outperforms $\ell_0TV$-PDA in most test cases. Interestingly, the performance gap between $\ell_0TV$-PADMM and $\ell_0TV$-PDA grows larger, as the noise level increases. \textbf{(iv)} For the same noise level,  $\ell_0TV$-PADMM achieves better recovery performance in the presence of salt-and-pepper impulse noise than random-valued impulse noise. This is primarily due to the fact that random-valued noise can take any value between 0 and 1, thus, making it more difficult to detect which pixels are corrupted.

\begin{figure}[!t]
\begin{center}
\subfloat[\figsizetwo Random-Value]{\includegraphics[width=\figurewidth, height=\figureheight]{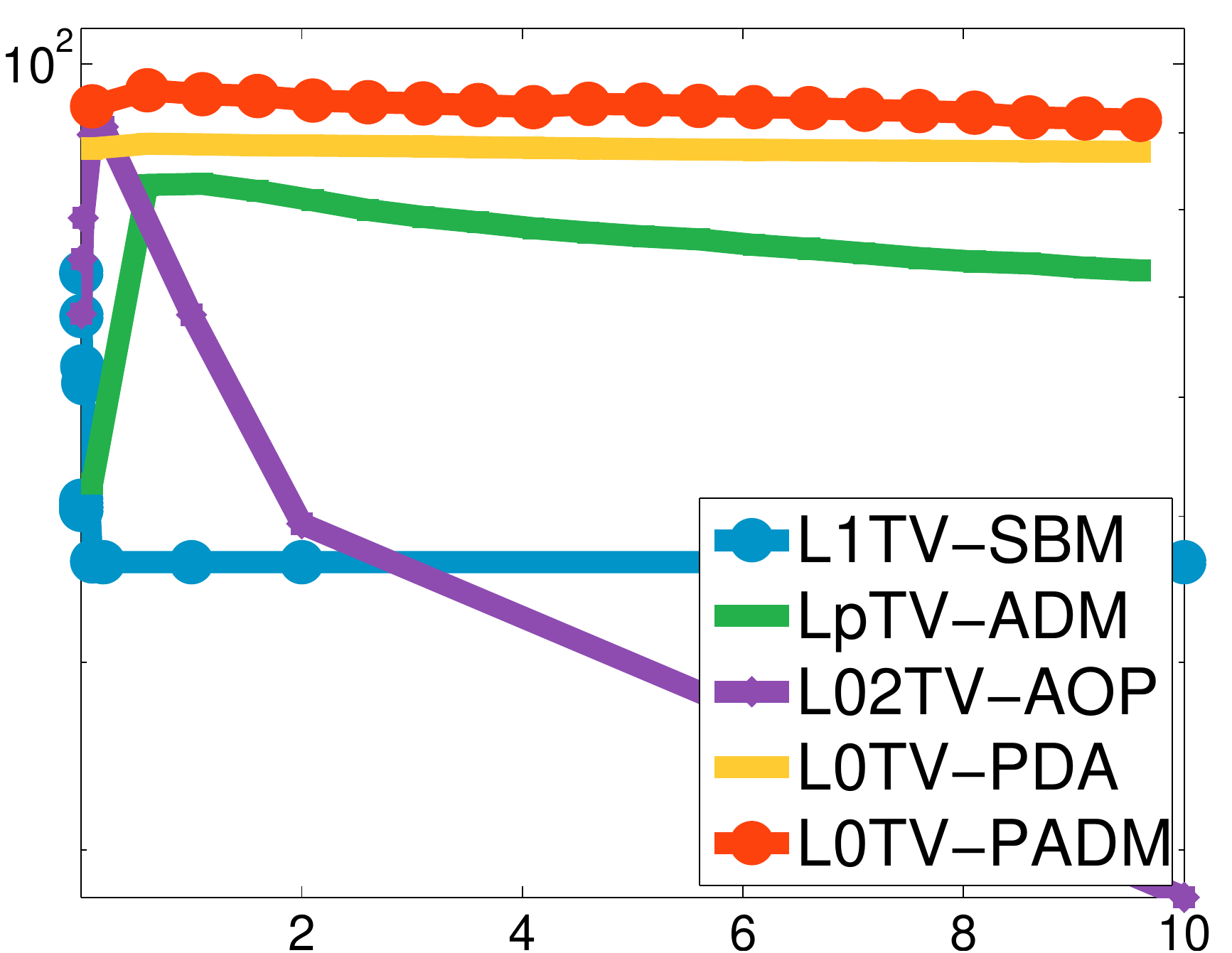}}\ghs
\subfloat[\figsizetwo Salt-and-Pepper]{\includegraphics[width=\figurewidth, height=\figureheight]{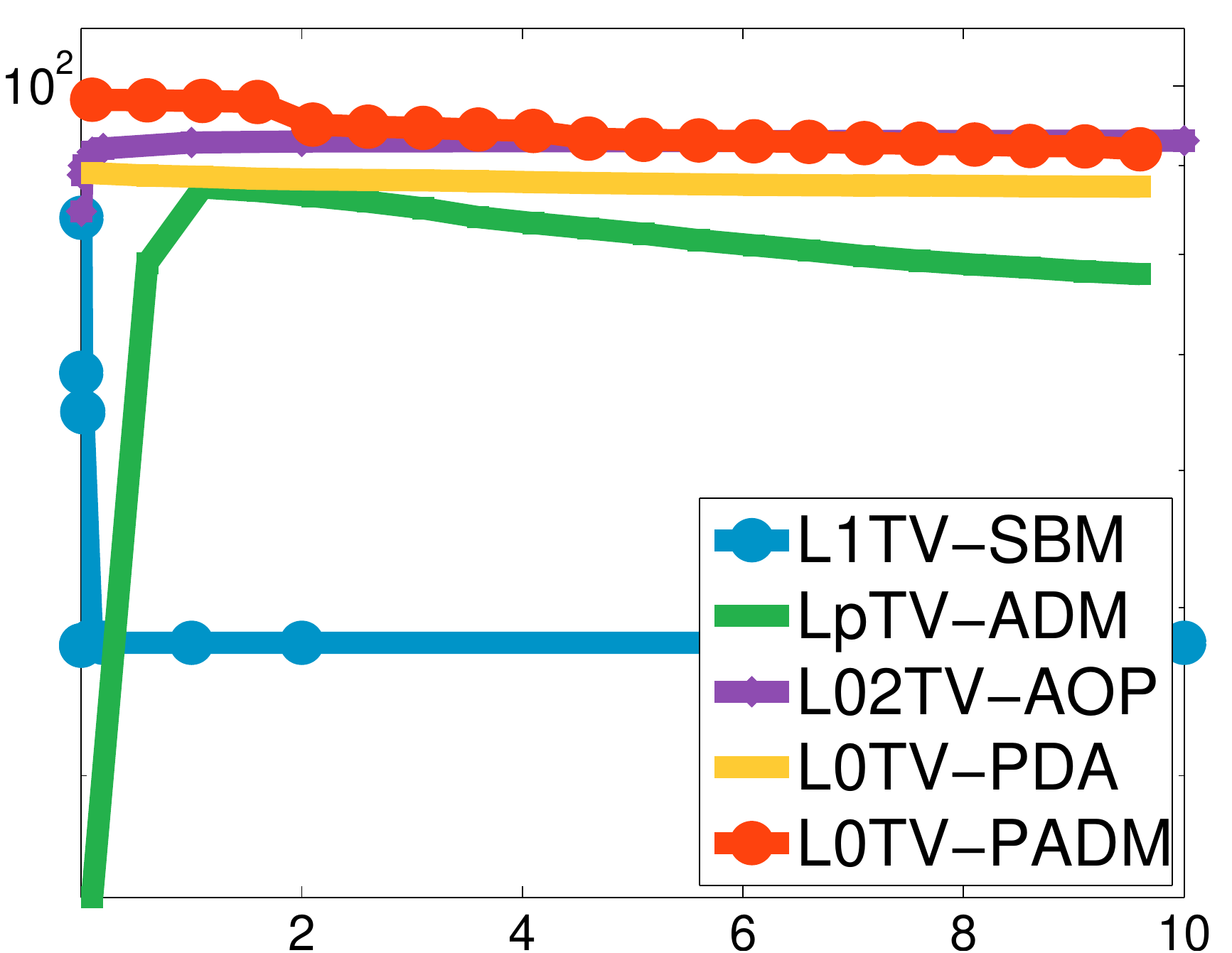}}\ghs
\subfloat[\figsizetwo Mixed]{\includegraphics[width=\figurewidth, height=\figureheight]{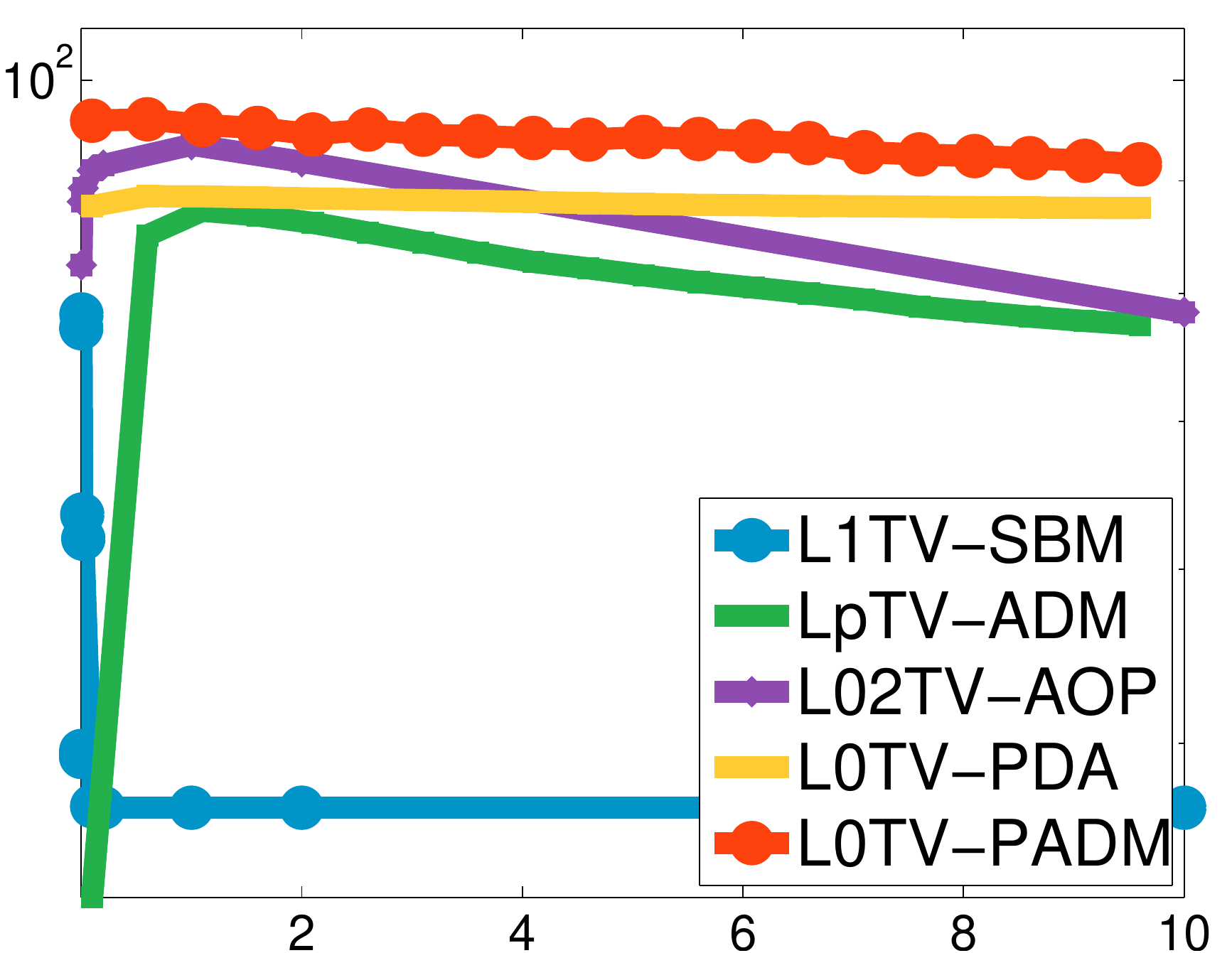}}
\\
\vspace{-8pt}
\subfloat[\figsizetwo Random-Value]{\includegraphics[width=\figurewidth, height=\figureheight]{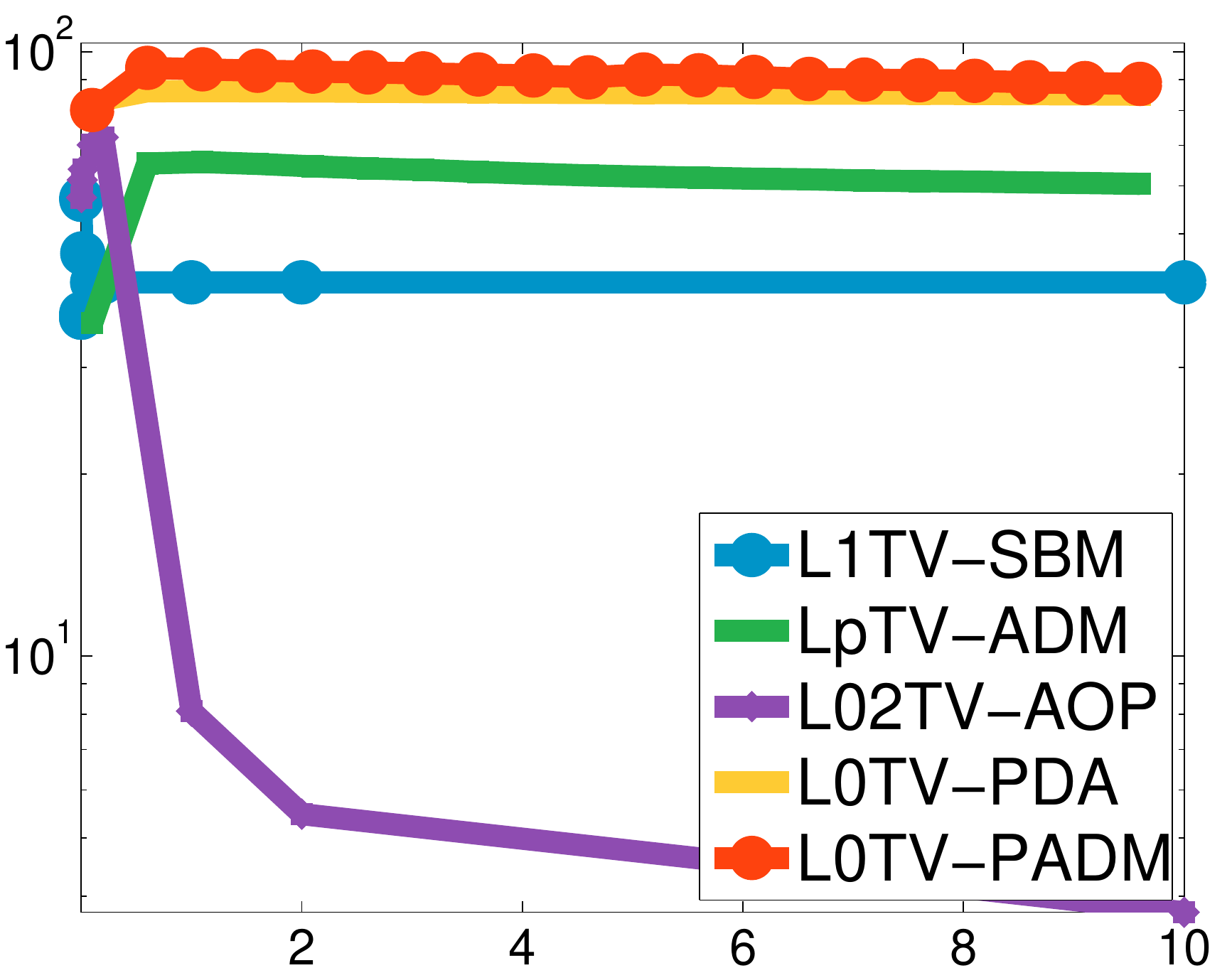}}\ghs
\subfloat[\figsizetwo Salt-and-Pepper]{\includegraphics[width=\figurewidth, height=\figureheight]{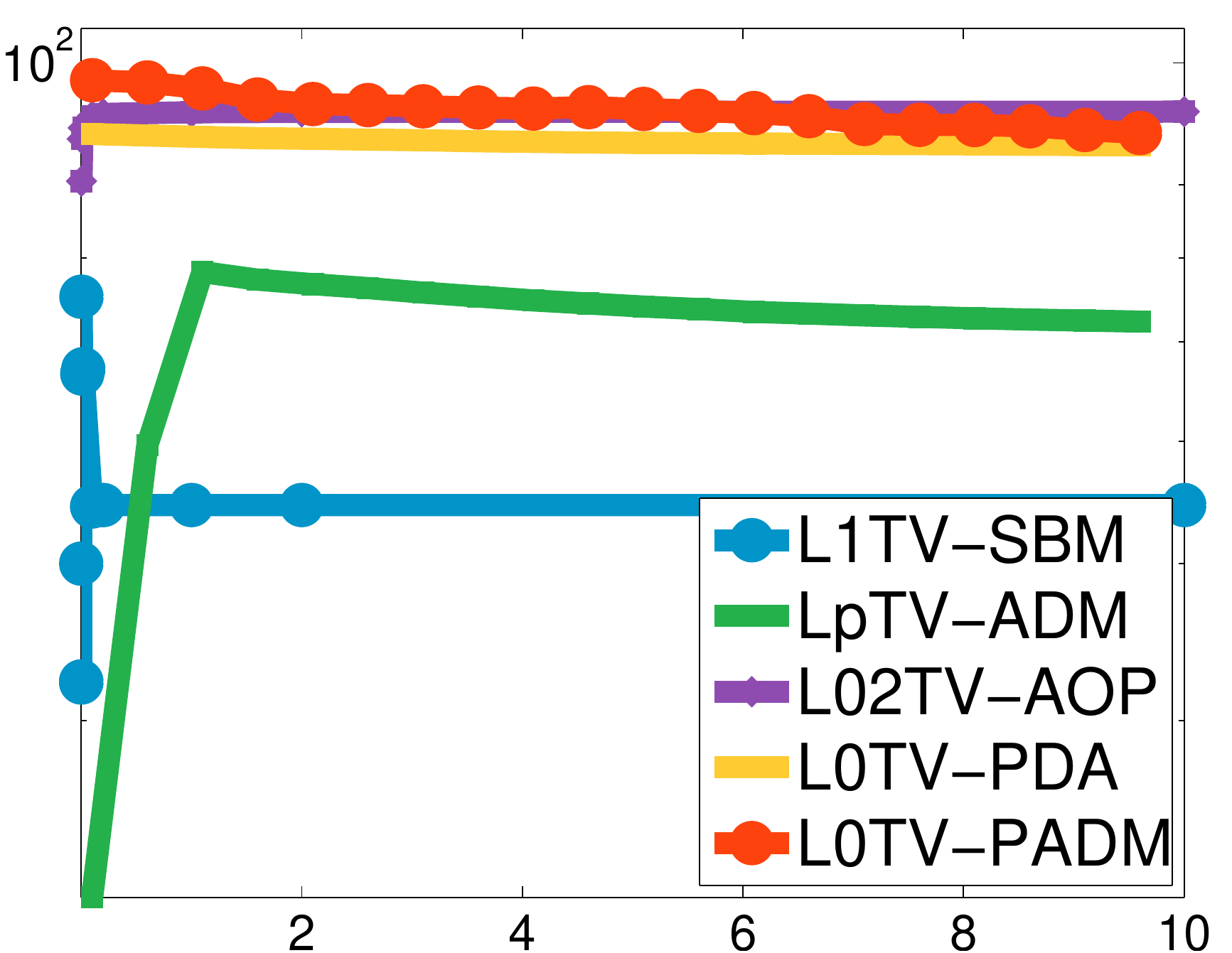}}\ghs
\subfloat[\figsizetwo Mixed]{\includegraphics[width=\figurewidth, height=\figureheight]{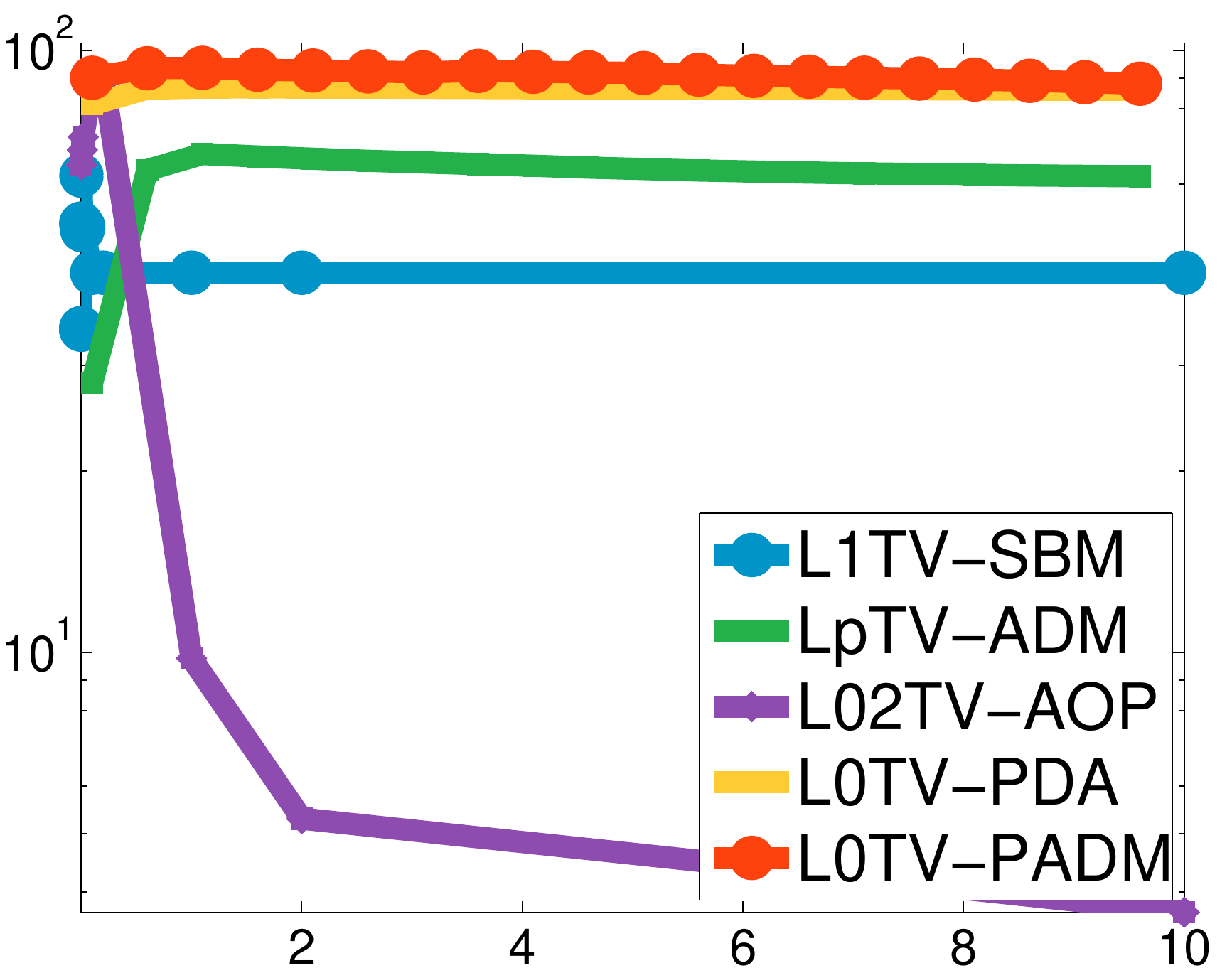}}
\\
\vspace{-8pt}
\subfloat[\figsizetwo Random-Value]{\includegraphics[width=\figurewidth, height=\figureheight]{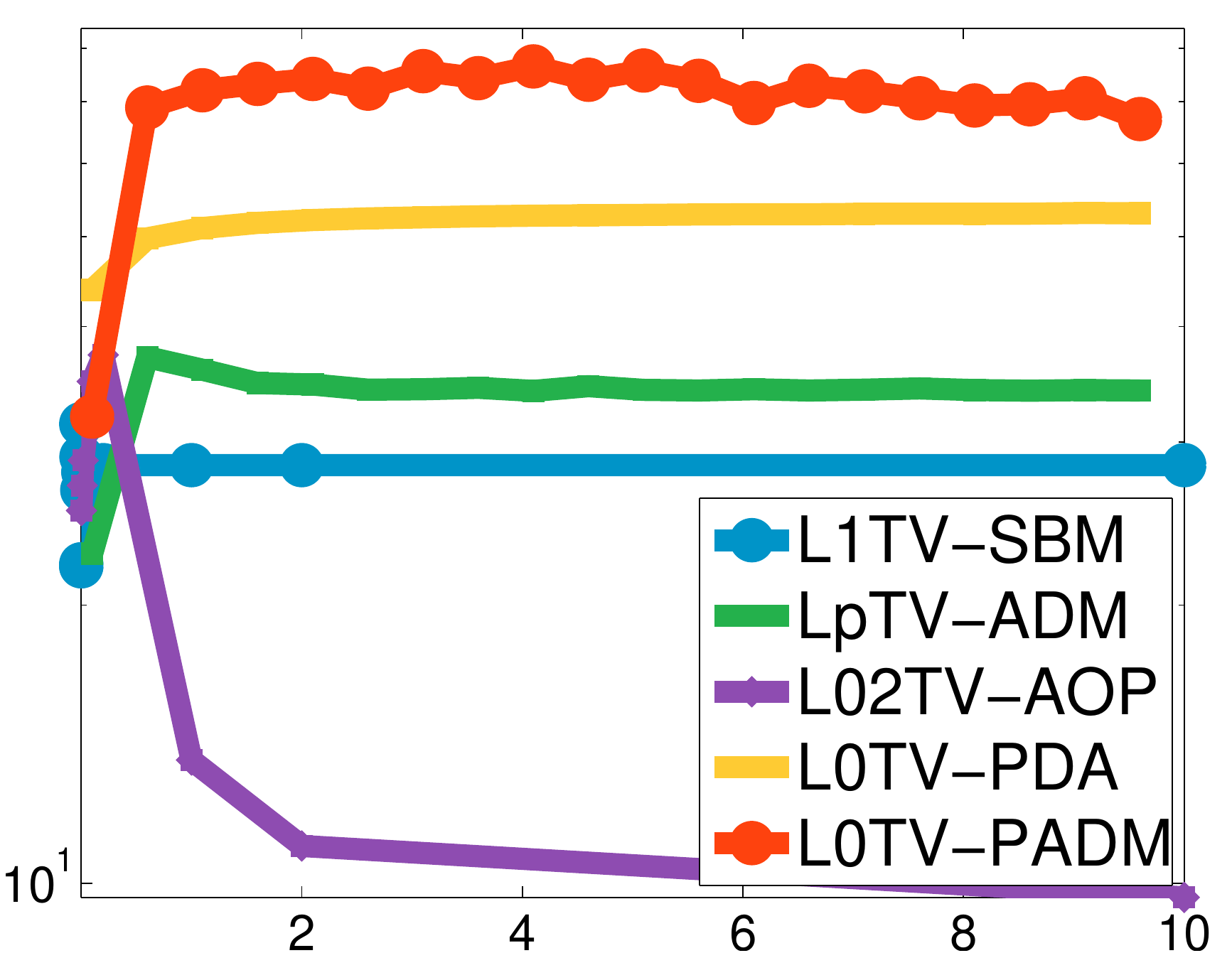}}\ghs
\subfloat[\figsizetwo Salt-and-Pepper]{\includegraphics[width=\figurewidth, height=\figureheight]{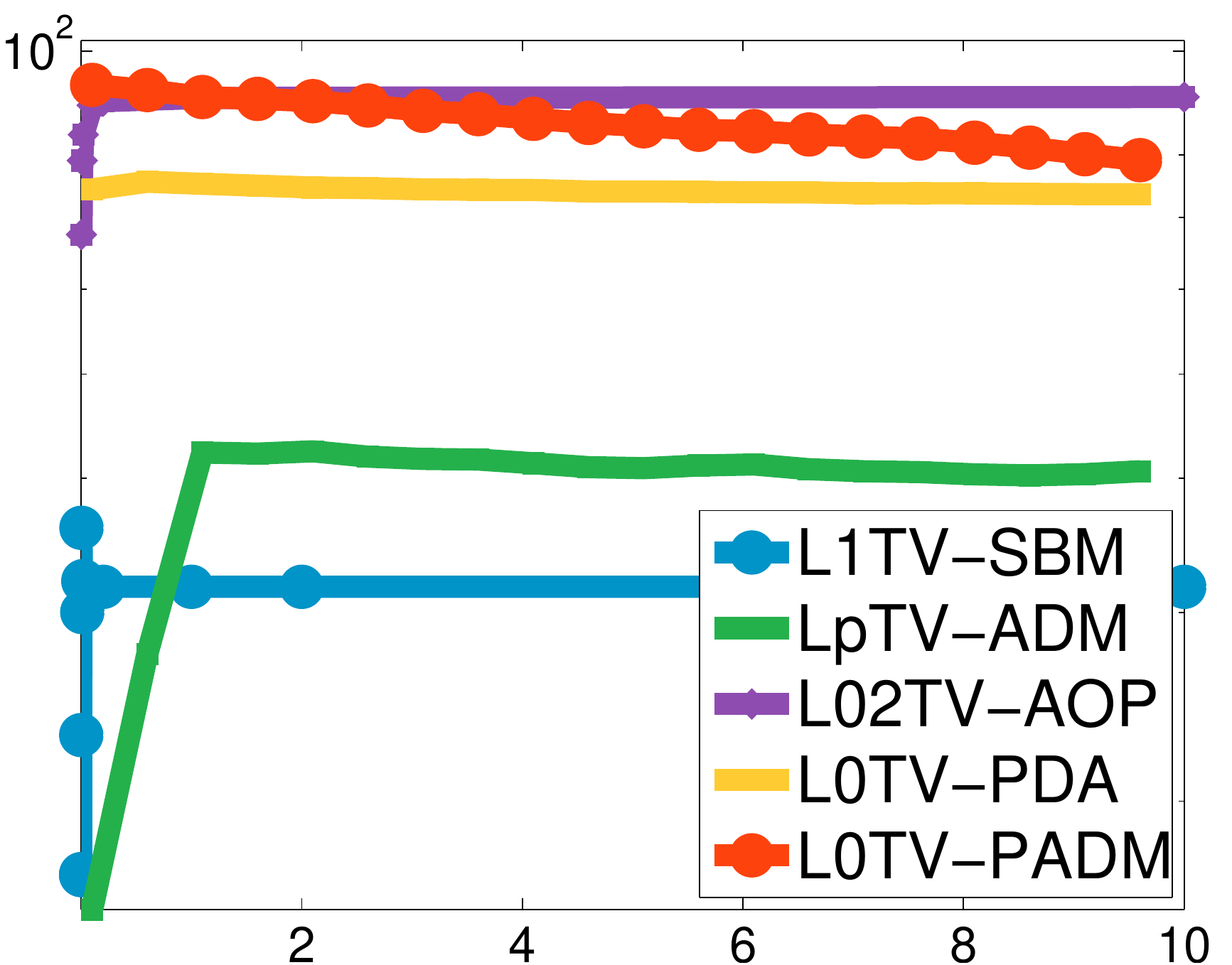}}\ghs
\subfloat[\figsizetwo Mixed]{\includegraphics[width=\figurewidth, height=\figureheight]{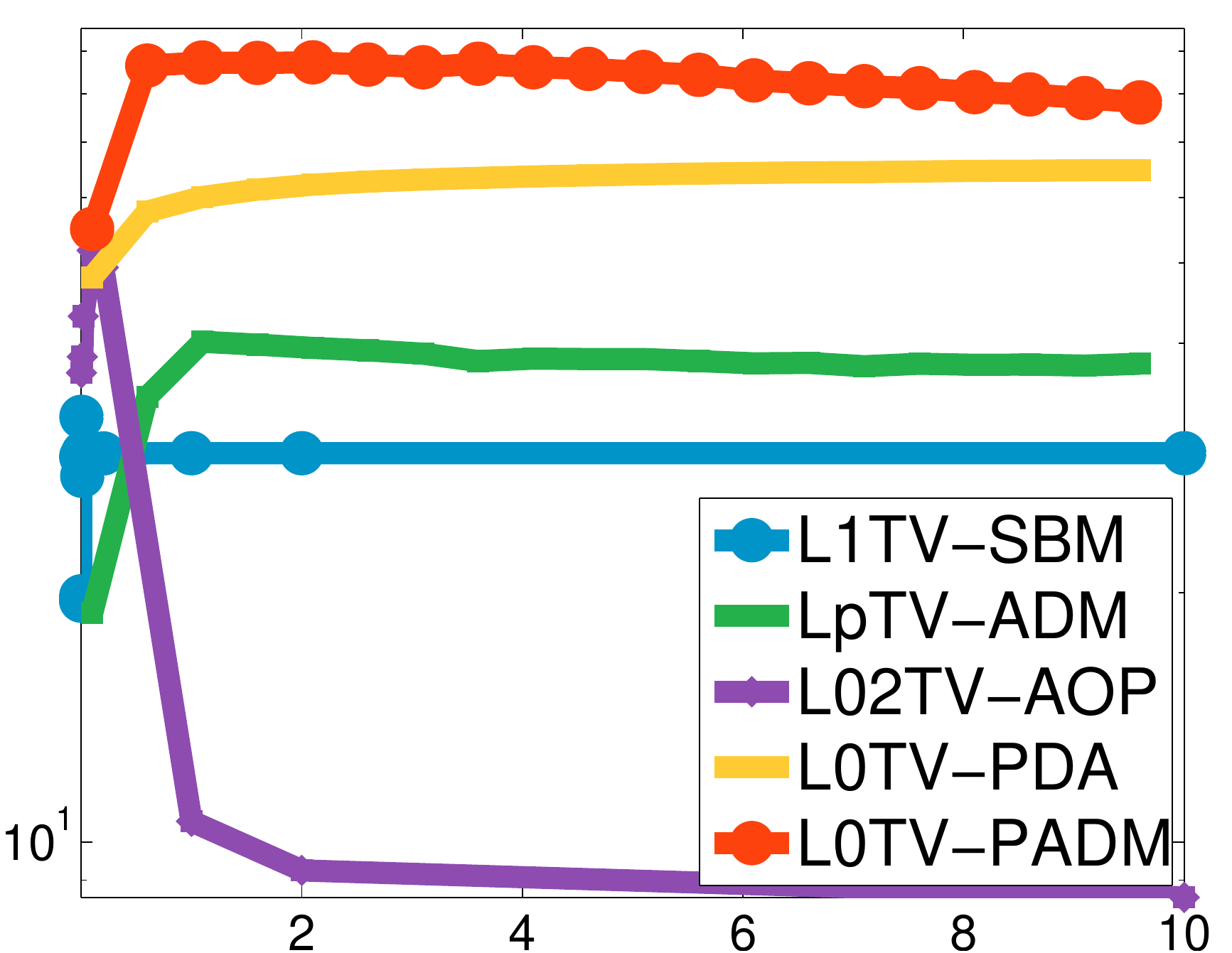}}

\caption{Image deblurring with varying the tuning parameter $\lambda$ in (\ref{eq:l0tv:2}) on ` cameraman' image. First row: noise level = 50\%. Second row: noise level = 70\%. Third row: noise level = 90\%.}\label{fig:deblurring:lambda}

\subfloat[\figsizetwo Random-Value]{\includegraphics[width=\figurewidth, height=\figureheight]{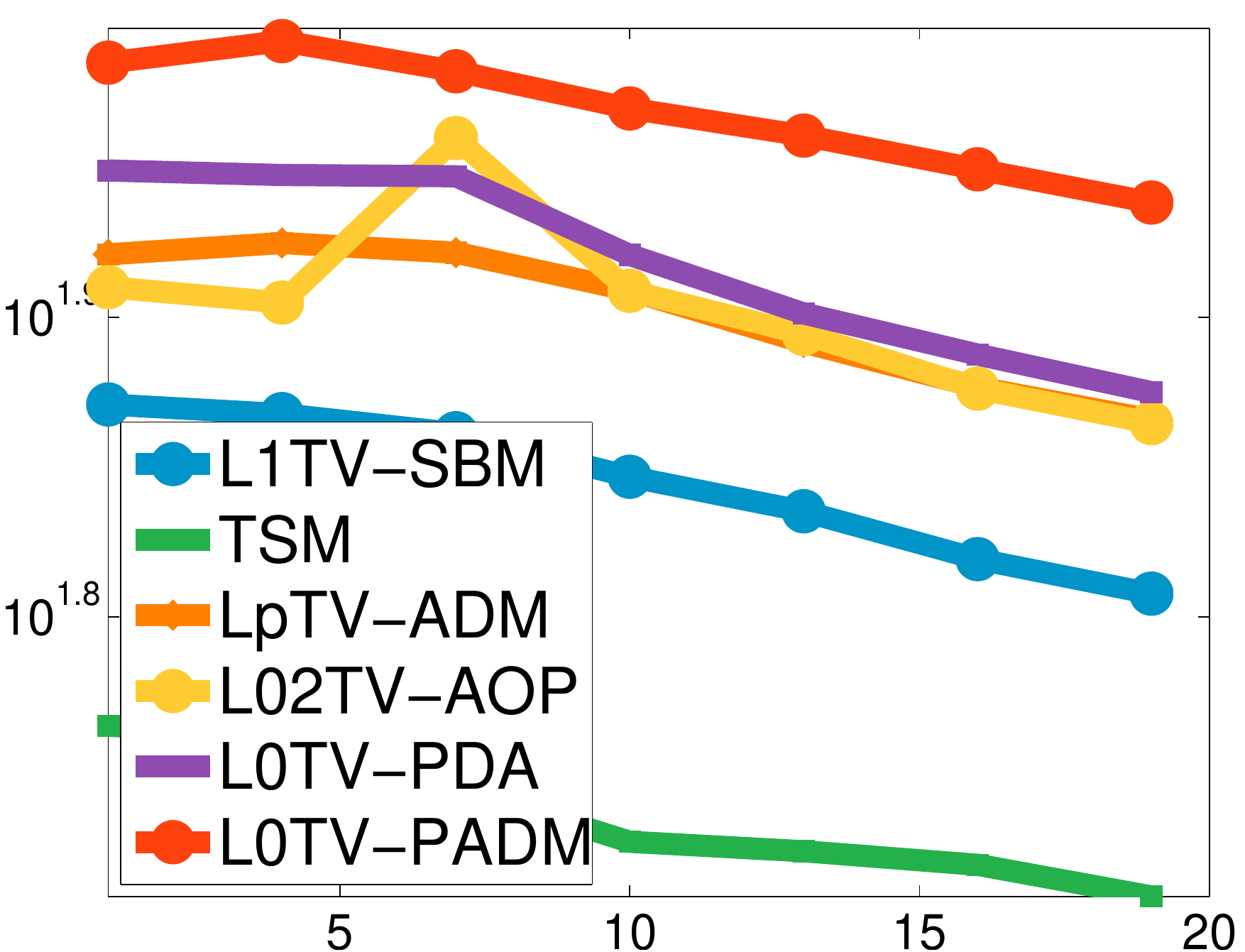}}\ghs
\subfloat[\figsizetwo Salt-and-Pepper]{\includegraphics[width=\figurewidth, height=\figureheight]{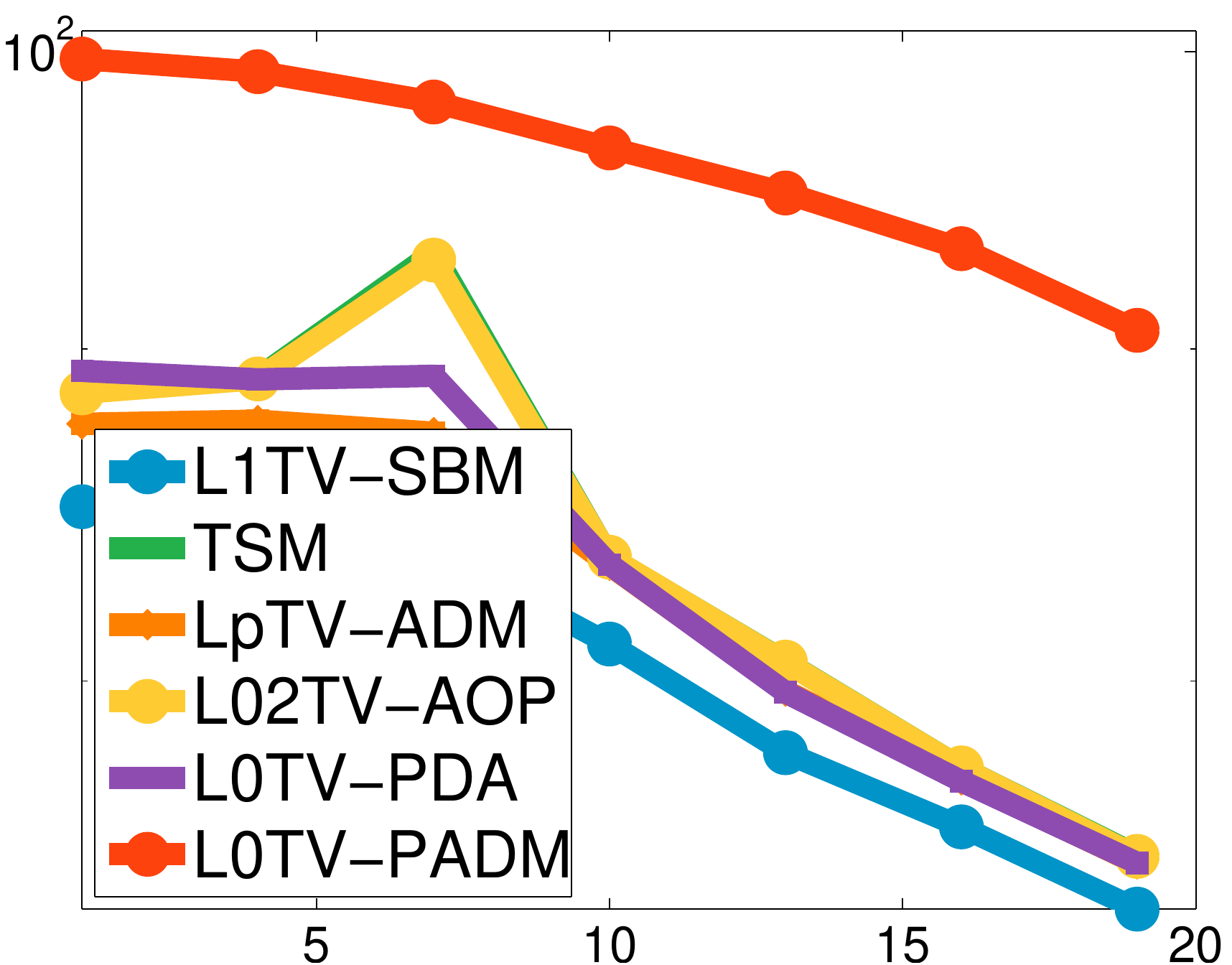}}\ghs
\subfloat[\figsizetwo Mixed]{\includegraphics[width=\figurewidth, height=\figureheight]{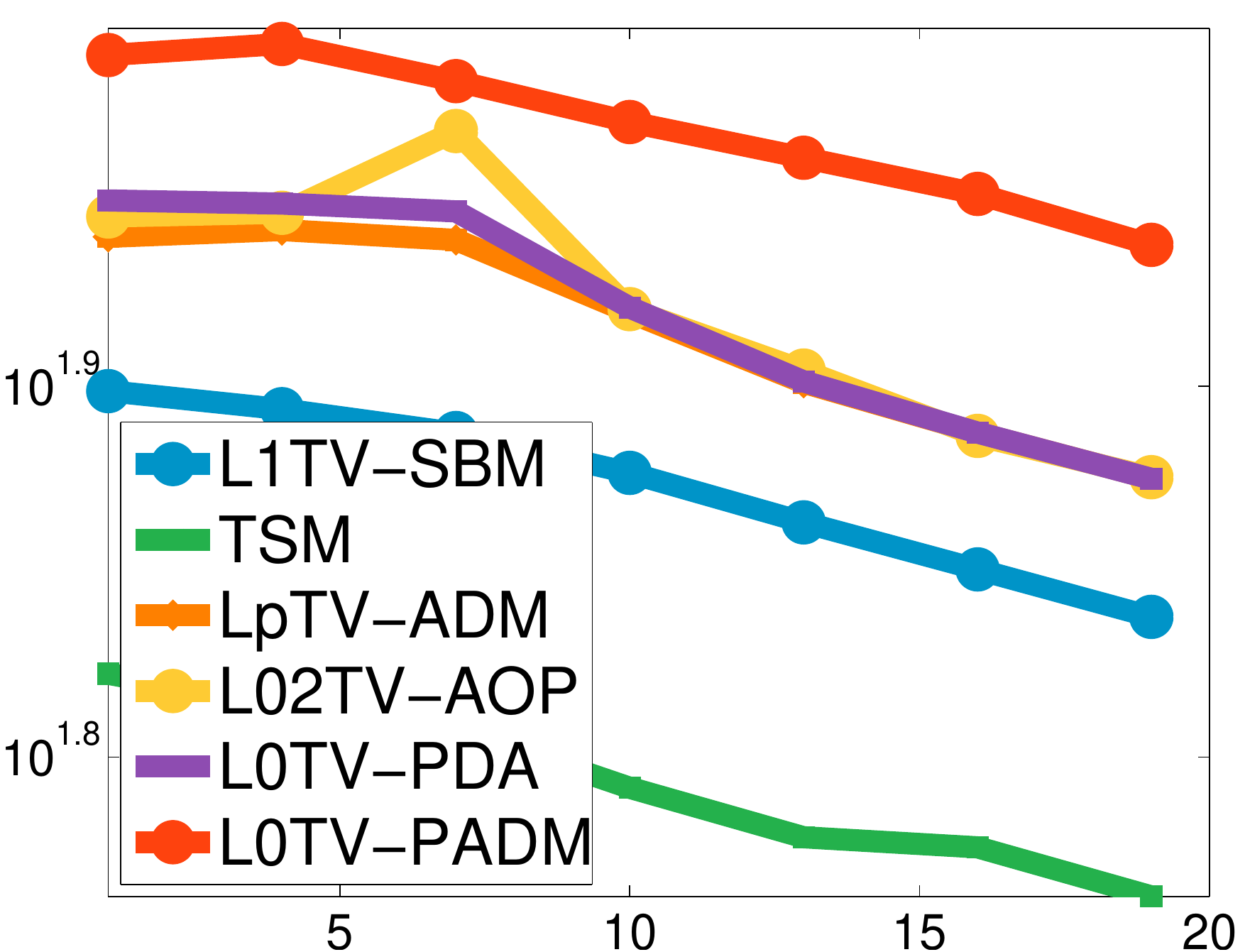}}
\\
\vspace{-8pt}
\subfloat[\figsizetwo Random-Value]{\includegraphics[width=\figurewidth, height=\figureheight]{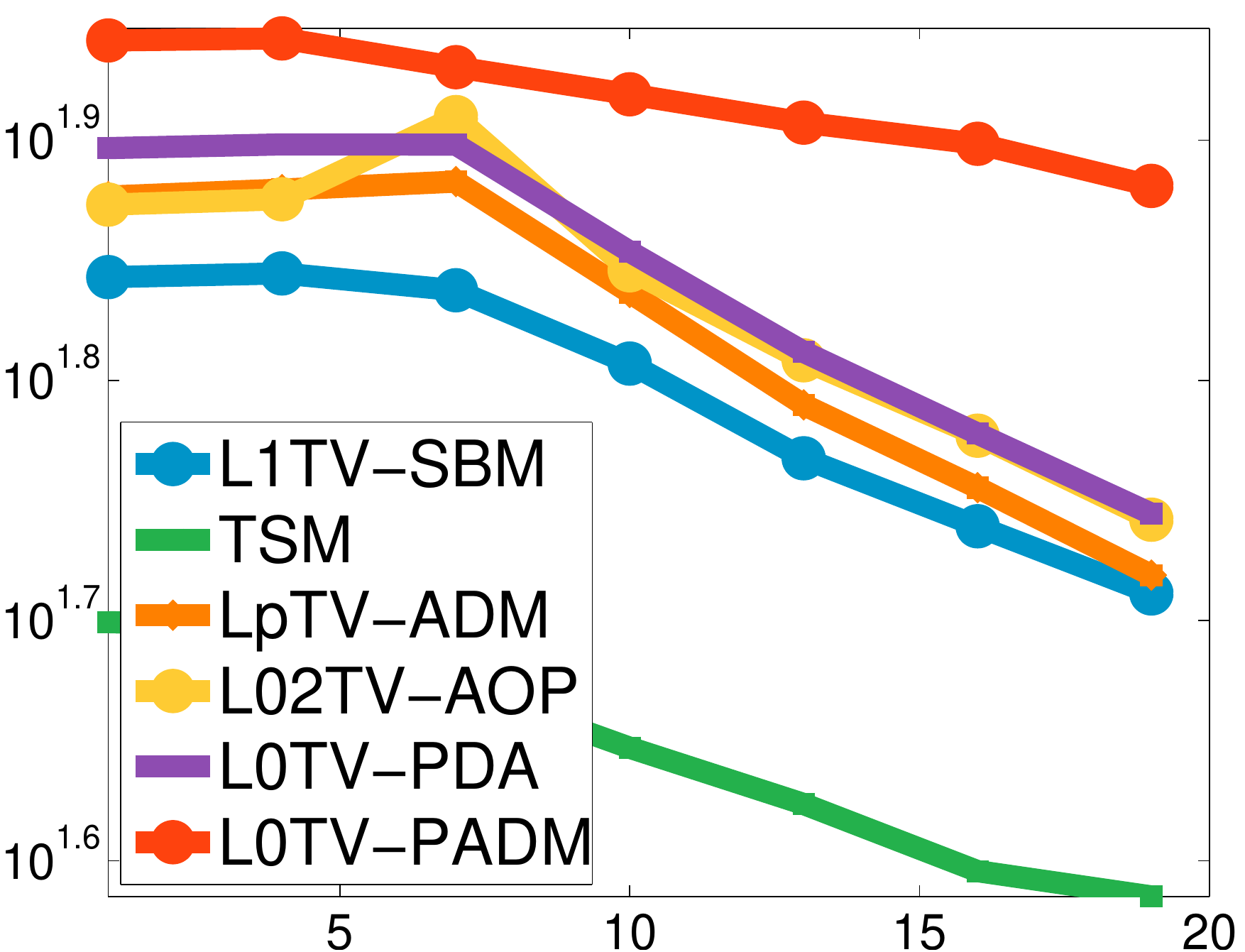}}\ghs
\subfloat[\figsizetwo Salt-and-Pepper]{\includegraphics[width=\figurewidth, height=\figureheight]{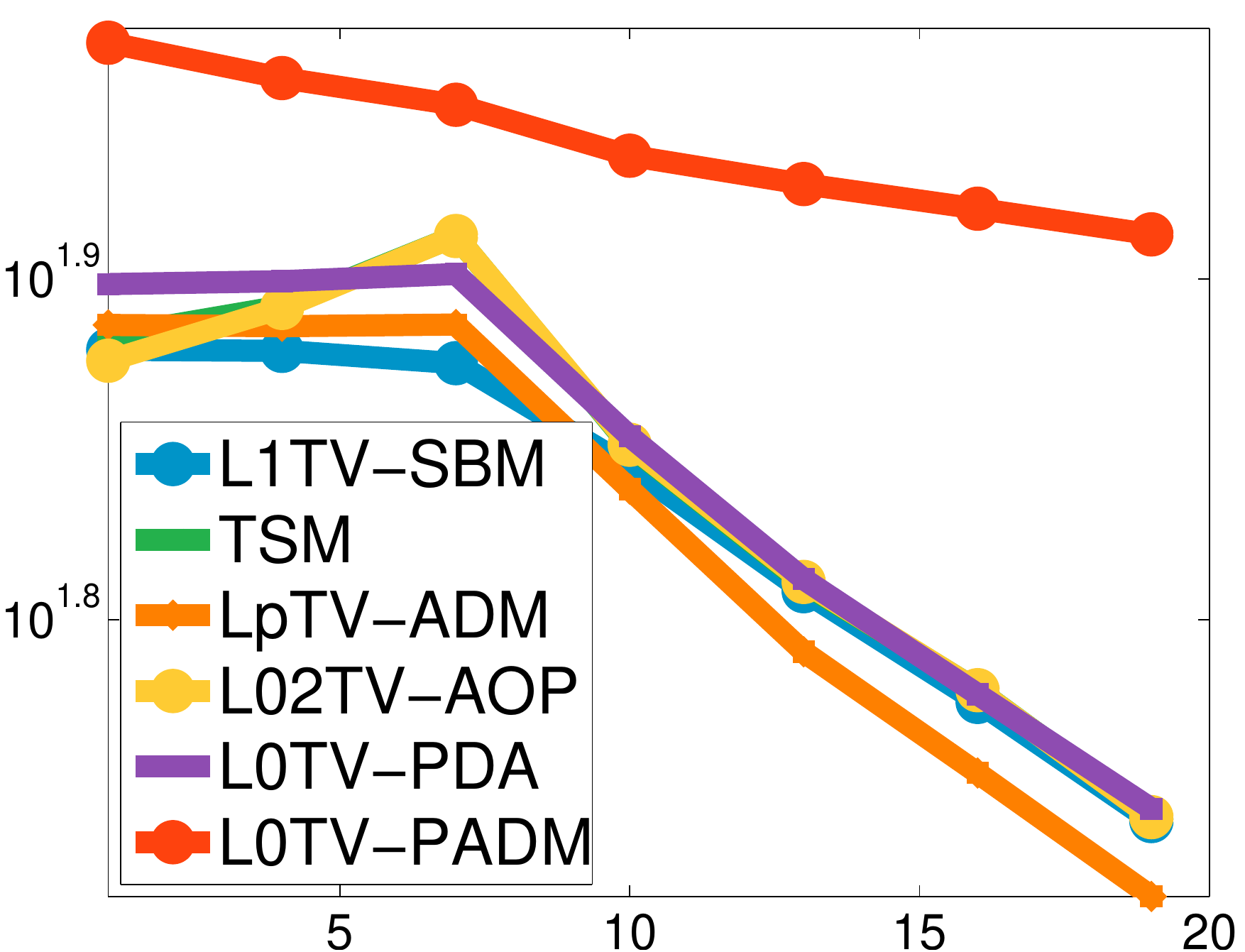}}\ghs
\subfloat[\figsizetwo Mixed]{\includegraphics[width=\figurewidth, height=\figureheight]{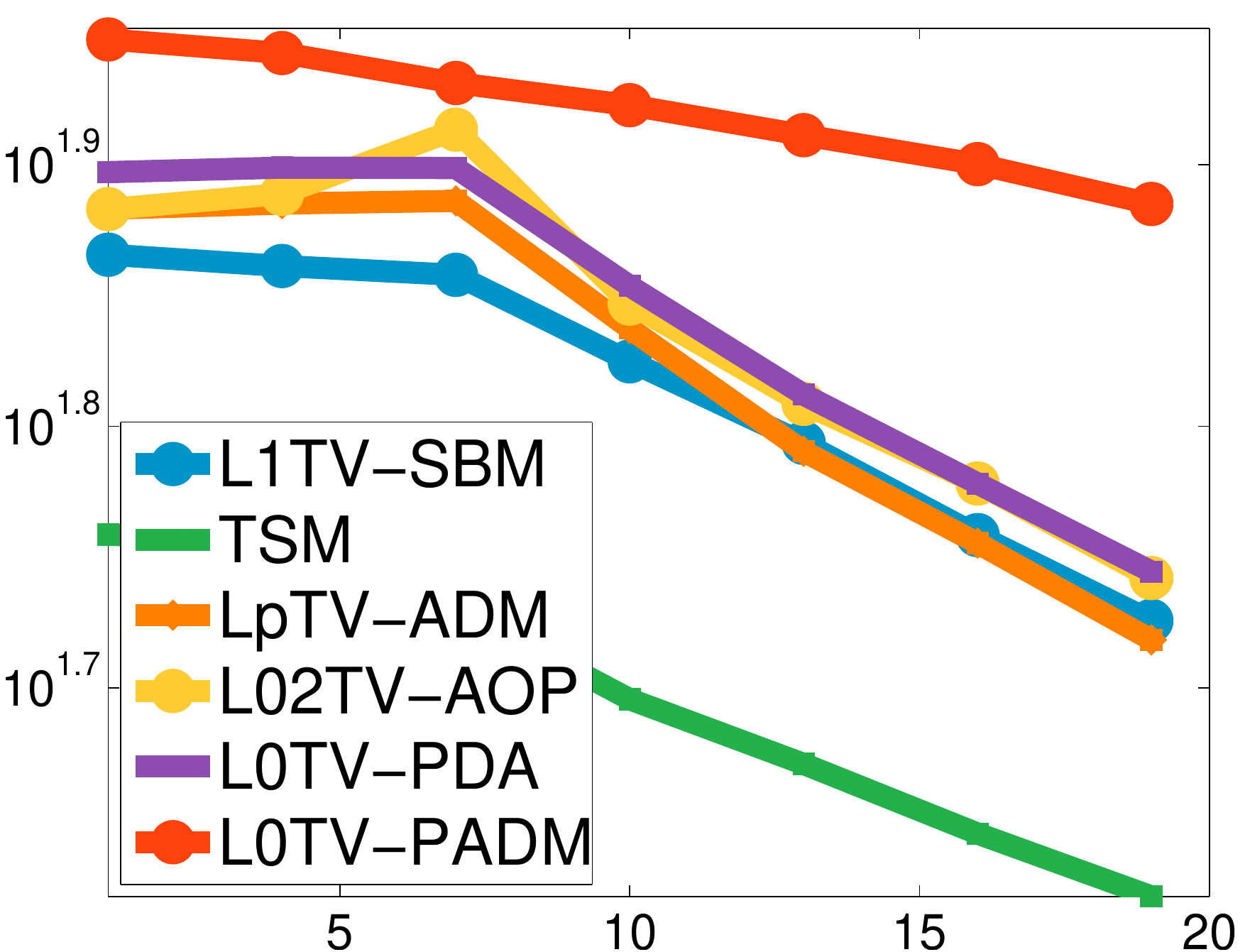}}

\caption{Image deblurring with varying the radius parameter $r$ in (\ref{eq:varyingr}). First row: `cameraman' image. Second row: `barbara' image.}
\label{fig:deblurring:R}

\end{center}
\end{figure}

\begin{table*} [!t]
\tiny
\begin{center}
\caption{General deblurring problems. The results separated by `/' are $SNR_{0}$, $SNR_{1}$ and $SNR_{2}$, respectively. The $1^{st}$, $2^{nd}$, and $3^{rd}$ best results are colored with \textbf{\cone{red}}, \ctwo{blue} and \cthree{green}, respectively. }
\label{deblurring:rv}
\scalebox{0.98}{\begin{tabular}{|p{1.8cm}|p{1.5cm}|p{1.5cm}|p{1.5cm}|p{1.5cm}|p{1.5cm}|p{1.5cm}|p{1.5cm}|}
\hline
\diagbox{Img.}{Alg.} & Corrupted &  $\ell_1TV$-$SBM$& $TSM$           &{ \tablefont $\ell_{p}TV$-$PADMM$}    & $\ell_{02}TV$-$AOP$             & $\ell_0TV$-$PDA$      & { \tablefont$\ell_0TV$-$PADM$}                 \\
\hline\hline
\multicolumn{8}{|c|}{Random-Valued Impulse Noise} \\
\hline
walkbridge+10\%  &  63/2.9/3.4  &  74/4.8/8.6  &  72/4.6/8.2  &  \cthree{77}/\cthree{5.1}/\cthree{9.2}  &  \ctwo{81}/\ctwo{5.6}/\ctwo{10.1}  &  76/5.0/9.0  &  \textbf{\cone{91}}/\textbf{\cone{7.0}}/\textbf{\cone{13.2}}  \\
walkbridge+30\%  &  52/1.1/0.0  &  72/4.6/8.1  &  61/3.7/6.8  &  \cthree{75}/\cthree{4.9}/\cthree{8.7}  &  \ctwo{79}/\ctwo{5.4}/\ctwo{9.7}  &  74/4.8/\cthree{8.7}  &  \textbf{\cone{86}}/\textbf{\cone{6.4}}/\textbf{\cone{11.7}}  \\
walkbridge+50\%  &  42/-0.2/-1.9  &  63/3.8/6.9  &  46/2.4/4.6  &  71/4.5/8.0  &  \ctwo{75}/\ctwo{4.9}/\ctwo{8.6}  &  \cthree{73}/\cthree{4.7}/\cthree{8.3}  &  \textbf{\cone{84}}/\textbf{\cone{6.0}}/\textbf{\cone{11.0}}  \\
walkbridge+70\%  &  31/-1.2/-3.2  &  46/2.1/3.8  &  33/1.1/2.3  &  55/2.9/\cthree{5.1}  &  \cthree{65}/\cthree{3.3}/4.8  &  \ctwo{69}/\ctwo{4.3}/\ctwo{7.7}  &  \textbf{\cone{81}}/\textbf{\cone{5.6}}/\textbf{\cone{10.1}}  \\
walkbridge+90\%  &  21/-2.0/-4.2  &  28/0.3/0.8  &  25/0.2/0.5  &  31/\cthree{0.6}/\cthree{1.2}  &  \cthree{33}/0.4/0.6  &  \ctwo{42}/\ctwo{1.7}/\ctwo{3.0}  &  \textbf{\cone{67}}/\textbf{\cone{3.7}}/\textbf{\cone{5.8}}  \\
\hline
pepper+10\%  &  81/4.9/4.5  &  \cthree{94}/\cthree{9.3}/\cthree{14.7}  &  93/8.3/13.6  &  70/5.1/10.1  &  \ctwo{96}/\ctwo{9.7}/\ctwo{15.8}  &  \cthree{94}/9.0/\cthree{14.7}  &  \textbf{\cone{99}}/\textbf{\cone{11.1}}/\textbf{\cone{19.8}}  \\
pepper+30\%  &  66/2.1/0.3  &  92/8.5/13.3  &  82/5.7/9.9  &  68/4.9/9.7  &  \ctwo{96}/\ctwo{9.7}/\ctwo{15.8}  &  \cthree{93}/\cthree{8.8}/\cthree{14.1}  &  \textbf{\cone{98}}/\textbf{\cone{10.7}}/\textbf{\cone{18.8}}  \\
pepper+50\%  &  52/0.4/-1.8  &  83/6.4/9.9  &  58/3.4/6.0  &  65/4.6/8.9  &  \ctwo{95}/\ctwo{9.3}/\ctwo{14.9}  &  \cthree{92}/\cthree{8.5}/\cthree{13.5}  &  \textbf{\cone{98}}/\textbf{\cone{10.4}}/\textbf{\cone{17.8}}  \\
pepper+70\%  &  37/-0.8/-3.2  &  58/3.1/4.7  &  37/1.6/2.9  &  52/3.0/\cthree{5.4}  &  \cthree{82}/\cthree{5.1}/\cthree{5.4}  &  \ctwo{90}/\ctwo{7.8}/\ctwo{12.1}  &  \textbf{\cone{97}}/\textbf{\cone{9.8}}/\textbf{\cone{16.4}}  \\
pepper+90\%  &  23/-1.8/-4.3  &  29/0.6/1.0  &  24/0.4/0.7  &  29/\cthree{0.9}/\cthree{1.3}  &  \cthree{38}/\cthree{0.9}/0.7  &  \ctwo{54}/\ctwo{2.5}/\ctwo{3.5}  &  \textbf{\cone{85}}/\textbf{\cone{6.1}}/\textbf{\cone{7.2}}  \\
\hline
mandrill+10\%  &  59/1.6/1.3  &  \cthree{67}/2.9/\cthree{4.7}  &  65/2.7/4.3  &  54/2.1/3.8  &  \ctwo{68}/\cthree{3.0}/4.5  &  \ctwo{68}/\ctwo{3.1}/\ctwo{5.0}  &  \textbf{\cone{78}}/\textbf{\cone{4.3}}/\textbf{\cone{7.3}}  \\
mandrill+30\%  &  50/0.0/-1.7  &  66/\cthree{2.9}/\cthree{4.6}  &  60/2.3/3.9  &  52/2.1/3.7  &  \ctwo{68}/\ctwo{3.0}/\cthree{4.6}  &  \cthree{67}/\ctwo{3.0}/\ctwo{4.8}  &  \textbf{\cone{76}}/\textbf{\cone{4.0}}/\textbf{\cone{6.8}}  \\
mandrill+50\%  &  40/-1.1/-3.4  &  64/\cthree{2.7}/4.3  &  50/1.6/2.9  &  51/2.0/3.5  &  \ctwo{68}/\ctwo{2.9}/\cthree{4.5}  &  \cthree{66}/\ctwo{2.9}/\ctwo{4.6}  &  \textbf{\cone{73}}/\textbf{\cone{3.6}}/\textbf{\cone{6.0}}  \\
mandrill+70\%  &  30/-2.0/-4.7  &  53/1.8/3.1  &  40/0.9/1.7  &  46/1.6/2.9  &  \cthree{64}/\cthree{2.5}/\cthree{3.6}  &  \ctwo{65}/\ctwo{2.7}/\ctwo{4.4}  &  \textbf{\cone{70}}/\textbf{\cone{3.3}}/\textbf{\cone{5.4}}  \\
mandrill+90\%  &  21/-2.7/-5.6  &  38/0.5/\cthree{0.9}  &  36/0.3/0.6  &  34/0.4/0.7  &  \cthree{42}/\cthree{0.6}/0.8  &  \ctwo{49}/\ctwo{1.5}/\ctwo{2.5}  &  \textbf{\cone{65}}/\textbf{\cone{2.7}}/\textbf{\cone{4.2}}  \\
\hline
lake+10\%  &  71/4.8/4.9  &  \cthree{84}/7.6/11.6  &  83/7.3/11.3  &  83/6.7/11.3  &  \ctwo{89}/\ctwo{8.6}/\ctwo{13.8}  &  \cthree{84}/\cthree{7.7}/\cthree{12.1}  &  \textbf{\cone{96}}/\textbf{\cone{10.0}}/\textbf{\cone{17.4}}  \\
lake+30\%  &  59/2.6/1.2  &  81/7.1/10.8  &  65/5.2/8.9  &  80/6.4/10.7  &  \ctwo{89}/\ctwo{8.5}/\ctwo{13.2}  &  \cthree{83}/\cthree{7.4}/\cthree{11.6}  &  \textbf{\cone{94}}/\textbf{\cone{9.5}}/\textbf{\cone{15.9}}  \\
lake+50\%  &  46/1.1/-0.7  &  68/5.5/8.8  &  35/3.2/5.6  &  76/6.0/9.8  &  \ctwo{86}/\ctwo{7.9}/\ctwo{11.9}  &  \cthree{82}/\cthree{7.2}/\cthree{11.1}  &  \textbf{\cone{92}}/\textbf{\cone{9.1}}/\textbf{\cone{15.1}}  \\
lake+70\%  &  34/0.0/-2.1  &  35/2.6/4.5  &  22/1.6/2.9  &  39/3.3/\cthree{5.6}  &  \cthree{66}/\cthree{4.3}/5.4  &  \ctwo{79}/\ctwo{6.7}/\ctwo{10.2}  &  \textbf{\cone{89}}/\textbf{\cone{8.5}}/\textbf{\cone{13.8}}  \\
lake+90\%  &  \cthree{22}/-0.9/-3.1  &  \cthree{22}/0.6/1.0  &  16/0.4/0.8  &  \cthree{22}/\cthree{0.7}/\cthree{1.3}  &  21/0.6/0.8  &  \ctwo{31}/\ctwo{2.1}/\ctwo{3.5}  &  \textbf{\cone{74}}/\textbf{\cone{5.6}}/\textbf{\cone{7.2}}  \\
\hline
jetplane+10\%  &  76/3.3/2.1  &  88/6.7/9.9  &  88/6.1/9.7  &  63/2.8/6.5  &  \ctwo{93}/\ctwo{7.9}/\ctwo{12.5}  &  \cthree{89}/\cthree{6.8}/\cthree{10.5}  &  \textbf{\cone{98}}/\textbf{\cone{9.1}}/\textbf{\cone{16.6}}  \\
jetplane+30\%  &  63/0.7/-1.9  &  86/6.2/9.1  &  68/3.2/6.3  &  66/2.7/6.2  &  \ctwo{93}/\ctwo{7.8}/\ctwo{12.0}  &  \cthree{88}/\cthree{6.6}/\cthree{10.0}  &  \textbf{\cone{97}}/\textbf{\cone{8.8}}/\textbf{\cone{15.6}}  \\
jetplane+50\%  &  49/-0.9/-3.9  &  74/3.9/6.6  &  34/0.9/2.6  &  55/2.5/5.6  &  \ctwo{91}/\ctwo{7.0}/\ctwo{9.7}  &  \cthree{87}/\cthree{6.3}/\cthree{9.4}  &  \textbf{\cone{95}}/\textbf{\cone{8.4}}/\textbf{\cone{14.2}}  \\
jetplane+70\%  &  36/-2.1/-5.3  &  37/0.3/1.3  &  22/-0.7/-0.3  &  35/-0.1/0.6  &  \cthree{64}/\cthree{1.5}/\cthree{1.9}  &  \ctwo{84}/\ctwo{5.8}/\ctwo{8.5}  &  \textbf{\cone{93}}/\textbf{\cone{7.8}}/\textbf{\cone{12.4}}  \\
jetplane+90\%  &  \cthree{23}/-3.0/-6.3  &  \cthree{23}/\cthree{-1.7}/\cthree{-2.3}  &  14/-1.9/-2.5  &  16/-2.2/-3.3  &  20/\cthree{-1.7}/-2.5  &  \ctwo{30}/\ctwo{0.0}/\ctwo{0.6}  &  \textbf{\cone{80}}/\textbf{\cone{4.5}}/\textbf{\cone{5.1}}  \\
\hline\hline
\multicolumn{8}{|c|}{Salt-and-Pepper Impulse Noise} \\
\hline
walkbridge+10\%  &  61/2.0/0.8  &  73/4.8/8.5  &  \ctwo{80}/\ctwo{5.6}/\ctwo{10.1}  &  \cthree{76}/\cthree{5.1}/\cthree{9.1}  &  \ctwo{80}/\ctwo{5.6}/\ctwo{10.1}  &  \cthree{76}/5.0/9.0  &  \textbf{\cone{94}}/\textbf{\cone{7.4}}/\textbf{\cone{14.3}}  \\
walkbridge+30\%  &  48/-0.5/-3.2  &  71/4.5/7.9  &  \ctwo{79}/\ctwo{5.4}/\ctwo{9.7}  &  74/4.8/8.5  &  \ctwo{79}/\ctwo{5.4}/\ctwo{9.7}  &  \cthree{75}/\cthree{4.9}/\cthree{8.8}  &  \textbf{\cone{92}}/\textbf{\cone{7.2}}/\textbf{\cone{13.7}}  \\
walkbridge+50\%  &  35/-2.1/-5.3  &  67/4.1/7.3  &  \ctwo{77}/\ctwo{5.2}/\ctwo{9.3}  &  72/4.5/8.1  &  \ctwo{77}/\ctwo{5.2}/\ctwo{9.3}  &  \cthree{73}/\cthree{4.8}/\cthree{8.5}  &  \textbf{\cone{90}}/\textbf{\cone{6.8}}/\textbf{\cone{12.9}}  \\
walkbridge+70\%  &  22/-3.3/-6.7  &  53/2.8/5.2  &  \ctwo{75}/\ctwo{5.0}/\ctwo{8.8}  &  61/3.5/6.4  &  \ctwo{75}/\cthree{4.9}/\ctwo{8.8}  &  \cthree{71}/4.5/\cthree{8.1}  &  \textbf{\cone{86}}/\textbf{\cone{6.4}}/\textbf{\cone{11.8}}  \\
walkbridge+90\%  &  8/-4.2/-7.7  &  31/0.6/1.0  &  \ctwo{73}/\ctwo{4.7}/\ctwo{8.3}  &  34/0.9/1.7  &  \ctwo{73}/\ctwo{4.7}/\ctwo{8.3}  &  \cthree{59}/\cthree{3.4}/\cthree{6.3}  &  \textbf{\cone{79}}/\textbf{\cone{5.4}}/\textbf{\cone{9.9}}  \\
\hline
pepper+10\%  &  79/3.6/1.3  &  \cthree{94}/8.9/14.2  &  \ctwo{96}/\ctwo{9.7}/\ctwo{15.8}  &  69/5.0/10.0  &  \ctwo{96}/\cthree{9.6}/\ctwo{15.8}  &  \cthree{94}/9.1/\cthree{14.8}  &  \textbf{\cone{99}}/\textbf{\cone{11.4}}/\textbf{\cone{20.3}}  \\
pepper+30\%  &  62/0.2/-3.2  &  92/8.5/13.2  &  \ctwo{96}/\ctwo{9.6}/\ctwo{15.7}  &  69/4.9/9.6  &  \ctwo{96}/\ctwo{9.6}/\ctwo{15.7}  &  \cthree{94}/\cthree{8.9}/\cthree{14.4}  &  \textbf{\cone{99}}/\textbf{\cone{11.2}}/\textbf{\cone{19.7}}  \\
pepper+50\%  &  45/-1.7/-5.4  &  87/7.3/11.2  &  \ctwo{95}/\ctwo{9.4}/\ctwo{15.4}  &  66/4.7/9.1  &  \ctwo{95}/\ctwo{9.4}/\ctwo{15.4}  &  \cthree{93}/\cthree{8.6}/\cthree{13.8}  &  \textbf{\cone{99}}/\textbf{\cone{10.9}}/\textbf{\cone{19.1}}  \\
pepper+70\%  &  28/-3.0/-6.8  &  70/4.3/6.5  &  \ctwo{95}/\ctwo{9.2}/\cthree{14.8}  &  56/3.7/6.8  &  \ctwo{95}/\ctwo{9.2}/\ctwo{14.9}  &  \cthree{91}/\cthree{8.3}/13.0  &  \textbf{\cone{98}}/\textbf{\cone{10.3}}/\textbf{\cone{18.2}}  \\
pepper+90\%  &  11/-4.1/-7.9  &  33/0.8/1.1  &  \ctwo{94}/\ctwo{8.8}/\ctwo{14.1}  &  32/1.1/1.8  &  \ctwo{94}/\ctwo{8.8}/\ctwo{14.1}  &  \cthree{79}/\cthree{5.6}/\cthree{8.8}  &  \textbf{\cone{96}}/\textbf{\cone{9.5}}/\textbf{\cone{15.8}}  \\
\hline
mandrill+10\%  &  58/0.7/-1.3  &  \cthree{67}/\cthree{2.9}/\cthree{4.7}  &  \cthree{67}/\cthree{2.9}/4.4  &  53/2.1/3.8  &  \cthree{67}/\cthree{2.9}/4.4  &  \ctwo{68}/\ctwo{3.1}/\ctwo{5.0}  &  \textbf{\cone{86}}/\textbf{\cone{5.2}}/\textbf{\cone{9.5}}  \\
mandrill+30\%  &  45/-1.7/-5.2  &  65/2.8/4.4  &  \cthree{67}/\cthree{2.9}/\cthree{4.5}  &  52/2.1/3.6  &  \cthree{67}/\cthree{2.9}/\cthree{4.5}  &  \ctwo{68}/\ctwo{3.0}/\ctwo{4.9}  &  \textbf{\cone{83}}/\textbf{\cone{4.9}}/\textbf{\cone{8.7}}  \\
mandrill+50\%  &  32/-3.2/-7.2  &  64/2.6/4.2  &  \cthree{66}/\cthree{2.8}/\cthree{4.4}  &  51/2.0/3.5  &  \cthree{66}/\cthree{2.8}/\cthree{4.4}  &  \ctwo{67}/\ctwo{3.0}/\ctwo{4.7}  &  \textbf{\cone{80}}/\textbf{\cone{4.5}}/\textbf{\cone{7.9}}  \\
mandrill+70\%  &  19/-4.4/-8.6  &  56/2.0/3.3  &  \cthree{65}/\cthree{2.7}/\cthree{4.2}  &  48/1.8/3.1  &  \cthree{65}/\cthree{2.7}/\cthree{4.2}  &  \ctwo{66}/\ctwo{2.8}/\ctwo{4.5}  &  \textbf{\cone{75}}/\textbf{\cone{4.0}}/\textbf{\cone{6.7}}  \\
mandrill+90\%  &  7/-5.2/-9.6  &  39/0.5/0.8  &  \ctwo{65}/\ctwo{2.7}/\ctwo{4.2}  &  35/0.5/1.0  &  \ctwo{65}/\ctwo{2.7}/\ctwo{4.2}  &  \cthree{60}/\cthree{2.4}/\cthree{3.9}  &  \textbf{\cone{70}}/\textbf{\cone{3.3}}/\textbf{\cone{5.3}}  \\
\hline
lake+10\%  &  69/3.9/2.4  &  83/7.4/11.4  &  \ctwo{90}/\ctwo{8.7}/\ctwo{13.8}  &  82/6.6/11.2  &  \ctwo{90}/\ctwo{8.7}/\ctwo{13.8}  &  \cthree{85}/\cthree{7.7}/\cthree{12.1}  &  \textbf{\cone{98}}/\textbf{\cone{10.3}}/\textbf{\cone{18.5}}  \\
lake+30\%  &  54/1.0/-1.8  &  81/7.1/10.6  &  \ctwo{89}/\ctwo{8.5}/\ctwo{13.4}  &  80/6.3/10.6  &  \ctwo{89}/\ctwo{8.5}/\ctwo{13.4}  &  \cthree{84}/\cthree{7.6}/\cthree{11.8}  &  \textbf{\cone{97}}/\textbf{\cone{10.1}}/\textbf{\cone{17.9}}  \\
lake+50\%  &  38/-0.7/-3.9  &  76/6.4/9.6  &  \ctwo{87}/\ctwo{8.2}/\ctwo{12.9}  &  77/6.0/9.8  &  \ctwo{87}/\ctwo{8.2}/\cthree{12.8}  &  \cthree{82}/\cthree{7.3}/11.3  &  \textbf{\cone{96}}/\textbf{\cone{9.8}}/\textbf{\cone{17.0}}  \\
lake+70\%  &  23/-1.9/-5.3  &  49/3.9/6.3  &  \ctwo{86}/\ctwo{7.9}/\ctwo{12.2}  &  56/4.4/7.3  &  \ctwo{86}/\ctwo{7.9}/\ctwo{12.2}  &  \cthree{81}/\cthree{7.0}/\cthree{10.7}  &  \textbf{\cone{94}}/\textbf{\cone{9.3}}/\textbf{\cone{15.9}}  \\
lake+90\%  &  8/-2.8/-6.4  &  24/0.9/1.4  &  \cthree{83}/\cthree{7.4}/\ctwo{11.2}  &  21/1.0/1.8  &  \ctwo{84}/\ctwo{7.5}/\cthree{11.1}  &  63/5.0/8.1  &  \textbf{\cone{88}}/\textbf{\cone{8.2}}/\textbf{\cone{13.3}}  \\
\hline
jetplane+10\%  &  75/2.3/-0.4  &  88/6.5/9.7  &  \ctwo{93}/\ctwo{8.0}/\ctwo{12.6}  &  67/2.8/6.5  &  \ctwo{93}/\ctwo{8.0}/\ctwo{12.6}  &  \cthree{89}/\cthree{6.9}/\cthree{10.6}  &  \textbf{\cone{99}}/\textbf{\cone{9.5}}/\textbf{\cone{17.8}}  \\
jetplane+30\%  &  58/-0.9/-4.8  &  86/6.2/9.0  &  \ctwo{93}/\ctwo{7.7}/\ctwo{11.9}  &  64/2.7/6.1  &  \cthree{92}/\cthree{7.6}/\cthree{11.8}  &  88/6.7/10.2  &  \textbf{\cone{99}}/\textbf{\cone{9.4}}/\textbf{\cone{17.2}}  \\
jetplane+50\%  &  42/-2.7/-7.0  &  82/5.4/7.8  &  \ctwo{91}/\ctwo{7.5}/\cthree{11.4}  &  54/2.5/5.7  &  \ctwo{91}/\ctwo{7.5}/\ctwo{11.5}  &  \cthree{87}/\cthree{6.5}/9.7  &  \textbf{\cone{98}}/\textbf{\cone{9.0}}/\textbf{\cone{16.2}}  \\
jetplane+70\%  &  25/-3.9/-8.4  &  48/1.9/3.8  &  \ctwo{90}/\ctwo{7.1}/\ctwo{10.7}  &  39/1.2/2.9  &  \ctwo{90}/\ctwo{7.1}/\cthree{10.6}  &  \cthree{86}/\cthree{6.1}/9.0  &  \textbf{\cone{96}}/\textbf{\cone{8.7}}/\textbf{\cone{14.9}}  \\
jetplane+90\%  &  8/-4.9/-9.5  &  24/-1.3/-1.8  &  \ctwo{89}/\ctwo{6.7}/\cthree{9.7}  &  21/-1.9/-2.8  &  \ctwo{89}/\ctwo{6.7}/\ctwo{9.9}  &  \cthree{72}/\cthree{3.6}/6.1  &  \textbf{\cone{92}}/\textbf{\cone{7.2}}/\textbf{\cone{11.8}}  \\
\hline\hline
\multicolumn{8}{|c|}{Mixed Impulse Noise (Half Random-Value Noise and Half Salt-and-Pepper Noise)} \\
\hline
walkbridge+10\%  &  62/2.4/1.9  &  74/4.8/8.5  &  72/4.6/8.2  &  \cthree{77}/\cthree{5.1}/\cthree{9.2}  &  \ctwo{81}/\ctwo{5.6}/\ctwo{10.1}  &  76/5.0/9.0  &  \textbf{\cone{93}}/\textbf{\cone{7.4}}/\textbf{\cone{14.0}}  \\
walkbridge+30\%  &  50/0.2/-1.9  &  71/4.5/7.9  &  65/3.9/7.2  &  \cthree{74}/\cthree{4.8}/8.6  &  \ctwo{79}/\ctwo{5.4}/\ctwo{9.6}  &  \cthree{74}/\cthree{4.8}/\cthree{8.7}  &  \textbf{\cone{87}}/\textbf{\cone{6.5}}/\textbf{\cone{12.0}}  \\
walkbridge+50\%  &  38/-1.3/-3.9  &  64/3.8/6.9  &  52/2.9/5.5  &  71/4.5/8.0  &  \ctwo{78}/\ctwo{5.2}/\ctwo{8.8}  &  \cthree{73}/\cthree{4.7}/\cthree{8.3}  &  \textbf{\cone{84}}/\textbf{\cone{6.1}}/\textbf{\cone{11.0}}  \\
walkbridge+70\%  &  27/-2.3/-5.3  &  48/2.4/4.4  &  38/1.6/3.2  &  59/3.3/6.0  &  \ctwo{74}/\ctwo{4.5}/\cthree{7.1}  &  \cthree{70}/\cthree{4.4}/\ctwo{7.8}  &  \textbf{\cone{81}}/\textbf{\cone{5.6}}/\textbf{\cone{10.1}}  \\
walkbridge+90\%  &  15/-3.2/-6.3  &  29/0.5/1.0  &  27/0.4/0.9  &  33/0.8/\cthree{1.6}  &  \cthree{43}/\cthree{1.1}/1.4  &  \ctwo{50}/\ctwo{2.3}/\ctwo{3.7}  &  \textbf{\cone{71}}/\textbf{\cone{4.3}}/\textbf{\cone{7.2}}  \\
\hline
pepper+10\%  &  80/4.2/2.6  &  \cthree{94}/\cthree{9.1}/14.5  &  93/8.5/13.7  &  69/5.1/10.0  &  \ctwo{96}/\ctwo{9.7}/\ctwo{15.9}  &  \cthree{94}/9.0/\cthree{14.7}  &  \textbf{\cone{99}}/\textbf{\cone{11.1}}/\textbf{\cone{19.8}}  \\
pepper+30\%  &  64/1.0/-1.8  &  91/8.4/13.0  &  87/6.4/10.9  &  68/4.9/9.6  &  \ctwo{96}/\ctwo{9.7}/\ctwo{15.8}  &  \cthree{93}/\cthree{8.8}/\cthree{14.1}  &  \textbf{\cone{99}}/\textbf{\cone{10.9}}/\textbf{\cone{19.3}}  \\
pepper+50\%  &  49/-0.8/-3.9  &  84/6.7/10.2  &  68/4.2/7.5  &  66/4.7/9.1  &  \ctwo{96}/\ctwo{9.4}/\ctwo{15.0}  &  \cthree{92}/\cthree{8.5}/\cthree{13.5}  &  \textbf{\cone{98}}/\textbf{\cone{10.5}}/\textbf{\cone{18.2}}  \\
pepper+70\%  &  33/-2.1/-5.4  &  61/3.5/5.2  &  43/2.3/4.0  &  54/3.4/6.3  &  \ctwo{94}/\ctwo{8.4}/\cthree{11.4}  &  \cthree{90}/\cthree{7.9}/\ctwo{12.2}  &  \textbf{\cone{97}}/\textbf{\cone{10.0}}/\textbf{\cone{16.8}}  \\
pepper+90\%  &  17/-3.1/-6.4  &  31/0.9/1.3  &  27/0.7/1.2  &  32/1.1/\cthree{1.8}  &  \cthree{55}/\cthree{2.0}/1.5  &  \ctwo{60}/\ctwo{3.2}/\ctwo{4.9}  &  \textbf{\cone{92}}/\textbf{\cone{8.2}}/\textbf{\cone{11.7}}  \\
\hline
mandrill+10\%  &  58/1.1/-0.2  &  \cthree{67}/\cthree{2.9}/\cthree{4.7}  &  65/2.7/4.3  &  53/2.1/3.8  &  \cthree{67}/\cthree{2.9}/4.6  &  \ctwo{68}/\ctwo{3.1}/\ctwo{5.0}  &  \textbf{\cone{85}}/\textbf{\cone{5.0}}/\textbf{\cone{9.4}}  \\
mandrill+30\%  &  47/-0.9/-3.7  &  66/\cthree{2.8}/4.5  &  62/2.4/4.0  &  52/2.1/3.7  &  \ctwo{68}/\ctwo{3.0}/\cthree{4.6}  &  \cthree{67}/\ctwo{3.0}/\ctwo{4.8}  &  \textbf{\cone{76}}/\textbf{\cone{4.0}}/\textbf{\cone{6.8}}  \\
mandrill+50\%  &  36/-2.3/-5.7  &  64/\cthree{2.6}/\cthree{4.2}  &  54/1.9/3.3  &  51/2.0/3.4  &  \ctwo{68}/\ctwo{2.9}/\ctwo{4.6}  &  \cthree{66}/\ctwo{2.9}/\ctwo{4.6}  &  \textbf{\cone{74}}/\textbf{\cone{3.7}}/\textbf{\cone{6.3}}  \\
mandrill+70\%  &  25/-3.3/-7.0  &  54/1.9/3.2  &  43/1.1/2.2  &  47/1.7/3.1  &  \ctwo{67}/\ctwo{2.8}/\cthree{4.2}  &  \cthree{65}/\cthree{2.7}/\ctwo{4.4}  &  \textbf{\cone{71}}/\textbf{\cone{3.4}}/\textbf{\cone{5.4}}  \\
mandrill+90\%  &  14/-4.2/-8.1  &  38/0.4/0.7  &  36/0.4/0.8  &  35/0.5/0.9  &  \ctwo{50}/\ctwo{1.3}/\ctwo{1.4}  &  \cthree{48}/\cthree{1.0}/\cthree{1.1}  &  \textbf{\cone{66}}/\textbf{\cone{2.8}}/\textbf{\cone{4.3}}  \\
\hline
lake+10\%  &  70/4.3/3.5  &  83/7.5/11.5  &  83/7.4/11.4  &  82/6.6/11.3  &  \ctwo{89}/\ctwo{8.6}/\ctwo{13.8}  &  \cthree{84}/\cthree{7.7}/\cthree{12.1}  &  \textbf{\cone{97}}/\textbf{\cone{10.0}}/\textbf{\cone{17.9}}  \\
lake+30\%  &  56/1.7/-0.5  &  80/7.0/10.6  &  74/5.8/9.6  &  80/6.3/10.6  &  \ctwo{88}/\ctwo{8.4}/\ctwo{13.3}  &  \cthree{83}/\cthree{7.5}/\cthree{11.6}  &  \textbf{\cone{94}}/\textbf{\cone{9.5}}/\textbf{\cone{16.2}}  \\
lake+50\%  &  42/0.1/-2.6  &  73/6.0/9.3  &  45/4.0/7.0  &  77/6.0/9.8  &  \ctwo{88}/\ctwo{8.1}/\ctwo{11.8}  &  \cthree{82}/\cthree{7.2}/\cthree{11.1}  &  \textbf{\cone{92}}/\textbf{\cone{9.1}}/\textbf{\cone{15.1}}  \\
lake+70\%  &  29/-1.0/-4.0  &  40/2.9/5.1  &  27/2.3/4.0  &  51/4.1/6.8  &  \ctwo{84}/\ctwo{7.4}/\ctwo{10.8}  &  \cthree{79}/\cthree{6.8}/\cthree{10.3}  &  \textbf{\cone{89}}/\textbf{\cone{8.5}}/\textbf{\cone{13.5}}  \\
lake+90\%  &  15/-2.0/-5.0  &  18/0.7/1.2  &  17/0.7/1.3  &  18/0.9/\cthree{1.6}  &  \cthree{32}/\cthree{1.4}/1.5  &  \ctwo{55}/\ctwo{3.8}/\ctwo{5.4}  &  \textbf{\cone{81}}/\textbf{\cone{6.8}}/\textbf{\cone{9.8}}  \\
\hline
jetplane+10\%  &  76/2.8/0.6  &  88/6.7/9.9  &  \cthree{89}/6.4/9.8  &  66/2.8/6.5  &  \ctwo{93}/\ctwo{7.9}/\ctwo{12.5}  &  \cthree{89}/\cthree{6.8}/\cthree{10.5}  &  \textbf{\cone{98}}/\textbf{\cone{9.1}}/\textbf{\cone{16.6}}  \\
jetplane+30\%  &  60/-0.2/-3.6  &  86/6.2/8.9  &  79/4.1/7.5  &  66/2.7/6.1  &  \ctwo{93}/\ctwo{7.8}/\ctwo{11.8}  &  \cthree{88}/\cthree{6.6}/\cthree{9.9}  &  \textbf{\cone{97}}/\textbf{\cone{8.8}}/\textbf{\cone{15.6}}  \\
jetplane+50\%  &  45/-1.9/-5.7  &  81/5.0/7.5  &  44/1.9/4.2  &  51/2.5/5.6  &  \ctwo{91}/\ctwo{7.1}/\ctwo{10.3}  &  \cthree{87}/\cthree{6.4}/\cthree{9.5}  &  \textbf{\cone{95}}/\textbf{\cone{8.4}}/\textbf{\cone{14.1}}  \\
jetplane+70\%  &  30/-3.1/-7.1  &  39/0.7/2.2  &  25/0.0/1.0  &  32/0.8/2.2  &  \ctwo{89}/\ctwo{6.4}/\cthree{8.4}  &  \cthree{85}/\cthree{5.9}/\ctwo{8.7}  &  \textbf{\cone{93}}/\textbf{\cone{7.7}}/\textbf{\cone{12.2}}  \\
jetplane+90\%  &  15/-4.1/-8.2  &  16/-1.6/-2.1  &  16/-1.6/-2.0  &  22/-2.0/-3.0  &  \cthree{30}/\cthree{-1.1}/\cthree{-1.8}  &  \ctwo{56}/\ctwo{2.0}/\ctwo{3.1}  &  \textbf{\cone{86}}/\textbf{\cone{5.8}}/\textbf{\cone{7.8}}  \\
\hline
\end{tabular}}
\end{center}
\end{table*}

\begin{figure}[!t]
\begin{center}

\subfloat{\fcolorbox{colorone}{colortwo}{\includegraphics[width=\figurewidth, height=\figureheight]{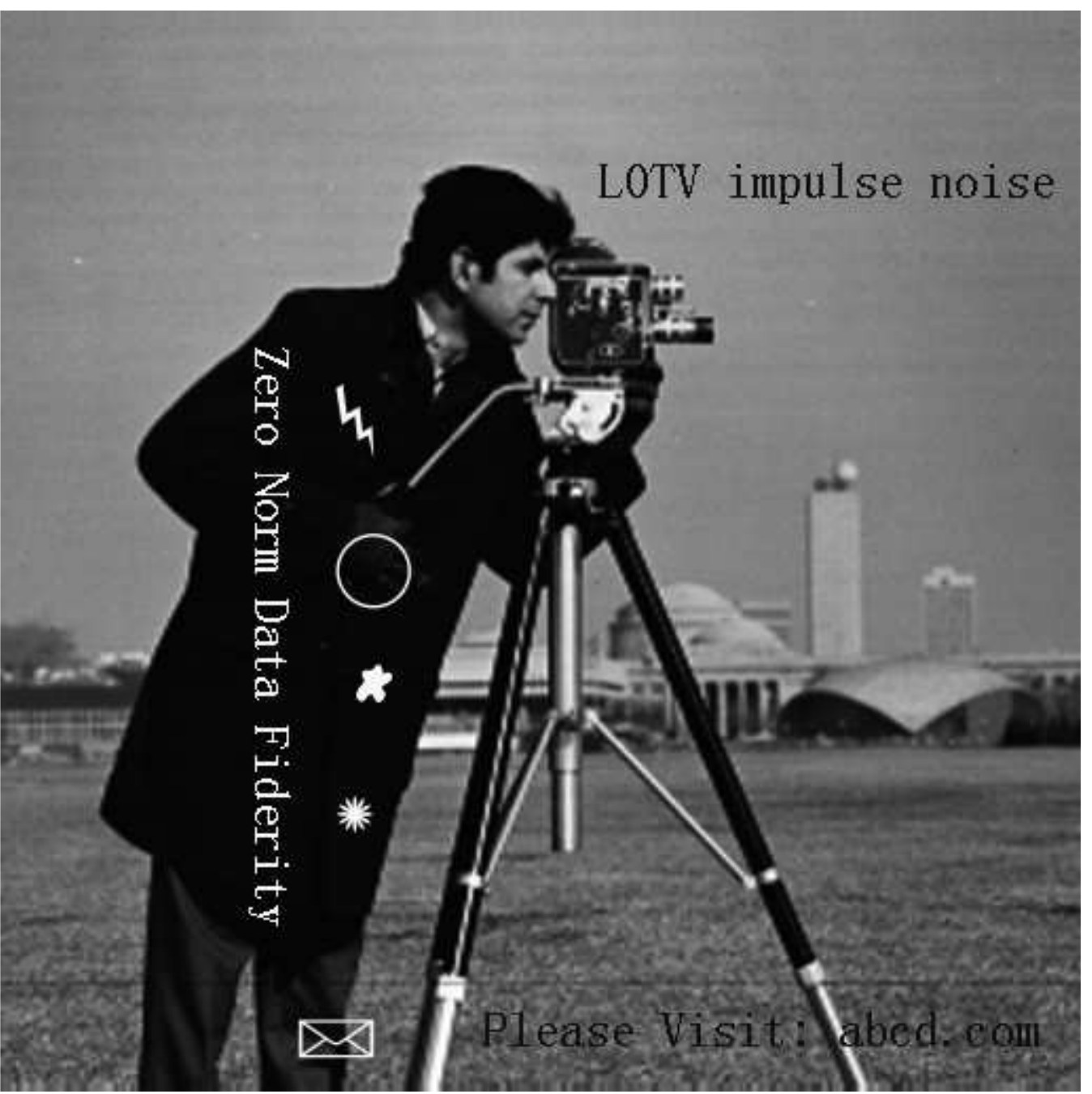}}}\ghs
\subfloat{\fcolorbox{colorone}{colortwo}{\includegraphics[width=\figurewidth, height=\figureheight]{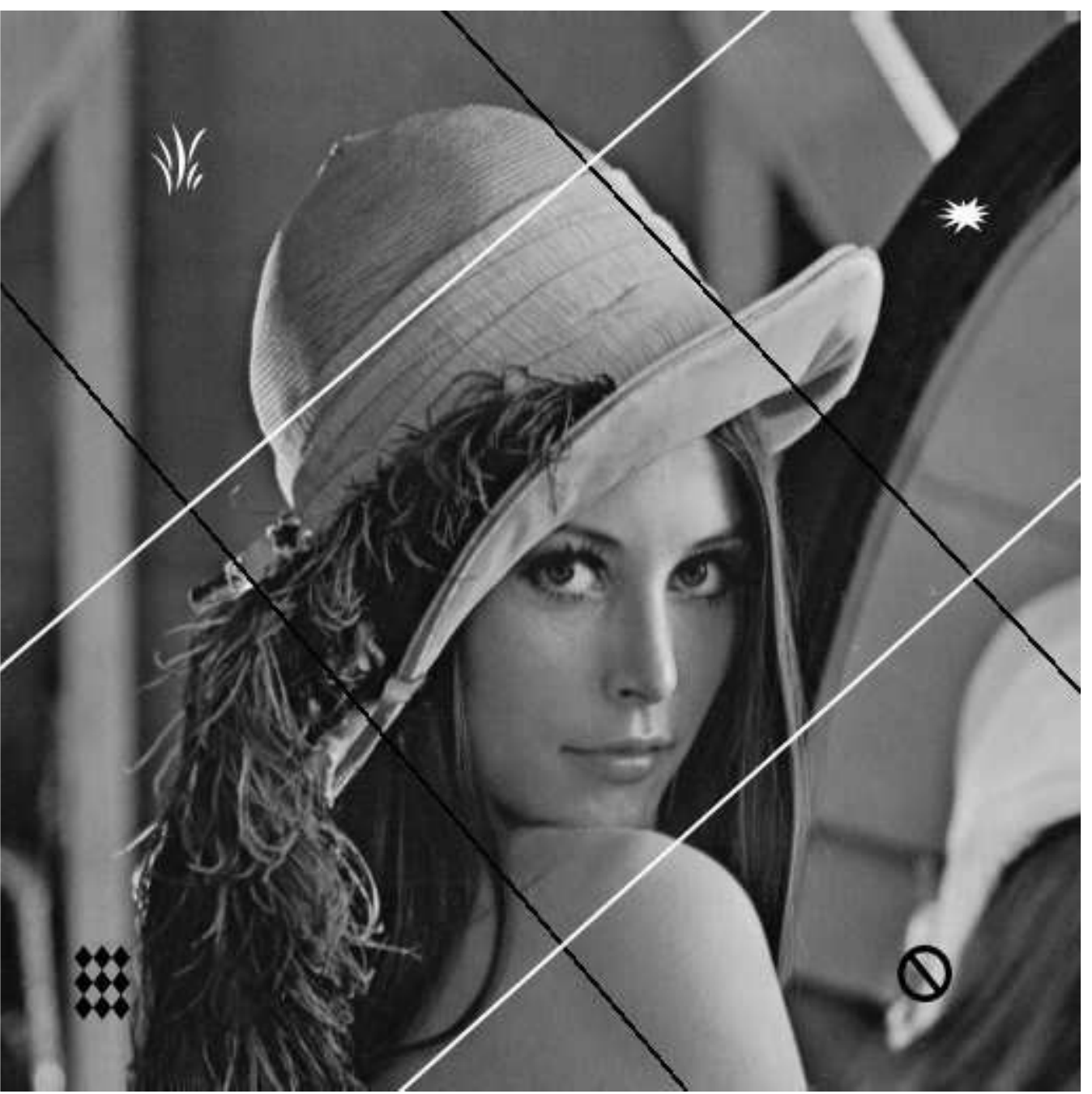}}}\ghs
\subfloat{\fcolorbox{colorone}{colortwo}{\includegraphics[width=\figurewidth, height=\figureheight]{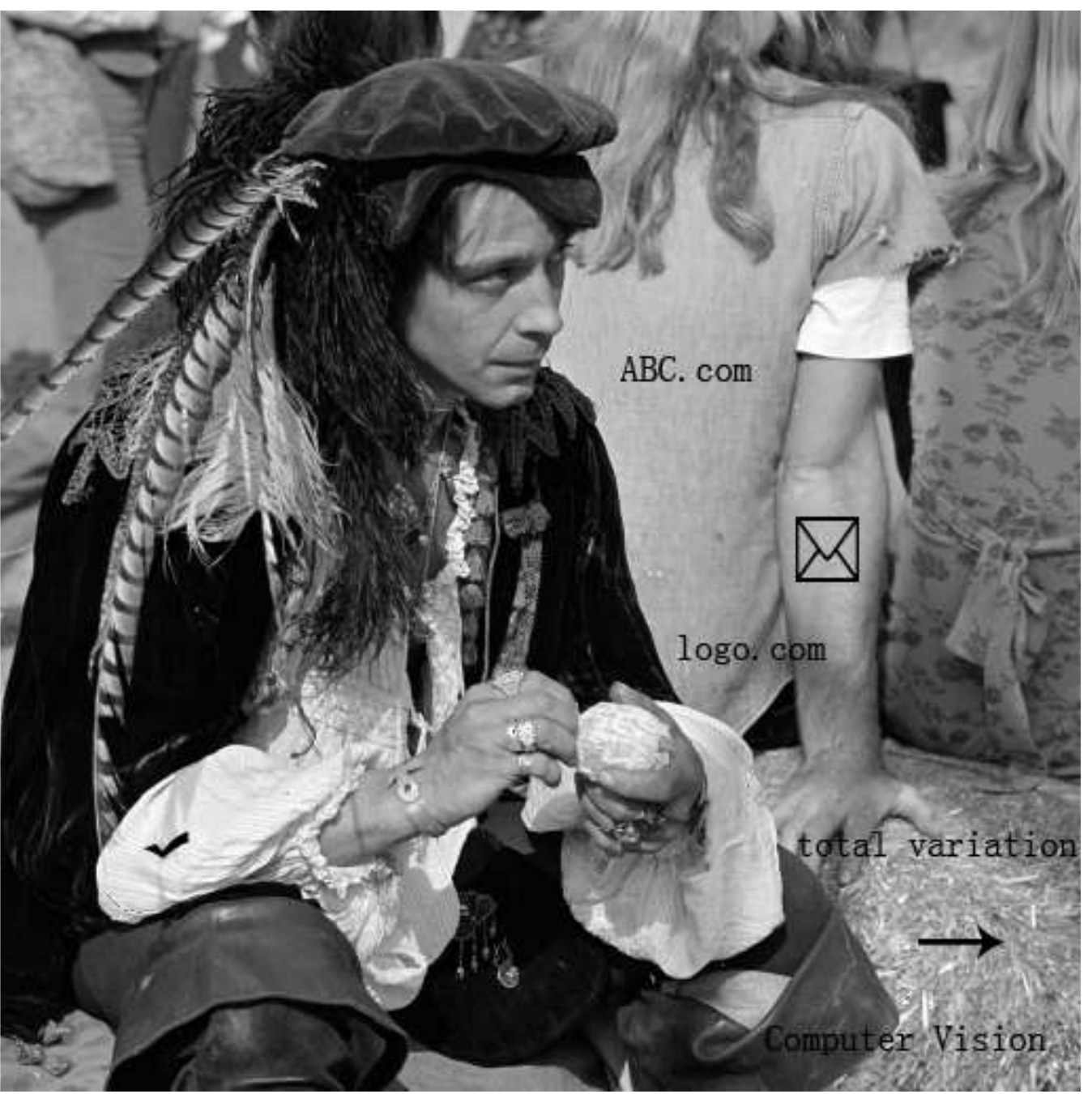}}}
\caption{Sample images in scratched image denoising problems.}
\label{fig:scratched:image}
\vspace{-5pt}
\subfloat{\fcolorbox{colorone}{colortwo}{\includegraphics[width=\figurewidth, height=\figureheight]{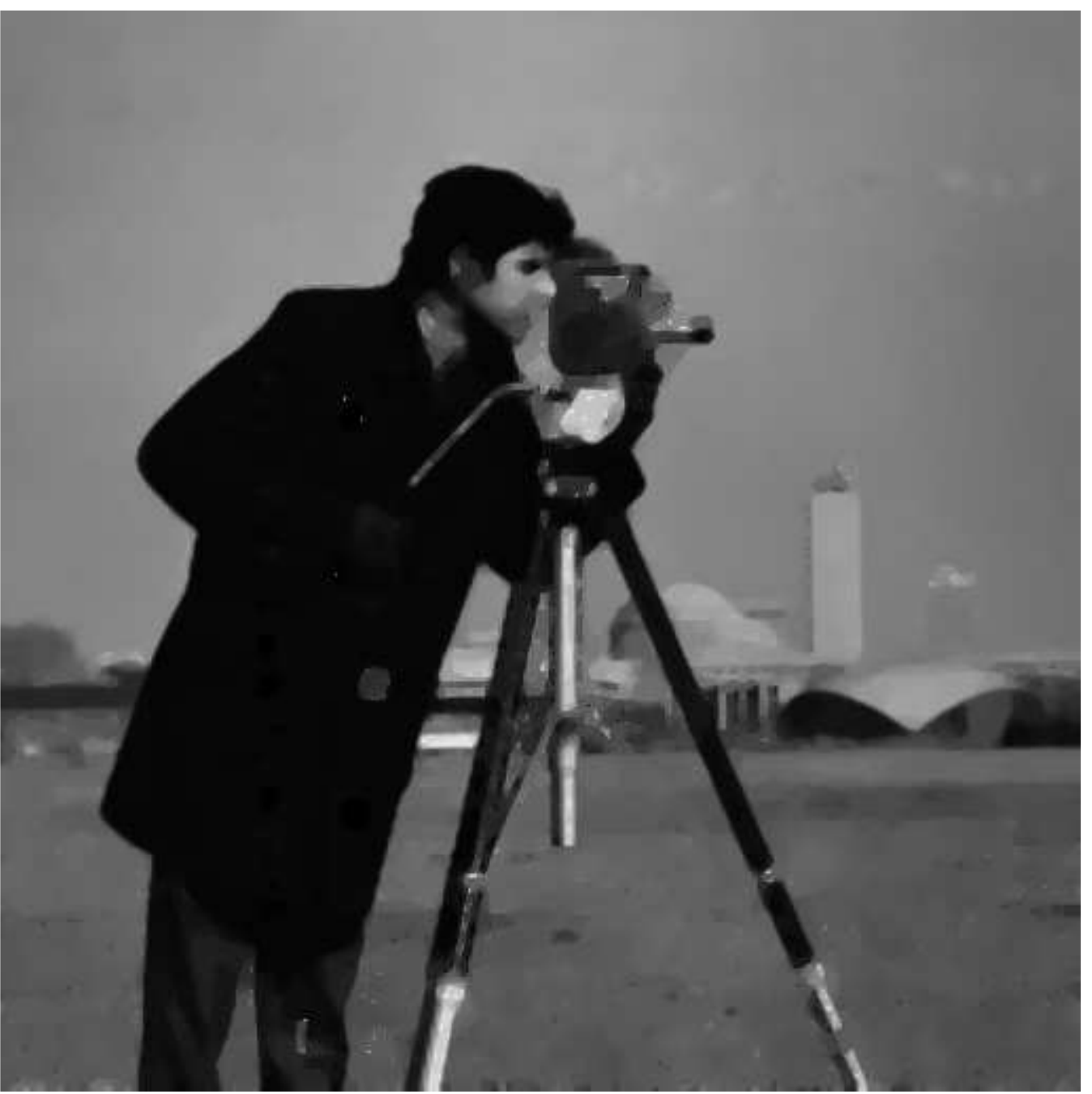}}}\ghs
\subfloat{\fcolorbox{colorone}{colortwo}{\includegraphics[width=\figurewidth, height=\figureheight]{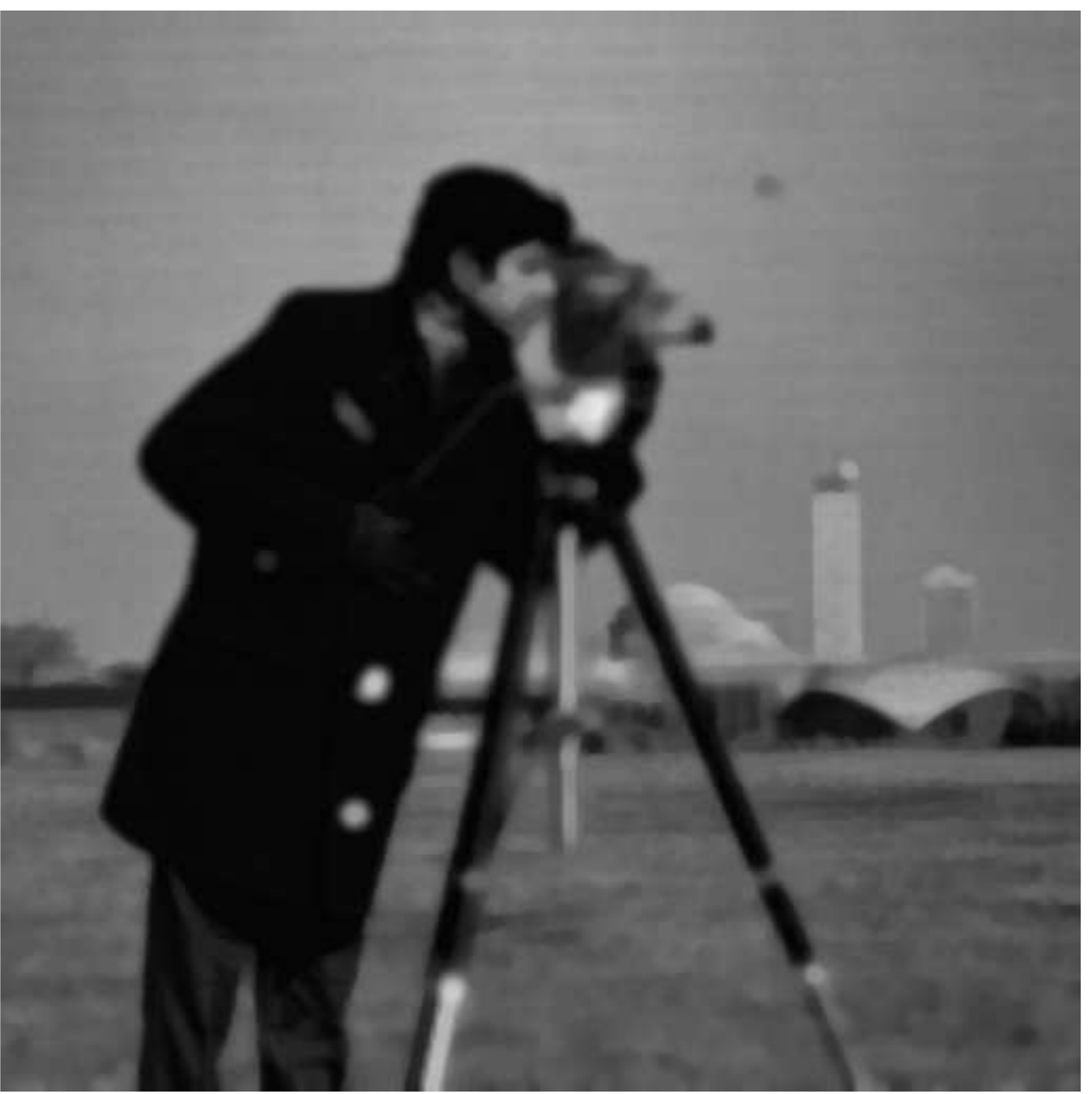}}}\ghs
\subfloat{\fcolorbox{colorone}{colortwo}{\includegraphics[width=\figurewidth, height=\figureheight]{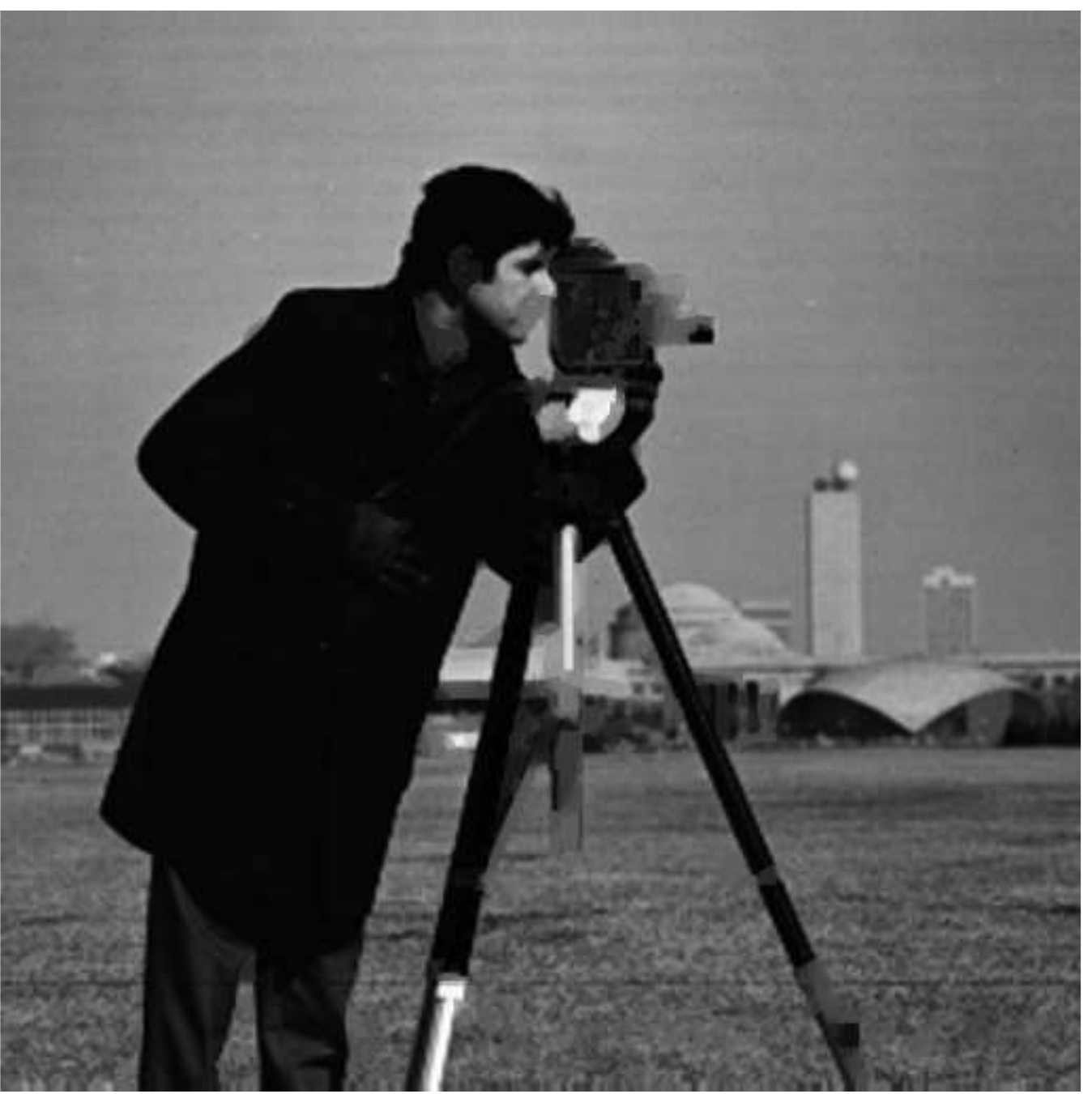}}}
\\
\vspace{-8pt}
\subfloat{\fcolorbox{colorone}{colortwo}{\includegraphics[width=\figurewidth, height=\figureheight]{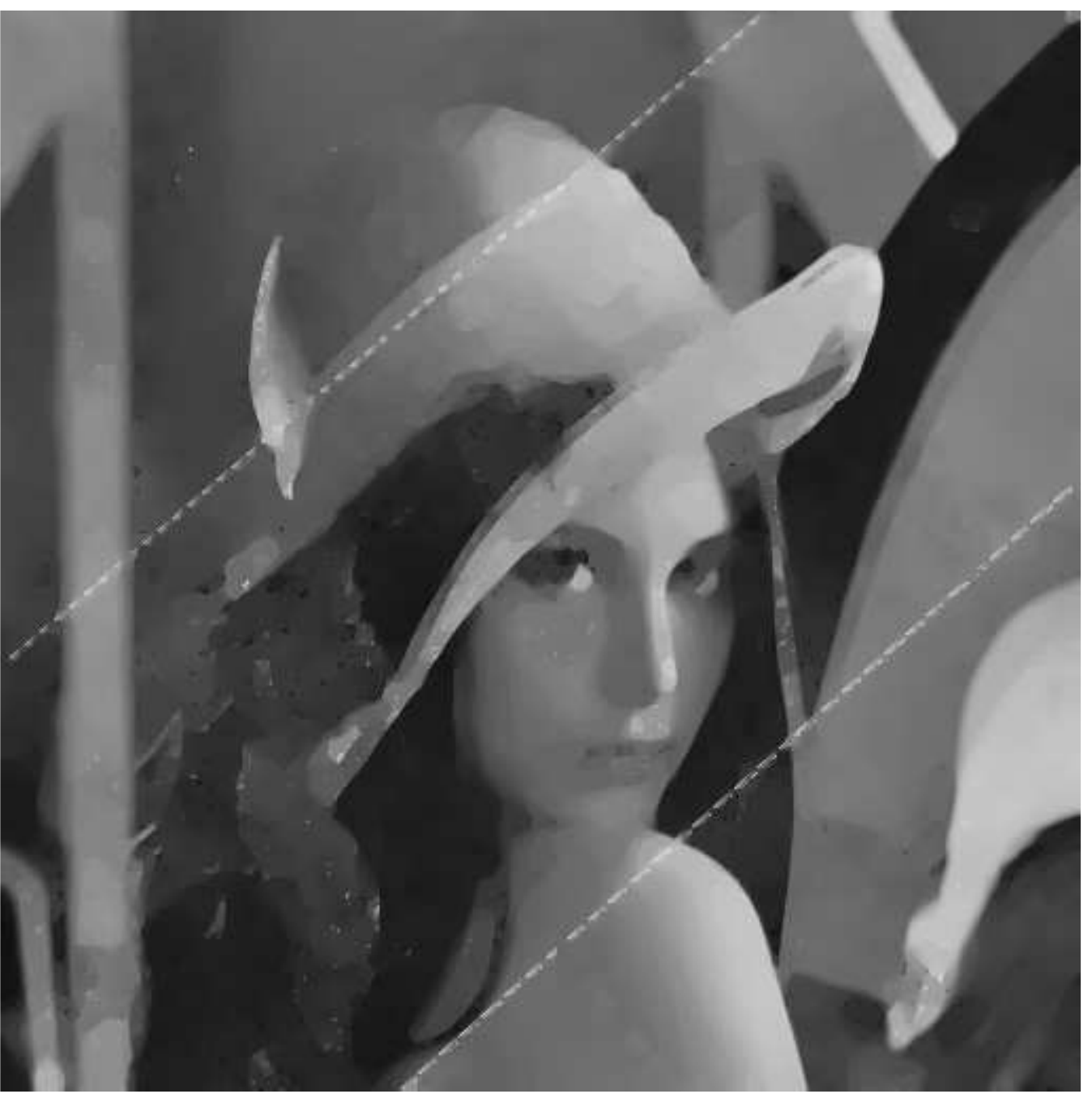}}}\ghs
\subfloat{\fcolorbox{colorone}{colortwo}{\includegraphics[width=\figurewidth, height=\figureheight]{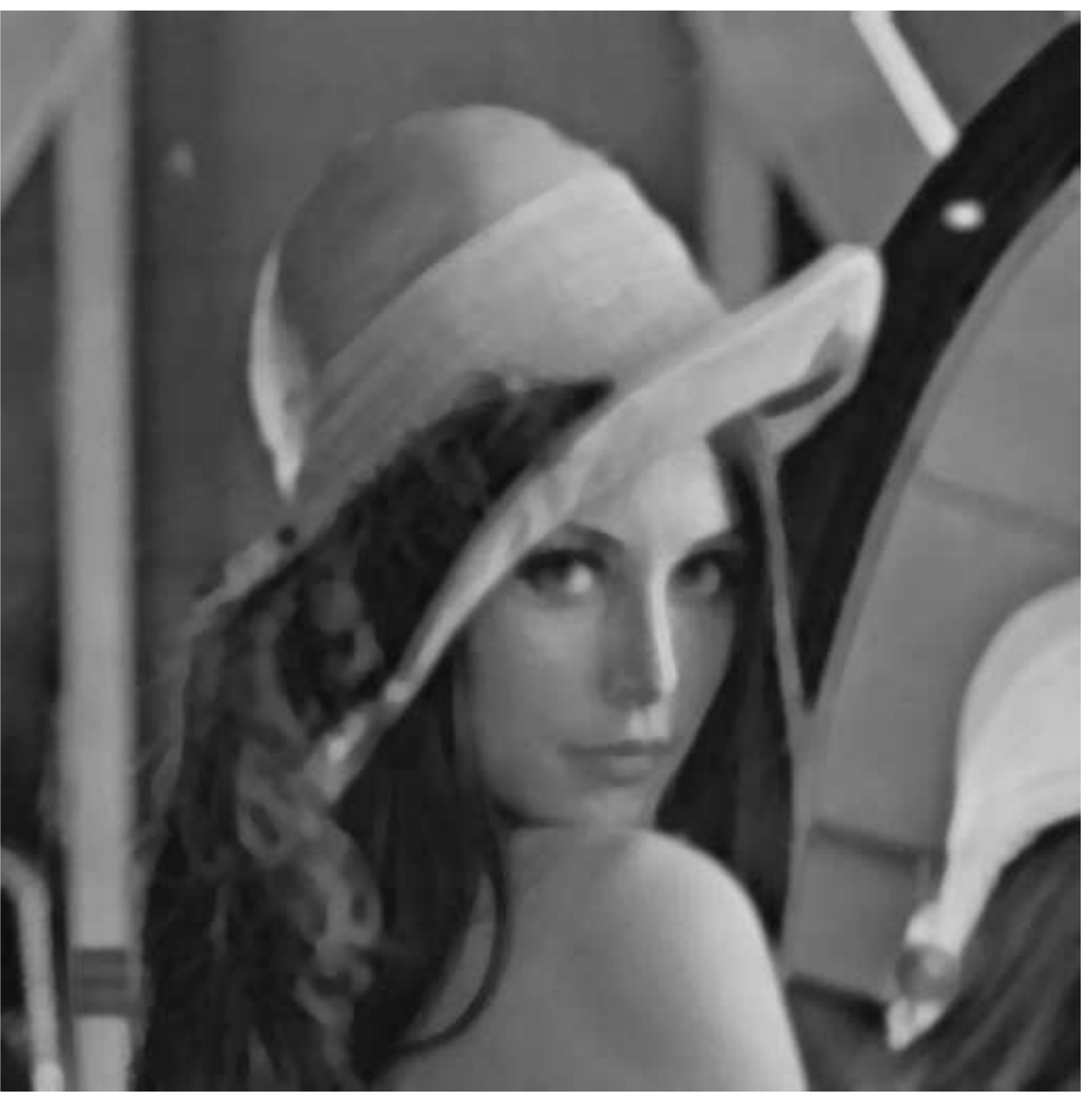}}}\ghs
\subfloat{\fcolorbox{colorone}{colortwo}{\includegraphics[width=\figurewidth, height=\figureheight]{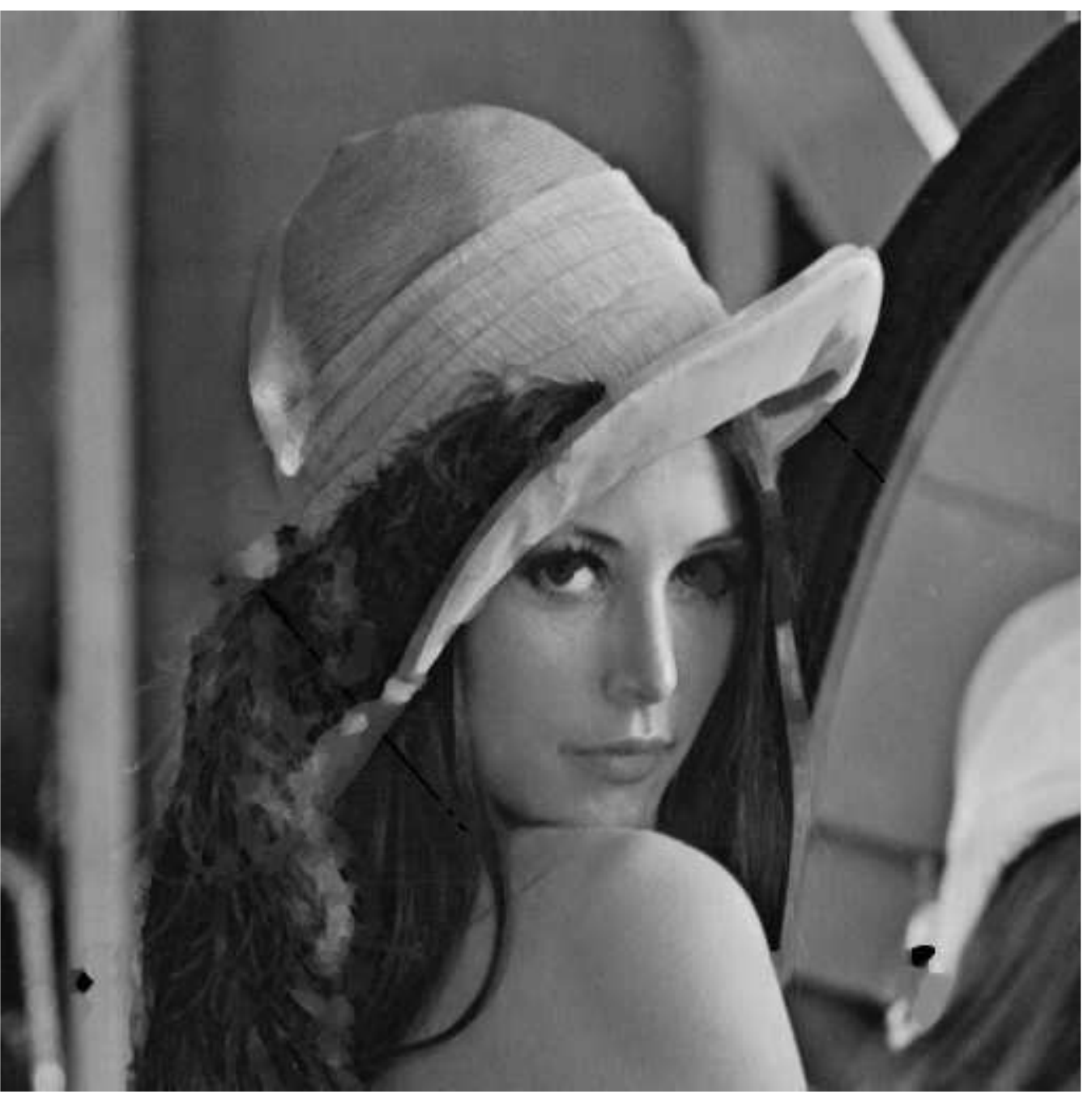}}}
\\
\vspace{-8pt}
\subfloat{\fcolorbox{colorone}{colortwo}{\includegraphics[width=\figurewidth, height=\figureheight]{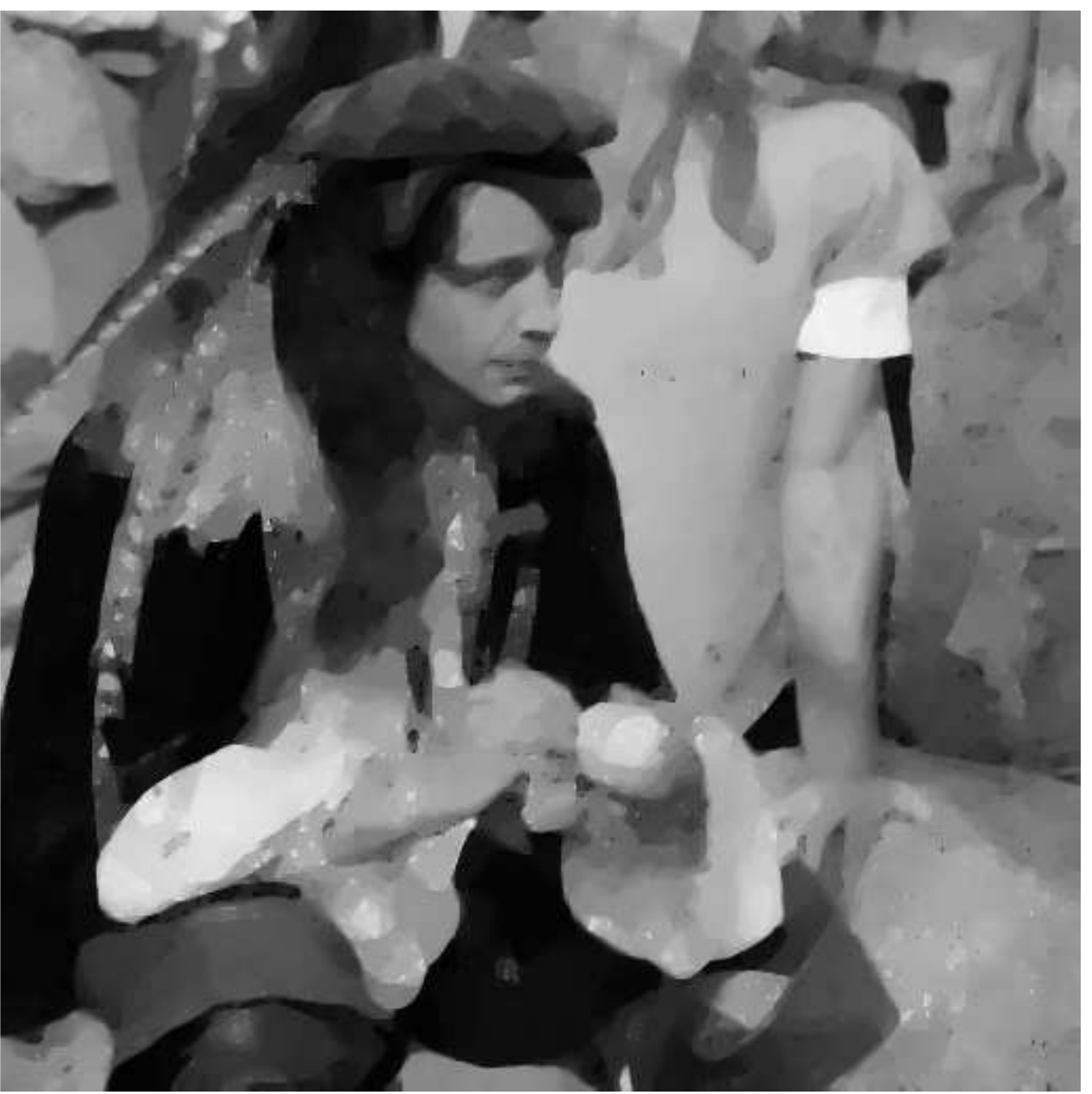}}}\ghs
\subfloat{\fcolorbox{colorone}{colortwo}{\includegraphics[width=\figurewidth, height=\figureheight]{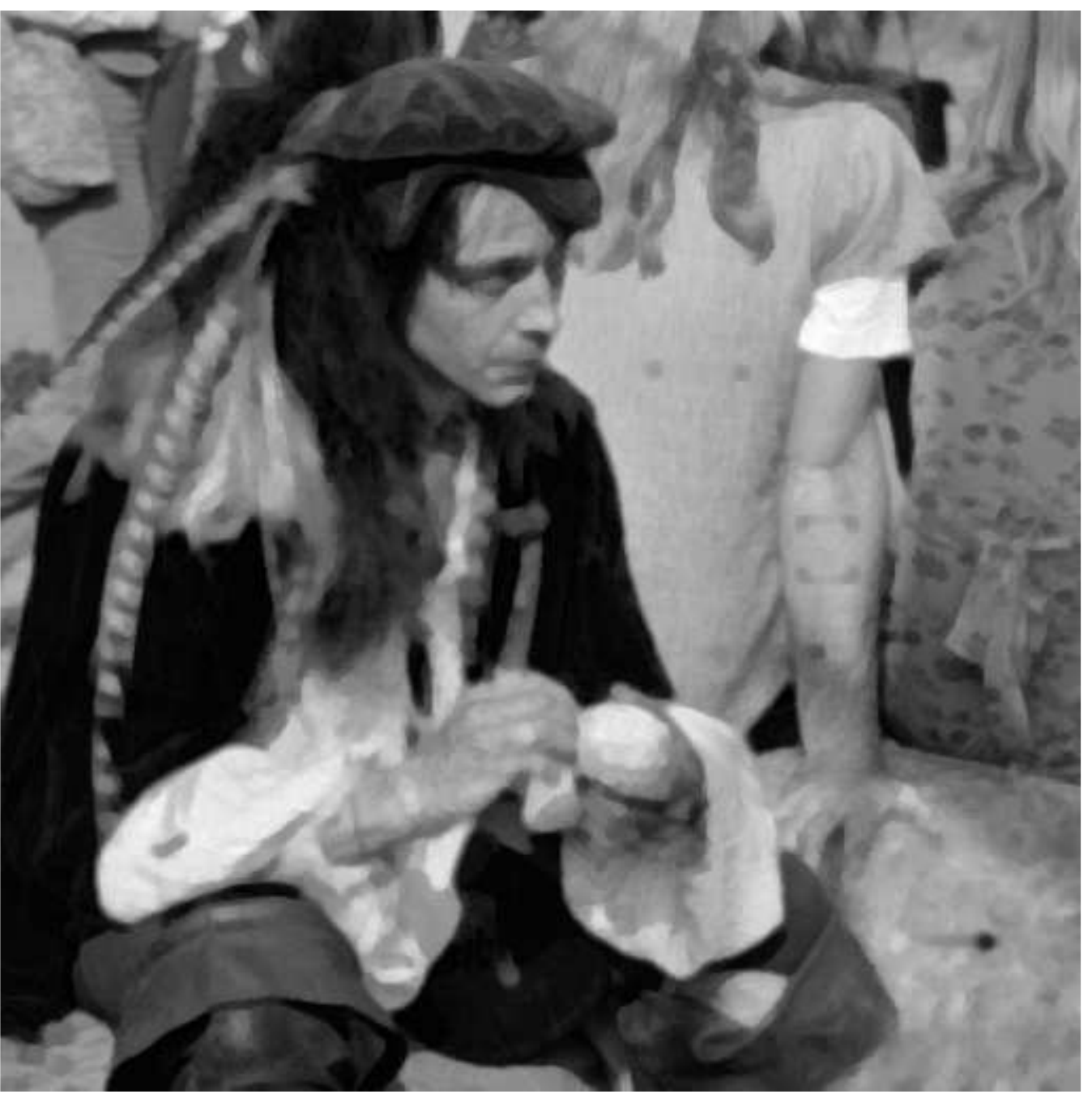}}}\ghs
\subfloat{\fcolorbox{colorone}{colortwo}{\includegraphics[width=\figurewidth, height=\figureheight]{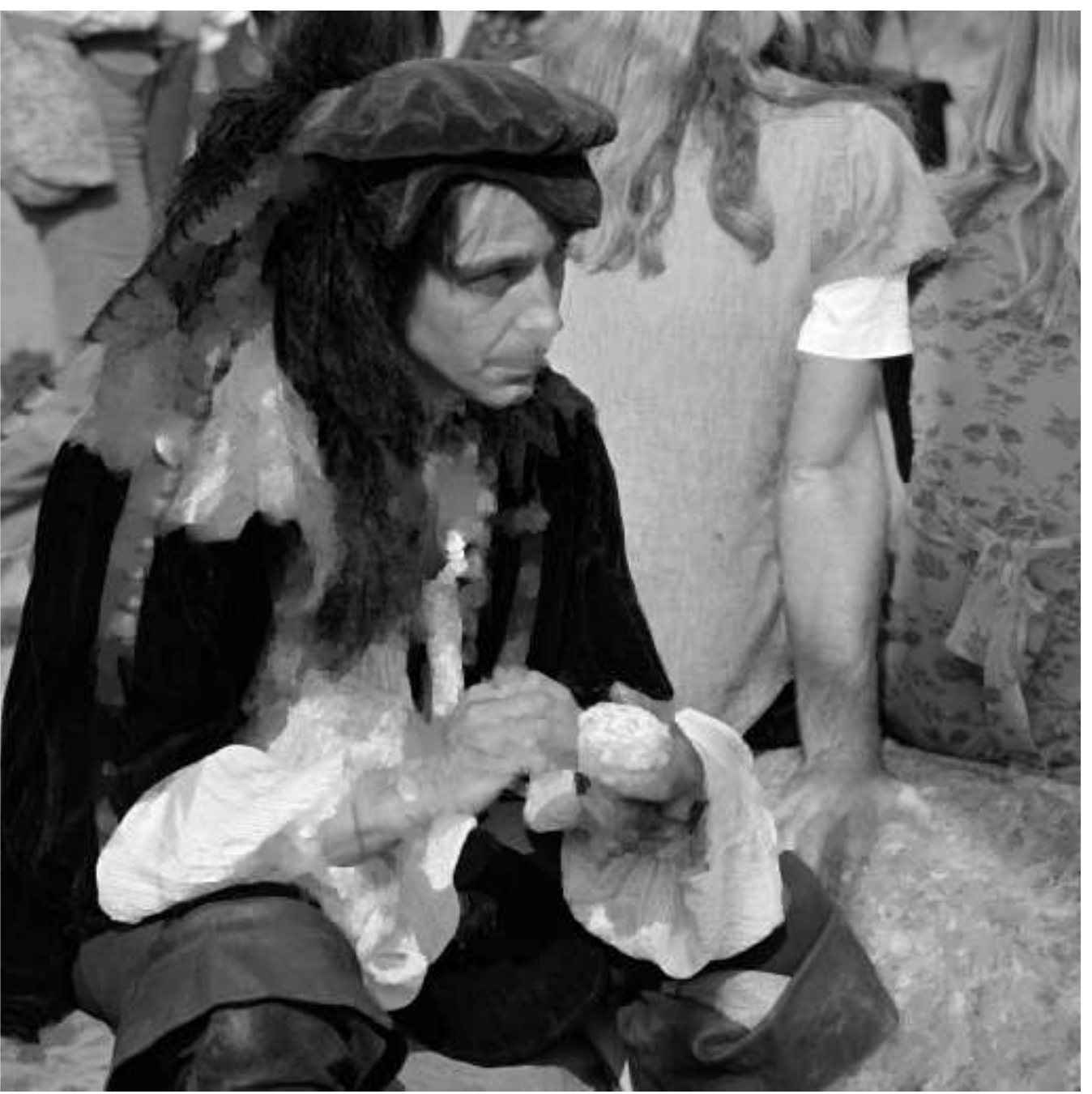}}}
\caption{Recovered images in scratched image denoising problems. First column: $\ell_{02}TV$-AOP, second column: $\ell_{0}TV$-PDA, third column: $\ell_{0}TV$-PADMM.}\label{fig:scratched:image1}
\end{center}
\end{figure}


\subsection{General Image Deblurring Problems}
In this subsection, we demonstrate the performance of all methods with their optimal regularization parameters on general deblurring problems. Table \ref{deblurring:rv} shows the recovery results for random-valued impulse noise, salt-and-pepper impulse noise, and mixed impulse noise, respectively. Figure \ref{fig:deblurring:lambda} shows image recovery results with varying the regularization parameter. We have the following interesting observations. (\textbf{i}) $\ell_{02}TV$-AOP significantly outperforms $\ell_1TV$-SBM, and the performance gap becomes larger as the noise level increases. This is because the key assumption in the $\ell_1$ model is that $Ku-b$ is sparse at the optimal solution $u^*$. This does not hold when the noise level is high. (\textbf{ii}) $\ell_0TV$-PDA outperforms $\ell_{02}TV$-AOP for high level ($\geq 30\%$) random-valued impulse noise. However, for salt-and-pepper impulse noise, $\ell_0TV$-PDA gives worse performance than $\ell_{02}TV$-AOP in most cases. This phenomenon indicates that the Penalty Decomposition Algorithm is not stable for deblurring problems. (\textbf{iii}) By contrast, our $\ell_0TV$-PADMM consistently outperforms all methods, especially when the noise level is large. We attribute this result to the ``lifting'' technique that is used in our optimization algorithm.

Finally, we also report the performance of all methods with sweeping the radius parameter $r$ as in (\ref{eq:varyingr}) over $\{1,4,7,...,20\}$ in Figure \ref{fig:deblurring:R}. We notice that the restoration quality degenerates as the radius of the kernel increases for all methods. However, our method consistently gives the best performance.

\subsection{Scratched Image Denoising Problems}
In this subsection, we demonstrate the superiority of the proposed $\ell_0TV$-PADMM in real-world image restoration problems. Specifically, we corrupt the images with scratches which can be viewed as impulse noise\footnote{Note that this is different from the classical image inpainting problem that assumes the mask is known. In our scratched image denoising problem, we assume the mask is unknown.}, see Figure \ref{fig:scratched:image}. We only consider recovering images using $\ell_{02}TV$-AOP, $\ell_0TV$-PDA and $\ell_0TV$-PADMM. We show the recovered results in Figure \ref{fig:scratched:image1}. For better visualization of the images recovered by all methods, we also show auxiliary images $\bbb{c}$ in Figure \ref{fig:scratched:image2}, which show the complement of the absolute residual between the recovered image $\bbb{u}$ and the corrupted image $\bbb{b}$ (i.e., $\bbb{c} = \{\bbb{1} - |\bbb{b}-\bbb{u}|\}$). Note that when $\bbb{c}_i$ is approximately equal to $1$, the color of the corresponding pixel at position $i$ in the image is white. A conclusion can be drawn that our method $\ell_0TV$-PADMM generates more `white' images $\bbb{c}$ than the other two methods, since it can identify the `right' outliers in the corrupted image and make the correction using their neighborhood information.

\subsection{Colored Image Denoising Problems}
\begin{figure}[!t]
\begin{center}
\subfloat{\fcolorbox{colorone}{colortwo}{\includegraphics[width=\figurewidth, height=\figureheight]{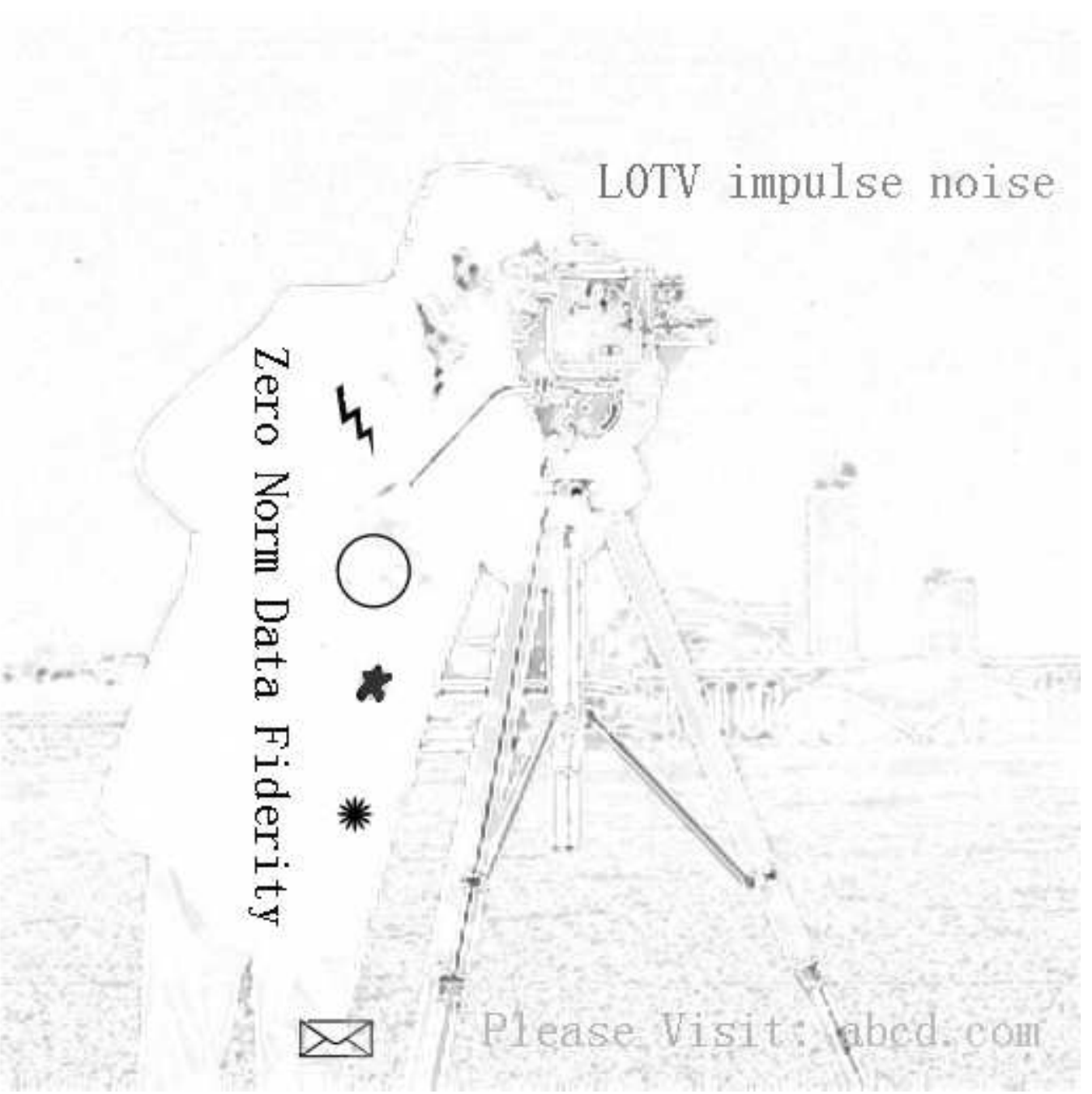}}}\ghs
\subfloat{\fcolorbox{colorone}{colortwo}{\includegraphics[width=\figurewidth, height=\figureheight]{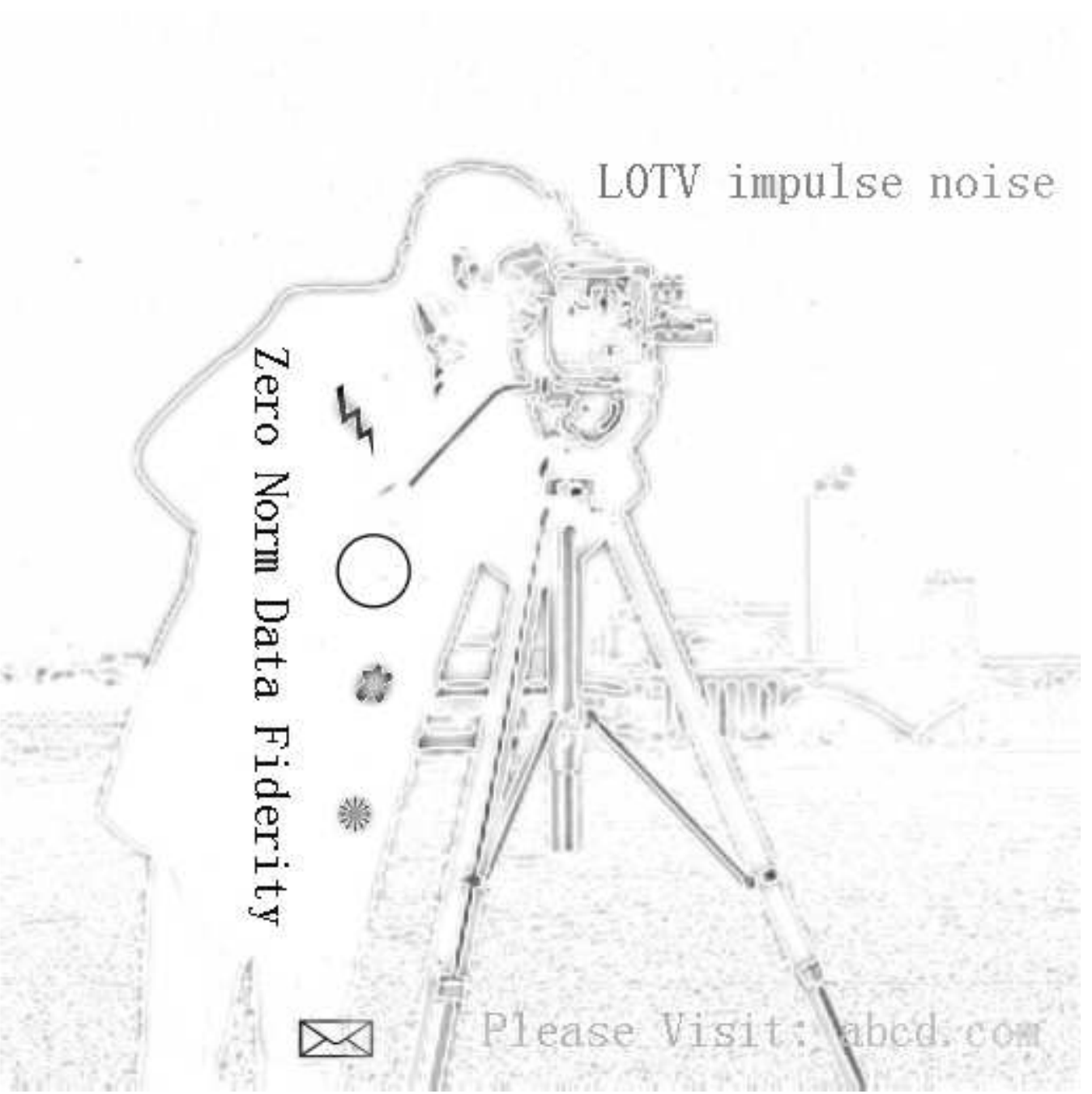}}}\ghs
\subfloat{\fcolorbox{colorone}{colortwo}{\includegraphics[width=\figurewidth, height=\figureheight]{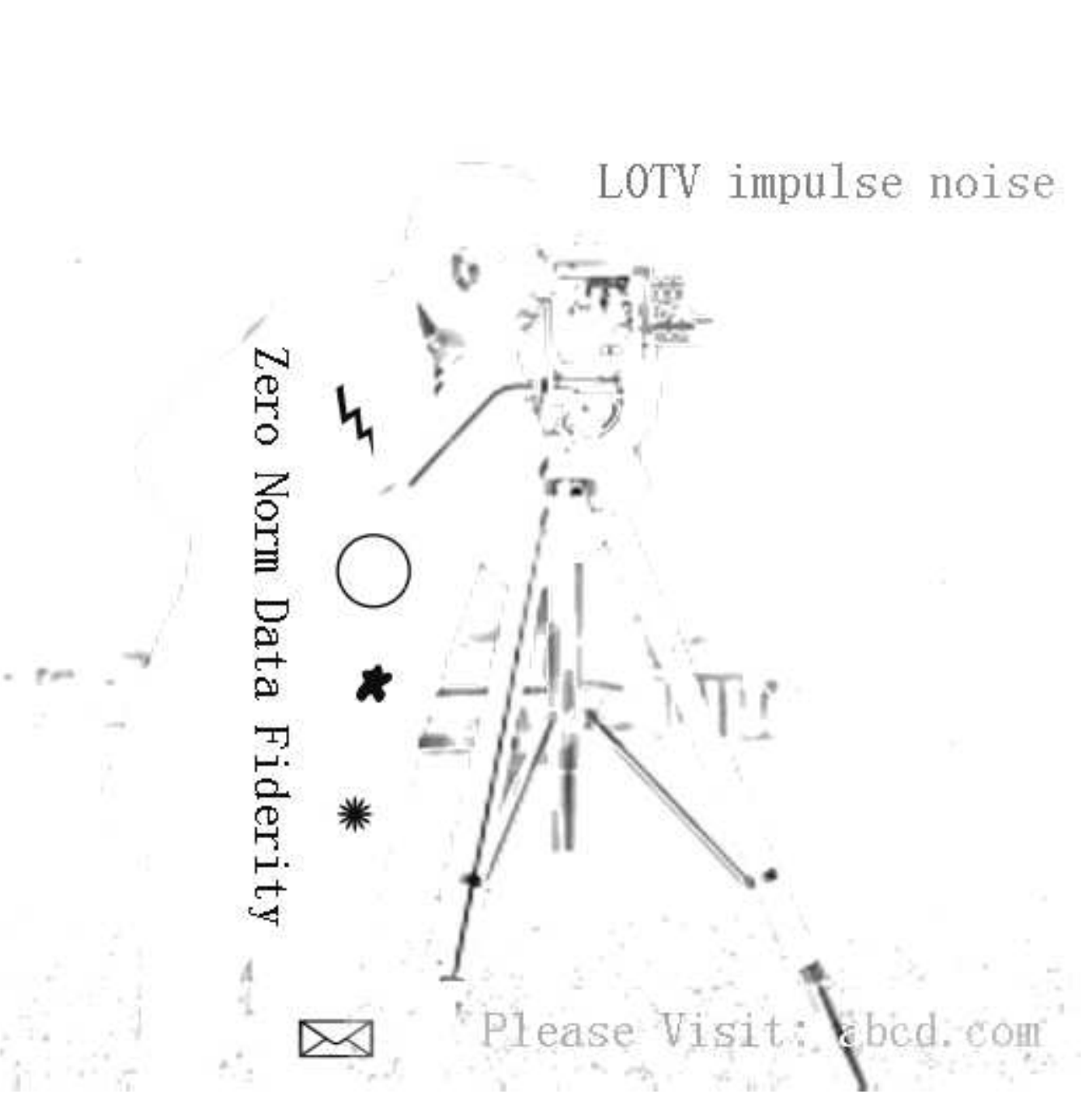}}}
\\
\vspace{-8pt}
\subfloat{\fcolorbox{colorone}{colortwo}{\includegraphics[width=\figurewidth, height=\figureheight]{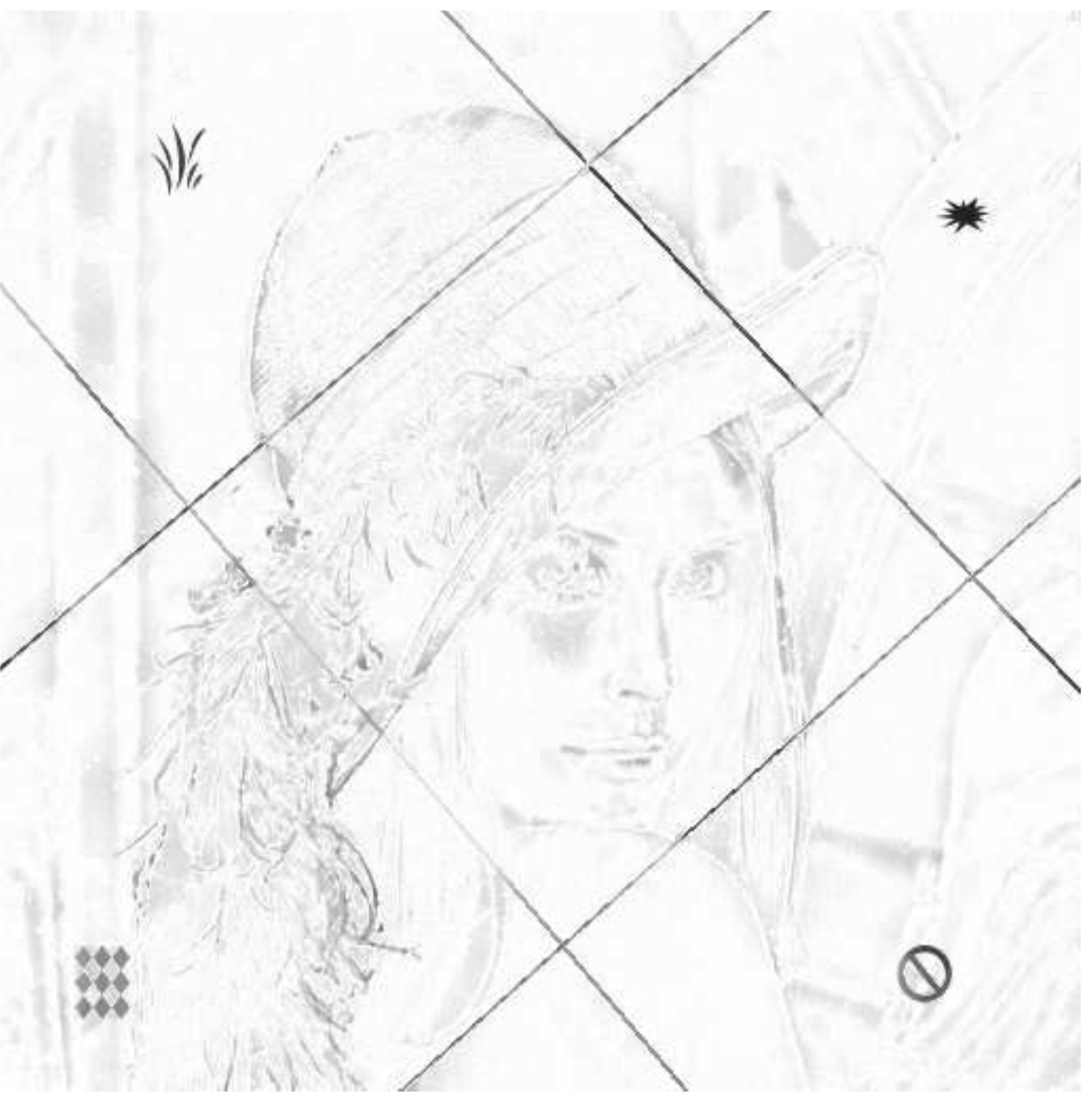}}}\ghs
\subfloat{\fcolorbox{colorone}{colortwo}{\includegraphics[width=\figurewidth, height=\figureheight]{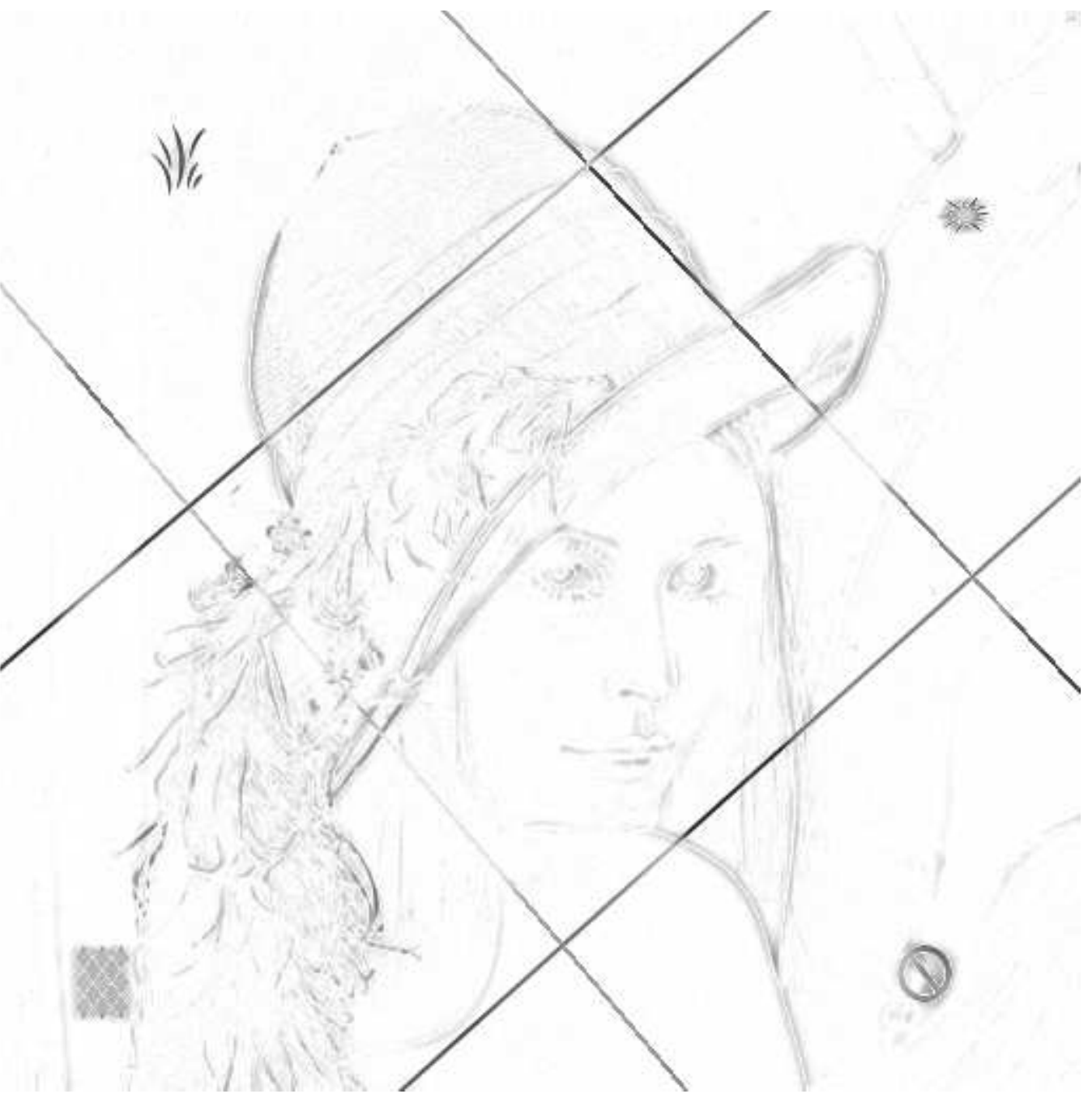}}}\ghs
\subfloat{\fcolorbox{colorone}{colortwo}{\includegraphics[width=\figurewidth, height=\figureheight]{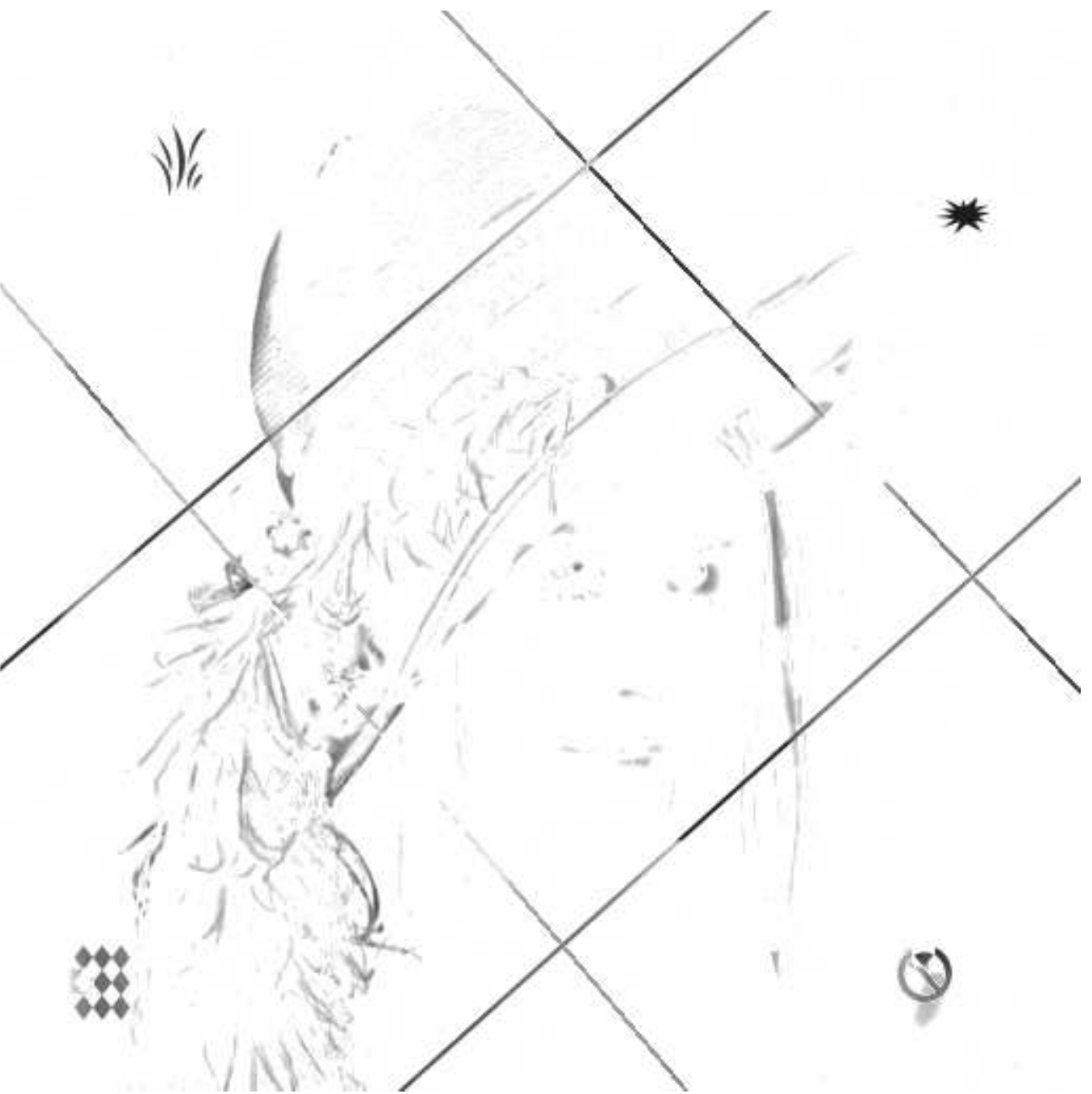}}}
\\
\vspace{-8pt}
\subfloat{\fcolorbox{colorone}{colortwo}{\includegraphics[width=\figurewidth, height=\figureheight]{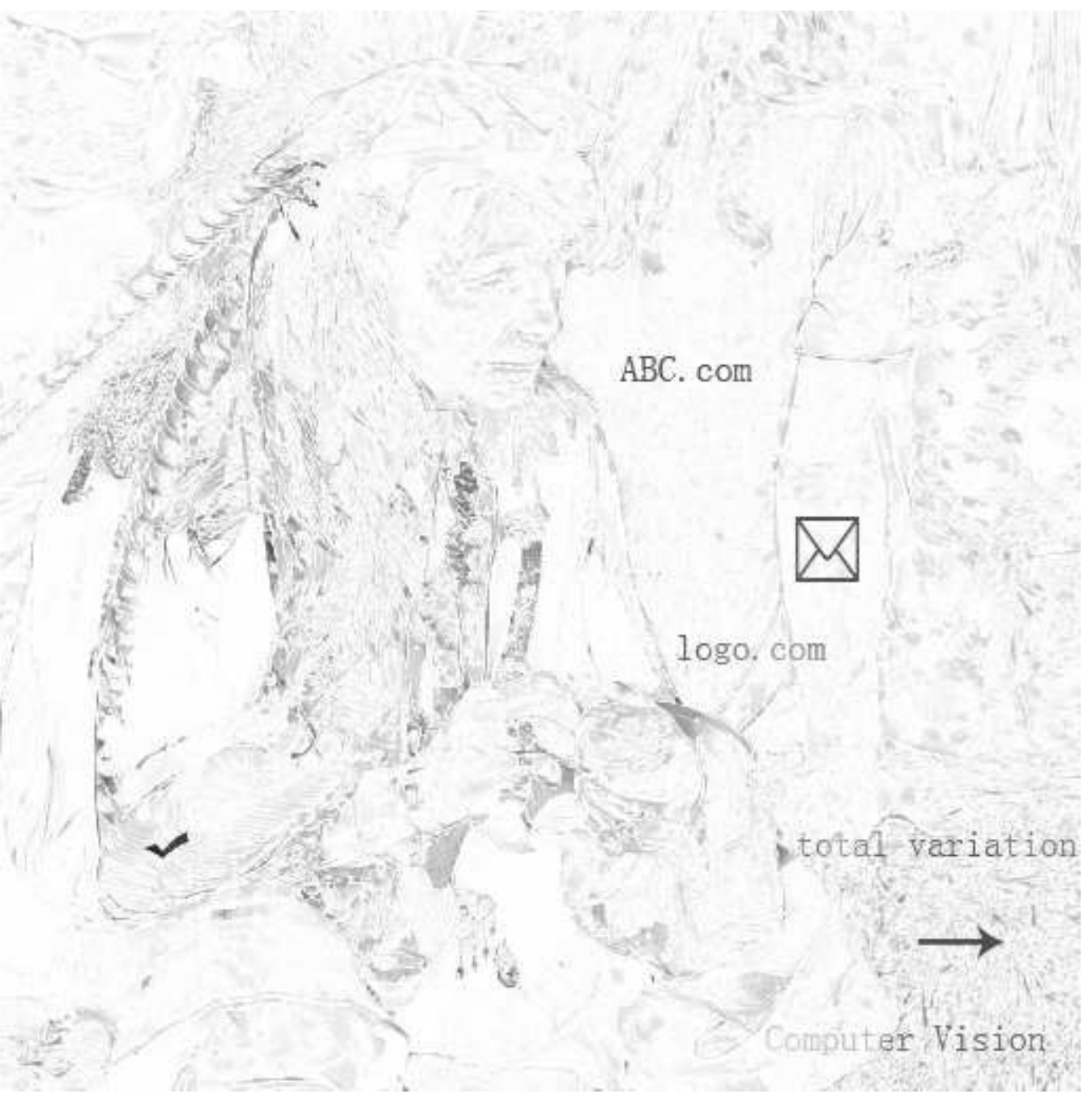}}}\ghs
\subfloat{\fcolorbox{colorone}{colortwo}{\includegraphics[width=\figurewidth, height=\figureheight]{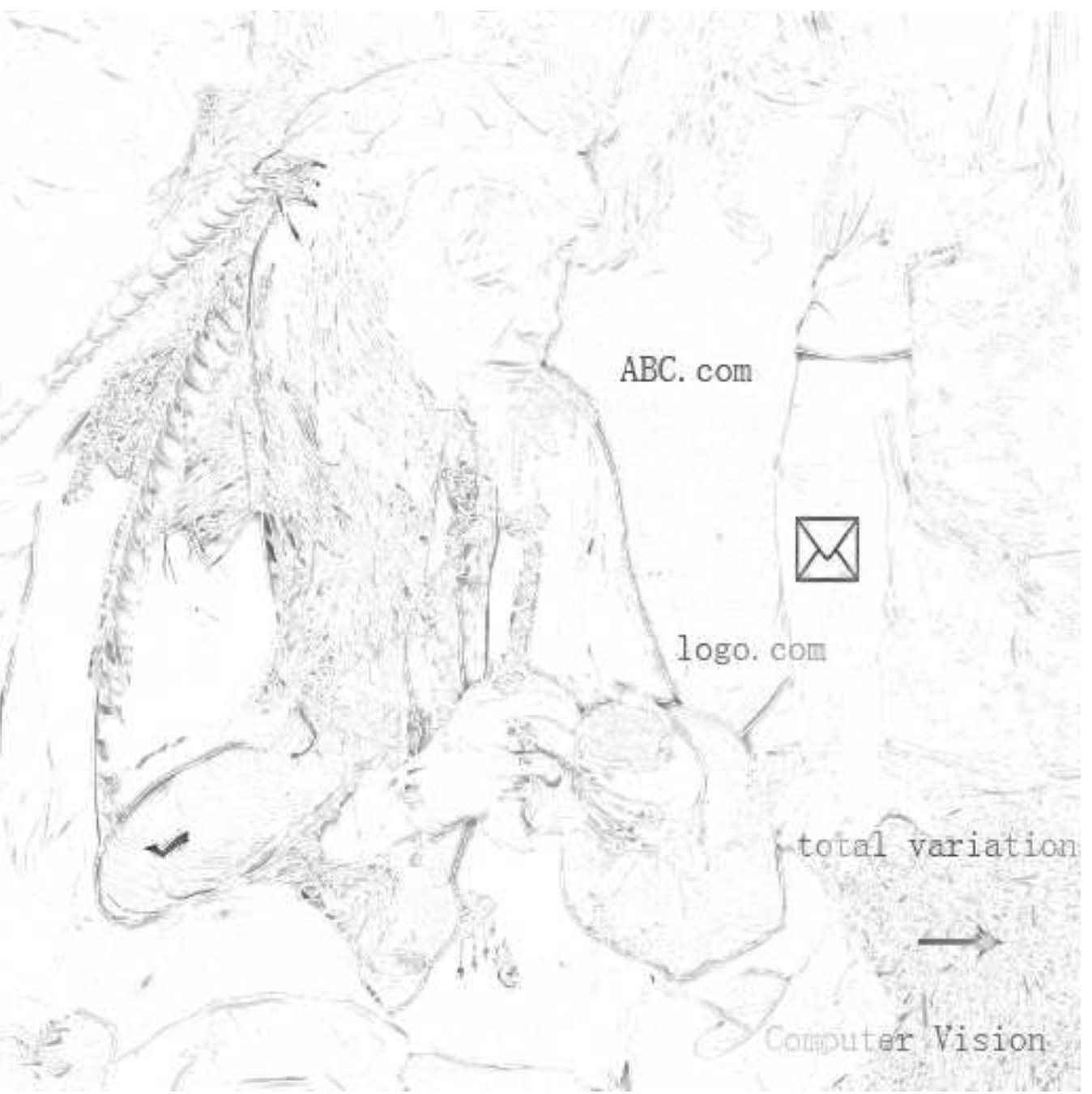}}}\ghs
\subfloat{\fcolorbox{colorone}{colortwo}{\includegraphics[width=\figurewidth, height=\figureheight]{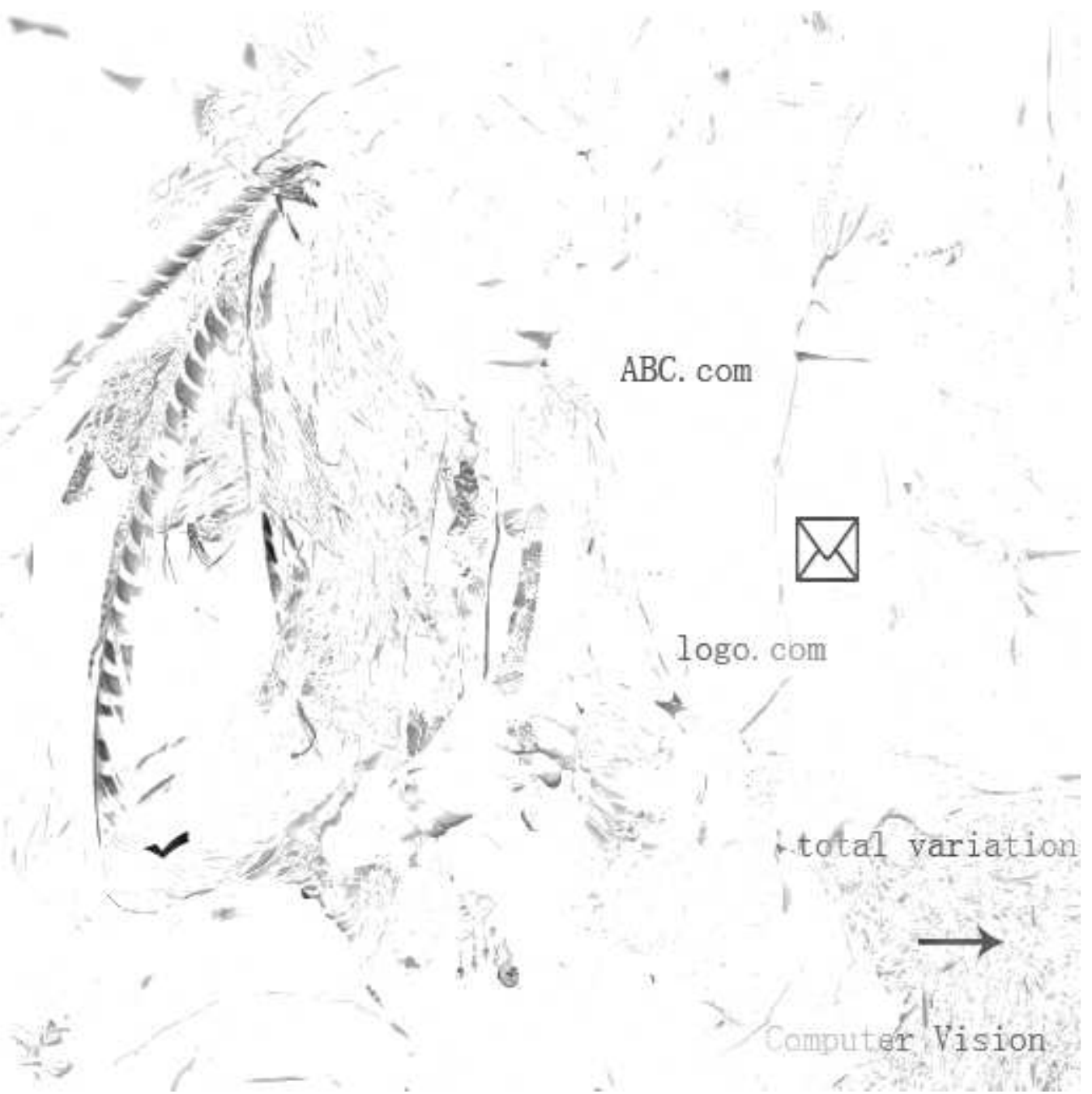}}}
\caption{Absolute residual (between scratched image and recovered image) in scratched image denoising problems. First column: $\ell_{02}TV$-AOP, second column: $\ell_{0}TV$-PDA, third column: $\ell_{0}TV$-PADMM.}\label{fig:scratched:image2}
\vspace{-5pt}
%
\subfloat[\figsizetwo clean `lenna']{\fcolorbox{colorone}{colortwo}{\includegraphics[width=\figurewidth, height=\figureheight]{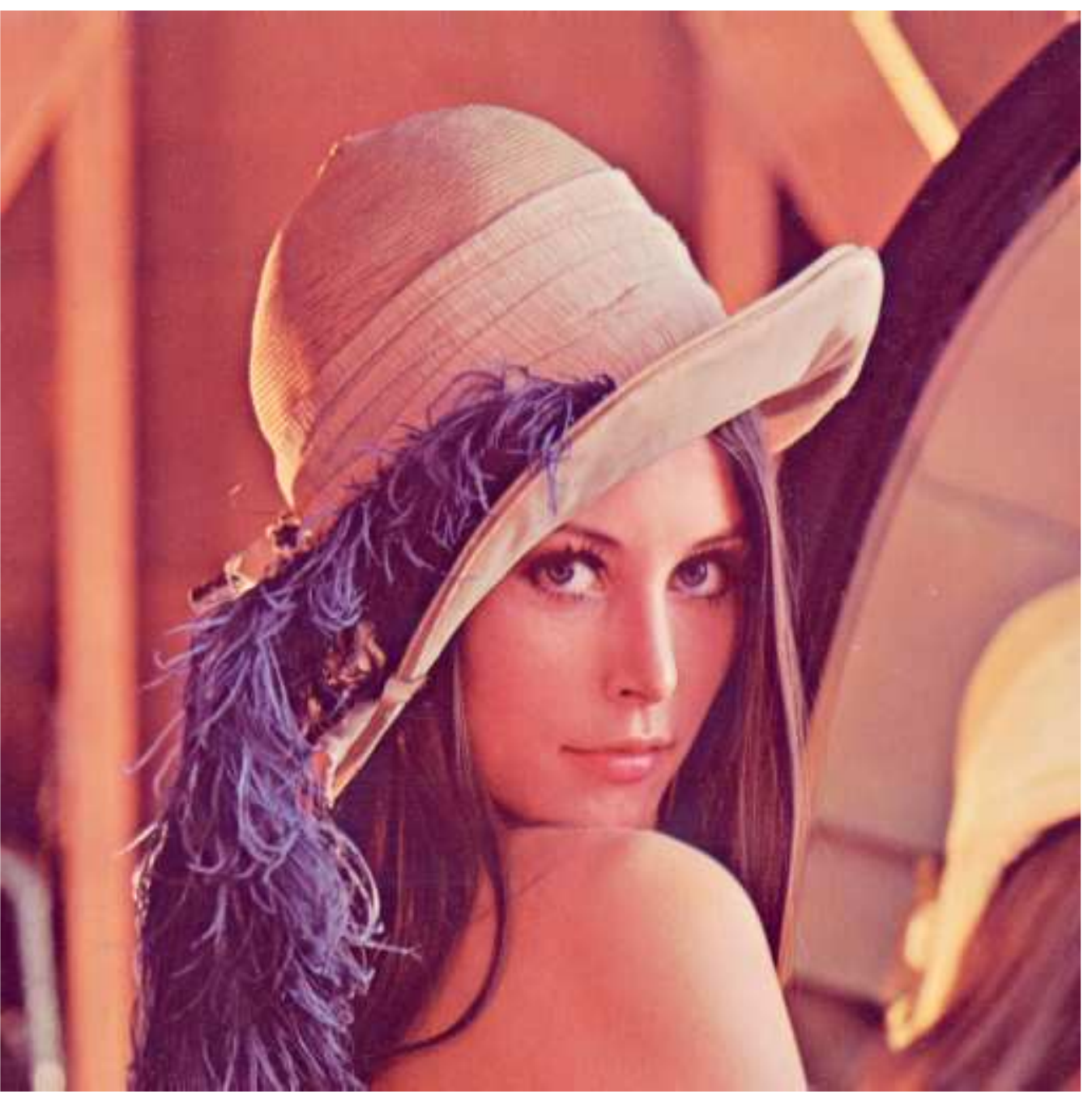}}}\ghs
\subfloat[\figsizetwo corrupted `lenna']{\fcolorbox{colorone}{colortwo}{\includegraphics[width=\figurewidth, height=\figureheight]{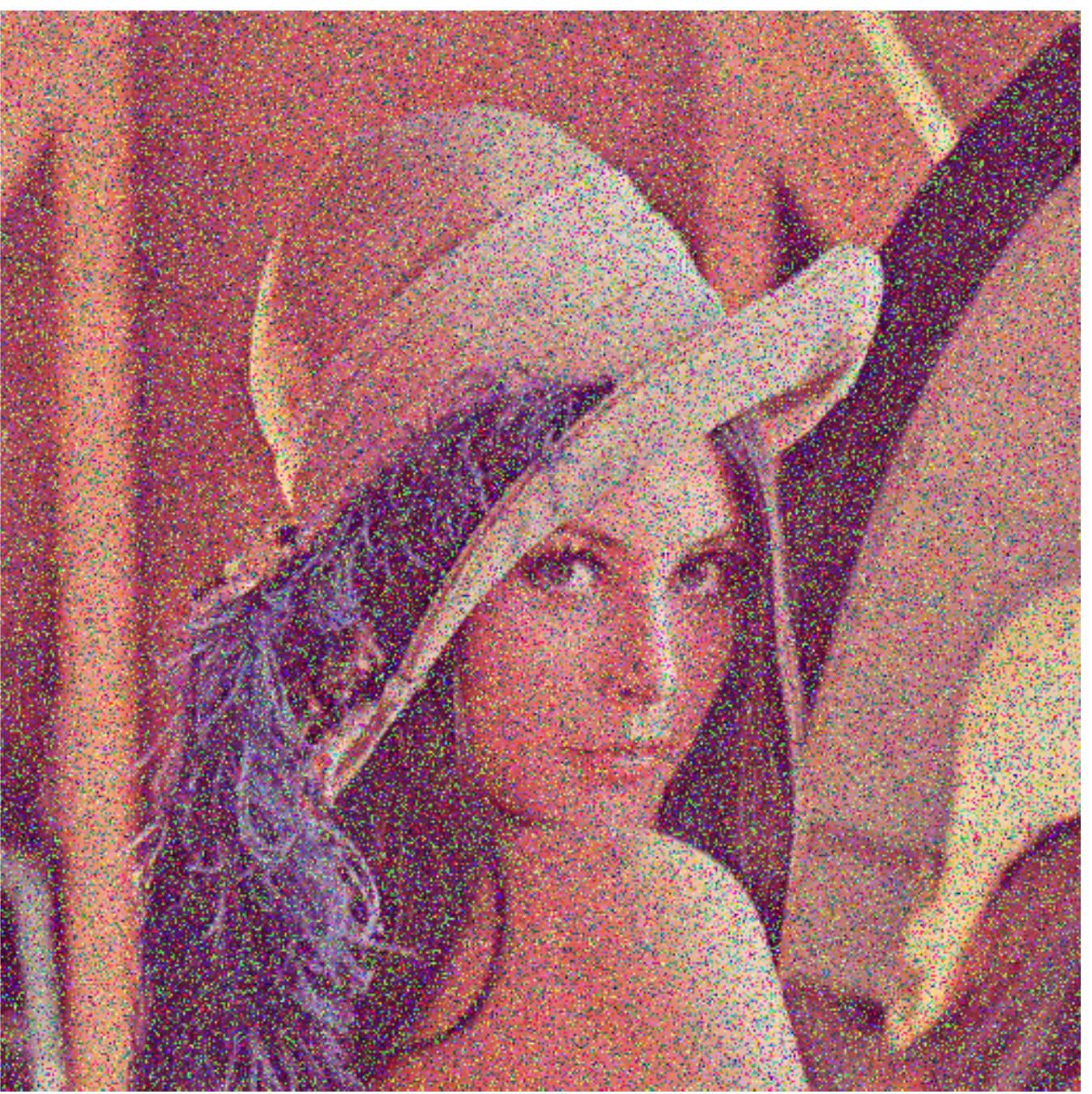}}}\ghs
\subfloat[\figsizetwo recovered `lenna']{\fcolorbox{colorone}{colortwo}{\includegraphics[width=\figurewidth, height=\figureheight]{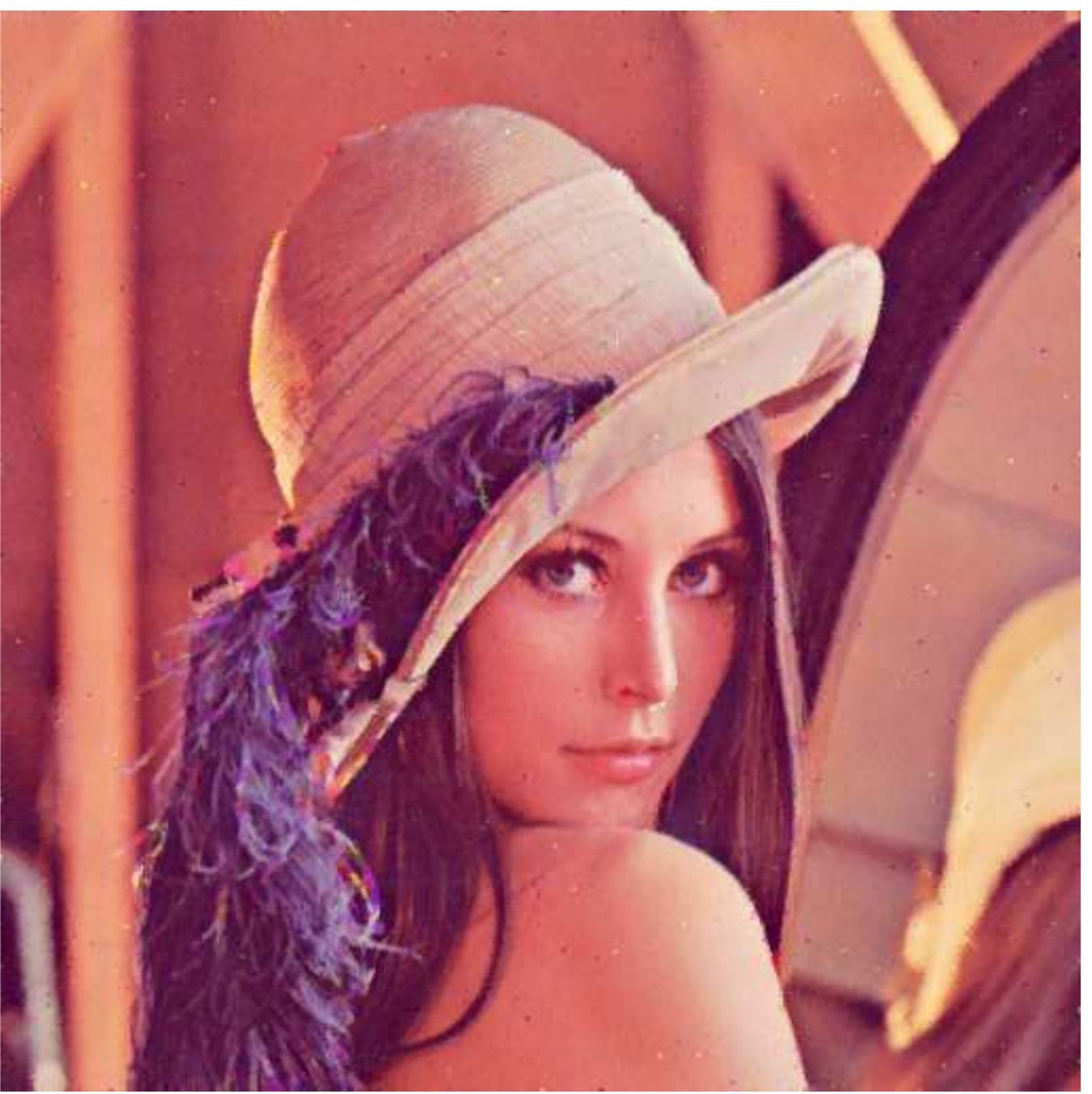}}}

\caption{Colored image denoising problems.}
\label{fig:color:image}
\end{center}
\end{figure}

Our proposed method can be directly extended to its color version. Since color total variation is not the main theme of this paper, we only provide a basic implementation of it. Specifically, we compute the color total variation channel-by-channel, and take a $\ell_1$-norm of the resulting vectors. Suppose we have RGB channels, then we have the following optimization problem:
\beq \label{eq:color:optimization}
 \min_{\bbb{0} \leq \bbb{u}^1,\bbb{u}^2,\bbb{u}^3 \leq \bbb{1}}~{\textstyle \sum_{k=1}^3 (\|\bbb{o}^k\odot (\bbb{Ku}^{k}-\bbb{b}^{k})\|_0 +   \lambda    \|\bbb{\nabla} \bbb{u}^k \|_{p,1})},\nn
\eeq
\noi where $\bbb{o}^k$ and $\bbb{u}^k$ are the prior and the solution of the $k${th} channel. The grayscale proximal ADM algorithm in Algorithm 1 can be directly extended to solve the optimization above. We demonstrate its applicability in colored image denoising problems in Figure \ref{fig:color:image}. The regularization parameter $\lambda$ is set to $8$ for the three images in our experiments.

\subsection{Running Time Comparisons}
We provide some running time comparisons for the methods $\ell_1TV$-SBM, TSM, $\ell_pTV$-ADMM, $\ell_{02}TV$-AOP, $\ell_0TV$-PDA, and $\ell_0TV$-PADMM on grayscale image `cameraman' corrupted by 50\% random-value impulse noise. For RGB color images, the running time is three times the amount of grayscale images since the colored image recovery problem can be decomposed into dependent subproblems. Table \ref{tab:time:comparision} shows the average CPU time for five runs. Generally, our method is efficient and comparable with existing solutions. This is expected since our method is an alternating optimization algorithm.

\begin{table}[!h]
\caption{CPU time (in seconds) comparisons. First row: image denoising; second row: image deblurring.}
\label{tab:time:comparision}
\begin{center}
\begin{tabular}{p{0.8cm}p{0.8cm}p{0.8cm}p{0.8cm}p{0.8cm}p{0.8cm}}
\hline
$\ell_1TV$-$SBM$ & $TSM$ & $\ell_pTV$-ADMM  & $\ell_{02}TV$-AOP & $\ell_0TV$-PDA & $\ell_0TV$-$PADMM$ \\
\hline
5$\pm$4 & 6$\pm$4 & 15$\pm$4 &30$\pm$5 & 17$\pm$3 & 14$\pm$4 \\
15$\pm 8$ & 16$\pm 7$ & 38$\pm$8 &62$\pm$4 & 39$\pm$7 & 35$\pm$8 \\
\hline
\end{tabular}
\end{center}
\end{table}

\section{Conclusions} \label{sec:conc}
In this paper, we propose a new method for image restoration based on total variation (TV) with $\ell_0$-norm data fidelity, which is particularly suitable for removing impulse noise. Although the resulting optimization model is non-convex, we design an efficient and effective proximal ADM method for solving the equivalent MPEC problem of the original $\ell_0$-norm minimization problem. Extensive numerical experiments indicate that the proposed $\ell_0$TV model significantly outperforms the state-of-the-art in the presence of impulse noise. In particular, our proposed proximal ADM solver is more effective than the penalty decomposition algorithm used for solving the $\ell_0$TV problem \cite{LuZ13}.



\vspace{10pt}
\noi \textbf{Acknowledgments.} We would like to thank Prof. Shaohua Pan for her helpful discussions on this paper. We also thank Prof. Ming Yan for sharing his code with us. This work was supported by the King Abdullah University of Science and Technology (KAUST) Office of Sponsored Research and, in part, by the NSF-China (61772570, 61402182).

\bibliographystyle{ieee}
\bibliography{my}
\begin{IEEEbiography}[{\includegraphics[width=1in,height=1.25in]{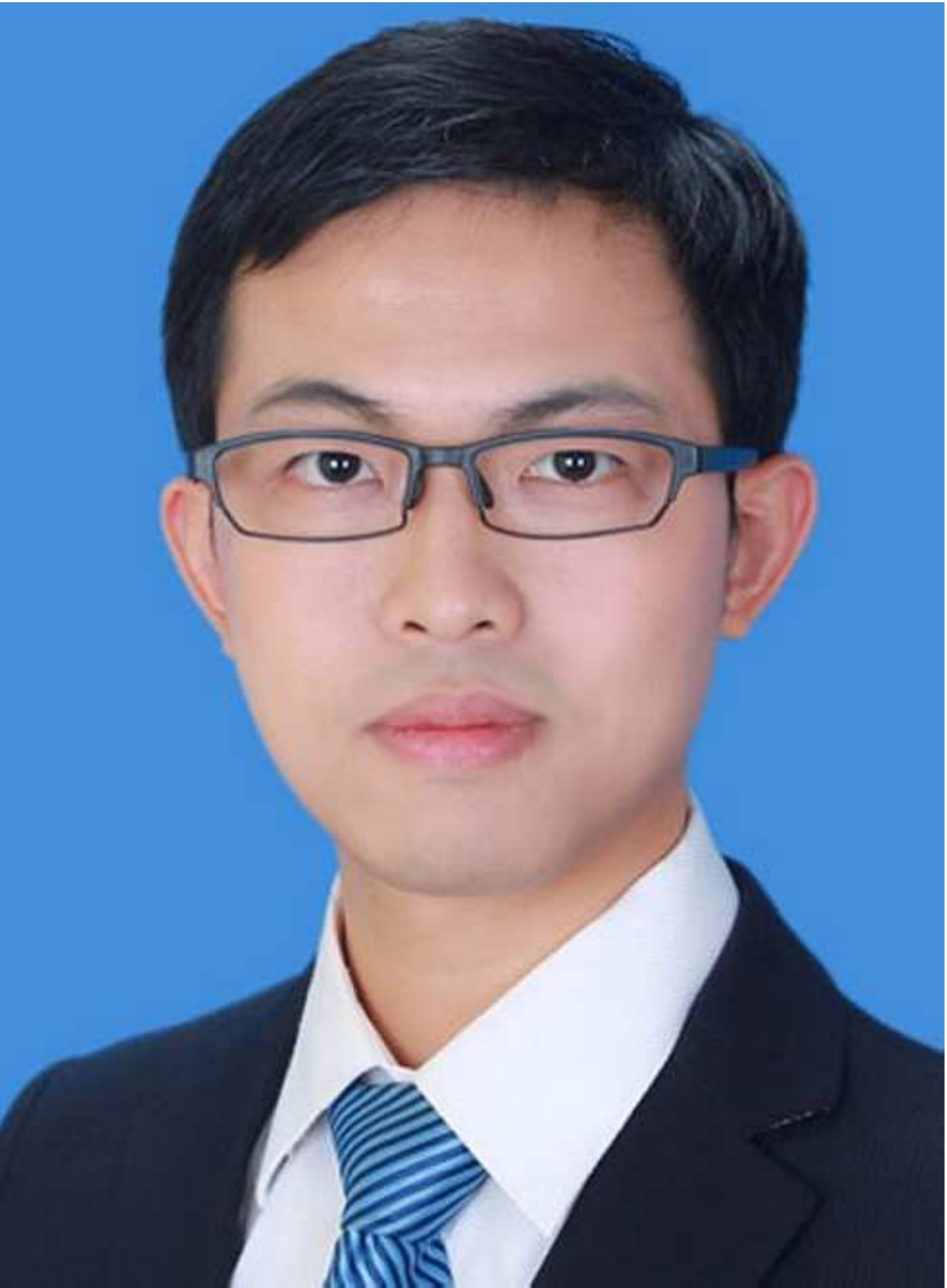}}]{Ganzhao Yuan} was born in Guangdong, China. He received his Ph.D. in School of Computer Science and Engineering, South China University of Technology (SCUT) in 2013. He is currently a research associate professor at School of Data and Computer Science in Sun Yat-sen University (SYSU). His research interests primarily center around large-scale nonlinear optimization and its applications in computer vision and machine learning. He has published papers in ICML, SIGKDD, AAAI, CVPR, VLDB, and ACM Transactions on Database System (TODS).
\end{IEEEbiography}

\begin{IEEEbiography}[{\includegraphics[width=1in,height=1.25in]{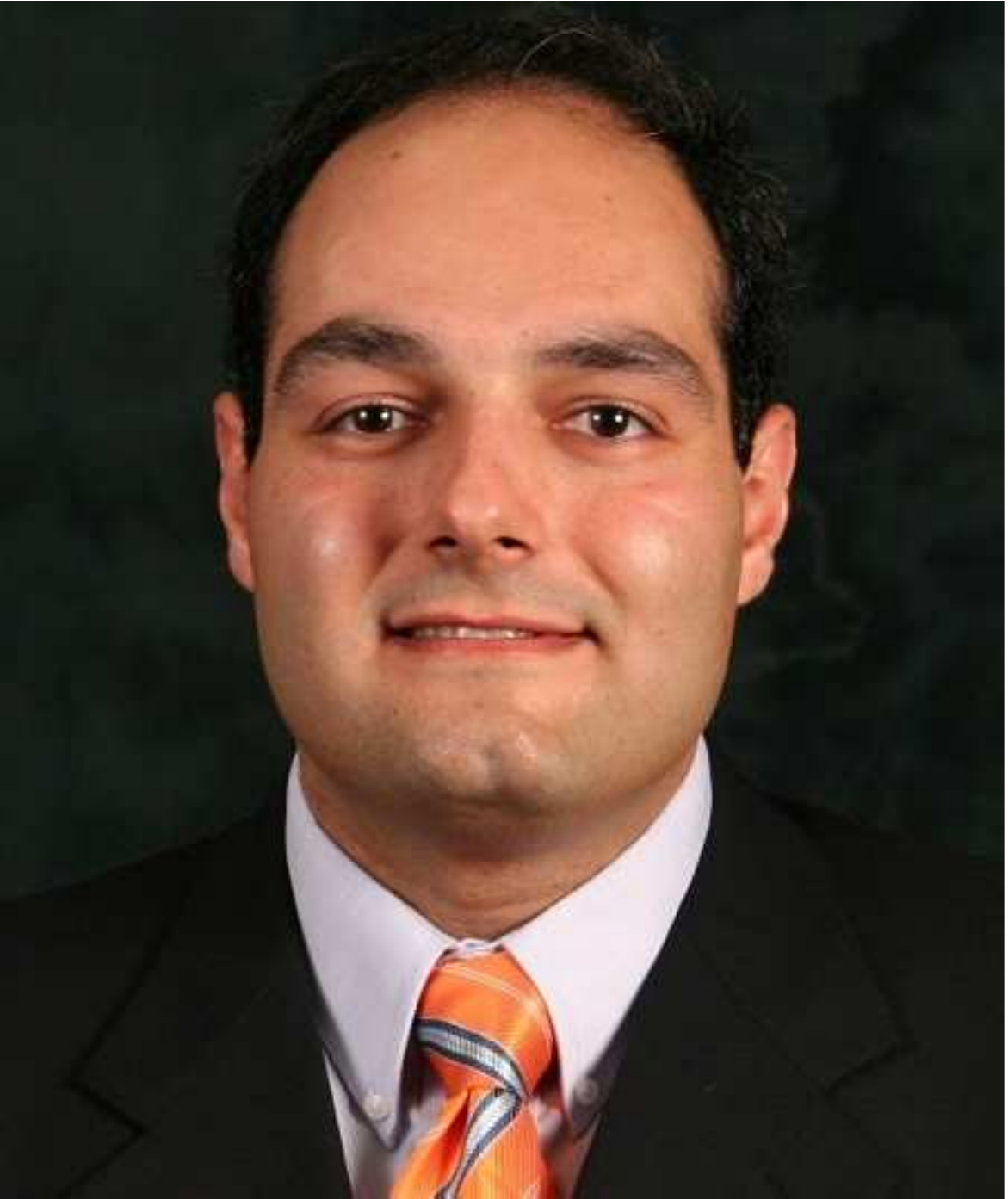}}]{Bernard Ghanem} was born in Betroumine, Lebanon. He received his Ph.D. in Electrical and Computer Engineering from the University of Illinois at Urbana-Champaign (UIUC) in 2010. He is currently an assistant professor at King Abdullah University of Science and Technology (KAUST), where he leads the Image and Video Understanding Lab (IVUL). His research interests focus on designing, implementing, and analyzing approaches to address computer vision problems (e.g. object tracking and action recognition/detection in video), especially at large-scale.
\end{IEEEbiography}

\appendices
\section{Proof of Theorem 1} \label{app:the1}

\begin{proof}
We define $Z\triangleq(X,Y)$ and denote $I(\cdot)$ as the indicator function on the constrained set $\Delta \triangleq \{z~|~\bbb{0}\leq \bbb{z}\leq \bbb{1}\}$. First of all, we present the first-order KKT conditions of the MPEC reformulation. Based on the augmented Lagrangian function $\L$, we naturally derive the following KKT conditions for $\{\bbb{u}^*, \bbb{v}^*, \bbb{x}^*, \bbb{y}^*, \bbb{\xi}^*, \bbb{\zeta}^*, \bbb{\pi}^*\}$:
\beq \label{eq:kkt:conditions}
 0 &\in& \bbb{\nabla}^T \bbb{\xi}^* + \bbb{K}^T\bbb{\zeta}^* + \partial I(\bbb{u}^*) \nn\\
 0 &\in & \bbb{\pi}^* \odot \bbb{o} \odot |\bbb{y}^*| -\bbb{1} + \partial I(\bbb{v}^*)\nn\\
 0 &\in& \partial \lambda\|\bbb{x}^*\|_{p,1} - \bbb{\xi}^* \nn\\
 0 &\in& \bbb{\pi}^* \odot \bbb{v}^*\odot \bbb{o}\odot \partial \|\bbb{y}^*\|_1-\bbb{\zeta}^*  \\
 0 &=& \bbb{\nabla} \bbb{u}^*-\bbb{x}^*~~~~~~~~~~~~~~~~~~~~~~~~~~~~~~~~~~ \nn\\
 0 &=&\bbb{Ku}^*-\bbb{b}-\bbb{y}^*~~~~~~~~~~~~~~~~~~~~~~~~~\nn\\
 0 &=&\bbb{o}\odot \bbb{v}^*\odot  |\bbb{y}^*|.~~~~~~~~~~~~~~~~~~~~~\nn
\eeq 

Secondly, we prove that the solution is convergent: $Z^{k+1}-Z^{k}\rightarrow0$. We observe that $\mathcal{L}$ can be rewritten as:
\beq
\mathcal{L}(Z) \triangleq \la \bbb{1}, \bbb{1}-\bbb{v}\ra + \lambda\|\bbb{x}\|_{p,1}+\tfrac{\beta}{2} \|\bbb{\nabla} \bbb{u} - \bbb{x} + \bbb{\xi}/\beta\|^2 \nn~~\\
- \tfrac{1}{2\beta}\|\bbb{\xi}\|^2 + \frac{\beta}{2} \|\bbb{Ku}-\bbb{b}-\bbb{y}+ \bbb{\zeta}/\beta\|^2 -  \tfrac{1}{2\beta}\|\bbb{\zeta}\|^2 ~\nn\\
+\tfrac{\beta}{2} \|\bbb{v}\odot \bbb{o}\odot |\bbb{y}| + \bbb{\pi}/\beta\|^2  -  \tfrac{1}{2\beta}\|\bbb{\pi}\|^2.~~~~~~~~~~~~~~~~\nn
\eeq
\noi Since $Y\triangleq(\bbb{\xi,\zeta,\pi})$ is bounded by assumption, $\mathcal{L}(Z)$ is bounded below for all $Z$. We now define $\mathcal{J}(Z)$ as:
\beq
\mathcal{J}(Z) = \mathcal{L}(Z) + \tfrac{1}{2}\|\bbb{u}-\bbb{u}'\|_{\bbb{D}}^2  + \tfrac{1}{2}\|\bbb{v}-\bbb{v}'\|_{\bbb{E}}^2,\nn
\eeq
\noi where $\bbb{u}'$ and $\bbb{v}'$ denote the values of $\bbb{u}$ and $\bbb{v}$ in the previous iteration. We define $Z^{-1}=Z^{0}$, and the variable $Z$ in $\mathcal{J}(Z)$ is in the range of $\{Z^0,~Z^1,~Z^2,...\}$. Since $\mathcal{J}(Z)$ is strongly and jointly convex with respect to $\{\bbb{u},\bbb{v}\}$ and $\{\bbb{u}^{k+1},~\bbb{v}^{k+1}\}$ is the minimizer of $\min_{\bbb{u},\bbb{v}}~\J(\bbb{u},~\bbb{v},~\bbb{x}^k,~\bbb{y}^k,~Y^k)$ which is based on $\{\bbb{u}^{k},~\bbb{v}^{k}\}$, using the second order growth condition, we have:
\beq \label{eq:J:u}
\begin{split}
\J(\bbb{u}^k,\bbb{v}^k,\bbb{x}^k,\bbb{y}^k,Y^k)-\J(\bbb{u}^{k+1},\bbb{v}^{k+1},\bbb{x}^k,\bbb{y}^k,Y^k) \\
\geq \tfrac{\mu}{2} \| \bbb{u}^k - \bbb{u}^{k+1}\|^2 + \tfrac{\mu}{2} \| \bbb{v}^k - \bbb{v}^{k+1}\|^2.~~~~~~~~~~~~
\end{split}
\eeq
\noi Using the same methodology for the variable $\bbb{x}$ and $\bbb{y}$, we have the following inequalities:
\beq \label{eq:J:vxy}
\begin{split}\J(\bbb{u}^{k+1},\bbb{v}^{k+1},\bbb{x}^k,\bbb{y}^k,Y^k)\\
- \J(\bbb{u}^{k+1},\bbb{v}^{k+1},\bbb{x}^{k+1},\bbb{y}^{k+1},Y^k) \\
\geq \tfrac{\beta}{2} \| \bbb{x}^k - \bbb{x}^{k+1}\|^2 + \tfrac{\beta}{2} \| \bbb{y}^k - \bbb{y}^{k+1}\|^2.
\end{split}\eeq
\noi Denoting $\rho=\tfrac{1}{2}\min(\mu,\beta)$ and combining (\ref{eq:J:u}) and (\ref{eq:J:vxy}), we obtain:
\beq  \label{eq:J:dec}
\J(X^k,Y^k) - \J(X^{k+1},Y^k) \geq  \rho \|X^k-X^{k+1}\|_F^2.
\eeq
\noi
Using the definition of $\J$ and the update rule of the multipliers, we have:
\beq \label{eq:dual:seq}
&&\J(X^{k+1},Y^{k+1}) - \J(X^{k+1},Y^{k}) \nn\\
&=&  \la \bbb{\nabla} \bbb{u}^{k+1} - \bbb{x}^{k+1}, \bbb{\xi}^{k+1} - \bbb{\xi}^k \ra  + \nn\\
 && \la \bbb{Ku}^{k+1}-\bbb{b}-\bbb{y}^{k+1} , \bbb{\zeta}^{k+1} - \bbb{\zeta}^{k}\ra  +\nn\\
  && \la \bbb{v}^{k+1} \odot \bbb{o}\odot |\bbb{y}^{k+1}|  , \bbb{\pi}^{k+1} - \bbb{\pi}^{k} \ra \nn\\
&=& \tfrac{1}{\gamma \beta} \| Y^{k+1} - Y^k\|^2.
\eeq

\noi Combining (\ref{eq:J:dec}) and (\ref{eq:dual:seq}), we have:
\beq
&&\J(X^k,Y^k) - \J(X^{k+1},Y^{k+1}) \nn\\
&\geq& \rho  \|X^k-X^{k+1}\|_F^2 - \tfrac{1}{\gamma \beta} \|Y^k-Y^{k+1}\|_F^2.\nn
\eeq
\noi Taking summation of the above inequality and using the boundedness of $\J(Z)$, we have that:
\beq
~~~~~~\textstyle \sum_{k=0}^{\infty} ( \rho \|X^k-X^{k+1}\|_F^2 - \tfrac{1}{\gamma \beta} \|Y^k-Y^{k+1}\|_F^2) \nn\\
\leq \J(X^0,Y^0) - \J(X^{\infty},Y^{\infty}) < \infty.~~~~~~~~~~~~~~~~~~~~~~\nn
\eeq
Since the second term in the inequality above is bounded, i.e. $\sum_{k=0}^{\infty} \lim_{k\rightarrow \infty} \|Y^k-Y^{k+1}\|_F^2 = 0$, we obtain that $\sum_{k=0}^{\infty} \lim_{k\rightarrow \infty} \|X^k-X^{k+1}\|_F^2 = 0$ and $X^k-X^{k+1}\rightarrow0$.

Finally, we are ready to prove the result of the theorem. By the update rule of $Y^k$, we have:
\beq
&&\bbb{\xi}^{k+1}-\bbb{\xi}^{k} = \gamma \beta(\bbb{\nabla} \bbb{u}^{k} - \bbb{x}^{k}) \nn \\
&&\bbb{\zeta}^{k+1} - \bbb{\zeta}^{k} = \gamma \beta(\bbb{K} \bbb{u}^{k} -\bbb{b}- \bbb{y}^{k}) \nn \\
&&\bbb{\pi}^{k+1}-\bbb{\pi}^{k} = \gamma \beta(\bbb{o}\odot \bbb{v}^{k}\odot |\bbb{y}^{k}|).\nn
\eeq
\noi Using the convergence of $Y$ that $Y^{k}-Y^{k+1}\rightarrow 0$ and the optimality of $X^{k+1}$ with respect to $\J(\cdot)$, we have:
\beq
0 &=& \bbb{\nabla}^T \bbb{\xi}^{k} + \bbb{K}^T\bbb{\zeta}^{k} + \partial I(\bbb{u}^{k+1}) + \mu( \bbb{u}^{k+1} - \bbb{u}^k) \nn\\
0 &=& \bbb{\pi}^{k} \odot \bbb{o} \odot |\bbb{y}^{k}| -\bbb{1} + \partial I(\bbb{v}^{k+1}) + \mu(\bbb{v}^{k+1} - \bbb{v}^k) \nn\\
0 &\in& \partial \lambda\|\bbb{x}^{k+1}\|_{p,1} - \bbb{\xi}^{k} \nn\\
0 &\in& \bbb{\pi}^{k} \odot \bbb{v}^{k+1}\odot \bbb{o}\odot \partial \|\bbb{y}^{k+1}\|_1-\bbb{\zeta}^{k}.\nn
\eeq
\noi Combining the convergence of $X$ that: $X^{k}-X^{k+1}\rightarrow 0$, we have
\beq
 0 &\in& \bbb{\nabla}^T \bbb{\xi}^{k+1} + \bbb{K}^T\bbb{\zeta}^{k+1} + \partial I(\bbb{u}^{k+1}) \nn\\
 0 &\in & \bbb{\pi}^{k+1} \odot \bbb{o} \odot |\bbb{y}^{k+1}| -\bbb{1} + \partial I(\bbb{v}^{k+1})\nn\\
 0 &\in& \partial \lambda\|\bbb{x}^{k+1}\|_{p,1} - \bbb{\xi}^{k+1} \nn\\
 0 &\in& \bbb{\pi}^{k+1} \odot \bbb{v}^{k+1}\odot \bbb{o}\odot \partial \|\bbb{y}^{k+1}\|_1-\bbb{\zeta}^{k+1}\nn  \\
 0 &=& \bbb{\nabla} \bbb{u}^{k+1}-\bbb{x}^{k+1}~~~~~~~~~~~~~~~~~~~~~~~~~~~~~~~~~~ \nn\\
 0 &=&\bbb{Ku}^{k+1}-\bbb{b}-\bbb{y}^{k+1}~~~~~~~~~~~~~~~~~~~~~~~~~\nn\\
 0 &=&\bbb{o}\odot \bbb{v}^{k+1}\odot  |\bbb{y}^{k+1}|,~~~~~~~~~~~~~~~~~~~~~\nn
\eeq
\noi which coincides with the KKT condition in (\ref{eq:kkt:conditions}). Therefore, $Z^{k+1}$ asymptotically converges to the KKT point.

\end{proof}


\end{document}